\documentclass[onecolumn,final]{elsarticle}

\makeatletter
\def\ps@pprintTitle{%
  \let\@oddhead\@empty
  \let\@evenhead\@empty
  \def\@oddfoot{\reset@font\hfil\thepage\hfil}
  \let\@evenfoot\@oddfoot
}
\makeatother

\usepackage[a4paper]{geometry}

\usepackage{algorithm} 
\usepackage{algpseudocode}
\usepackage{amsmath}
\usepackage{bbm}
\usepackage{bm}
\usepackage{booktabs}
\usepackage{braket}
\usepackage{cases}
\usepackage{chapterbib} 
\usepackage{chngpage}
\usepackage{color}
\usepackage{comment}
\usepackage[bottom]{footmisc} 
\usepackage{graphics}
\usepackage{graphicx}
\usepackage{lettrine}
\usepackage{makeidx}          
\usepackage{multicol}
\usepackage{multirow}
\usepackage{newtxmath} 
\usepackage{newtxtext}        
\usepackage{overpic}
\usepackage{pgfplots}
\usepackage{pgfplotstable}
\pgfplotsset{compat=newest}
\usepackage{placeins}
\usepackage{rotating}
\usepackage{setspace}
\usepackage{sidecap}
\usepackage{siunitx}
\usepackage{subfigure}
\usepackage{tabularx}
\usepackage{times}
\usepackage{tikz}
\usepackage{type1cm}

\usepackage{float}

\definecolor{MyDarkGreen}{rgb}{0,0.45,0}

\newcommand{\tommy}[1]{\noindent{\color{blue}\textbf{[T:~#1]}}}

\newtheorem{theorem}{Theorem}[section]
\newtheorem{lemma}[theorem]{Lemma}
\newtheorem{corollary}[theorem]{Corollary}

\newtheorem{definition}[theorem]{Definition}

\newtheorem{theo}{Theorem}[section]{\bfseries}{\itshape}
\newtheorem{gass}[theorem]{Assumption G\!}{\bfseries}{\itshape}
{\bfseries}{\itshape}

\newcommand{\Pinabla}{\Pi^{\nabla}}

\newcommand{\Vh}{V_h}

\newcommand{\dofiglob}[1]{\mathrm{dof}_i(#1)}

\newcommand{\perf}[1]{\mathbf{p}_{#1}}

\newcommand{\IC}{IC}
\newcommand{\CC}{CC}
\newcommand{\CR}{CR}
\newcommand{\AR}{AR}
\newcommand{\KE}{KE}
\newcommand{\KA}{KAR}
\newcommand{\AP}{APR}
\newcommand{\MA}{MA}
\newcommand{\SE}{SE}
\newcommand{\ER}{ER}
\newcommand{\MP}{MPD}
\newcommand{\MX}{MXA}

\newcommand{\SR}{SR}

\newcommand{\hmax}{h_{\text{max}}}
\newcommand{\hmin}{h_{\text{min}}}
\newcommand{\hav}{h_{\text{av}}}

\def\trait #1 #2 #3 {\vrule width #1pt height #2pt depth #3pt}
\def\fin{\hfill
        \trait .3 5 0
        \trait 5 .3 0
        \kern-5pt
        \trait 5 5 -4.7
        \trait 0.3 5 0
\medskip}

\newcommand{\KER} {\textrm{ker}}
\newcommand{\INTP}{I}

\newcommand{\REAL}{\mathbbm{R}}

\newcommand{\PGRAPH}[1]{\medskip\noindent\textbf{#1}}

\newcommand{\kv}{\mathbf{k}}

\newcommand{\xv}{\mathbf{x}}
\newcommand{\yv}{\mathbf{y}}

\newcommand{\calB}{\mathcal{B}}

\newcommand{\calD}{\mathcal{D}}
\newcommand{\calE}{\mathcal{E}}

\newcommand{\calH}{\mathcal{H}}
\newcommand{\calI}{\mathcal{I}}

\newcommand{\calM}{\mathcal{M}}
\newcommand{\calN}{\mathcal{N}}

\newcommand{\calT}{\mathcal{T}}

\newcommand{\calV}{\mathcal{V}}

\newcommand{\as}{a}

\newcommand{\cs}{c}

\newcommand{\fs}{f}

\newcommand{\ks}{k}

\newcommand{\qs}{q}

\renewcommand{\ss}{s}

\newcommand{\us}{u}
\newcommand{\vs}{v}
\newcommand{\ws}{w}
\newcommand{\xs}{x}
\newcommand{\ys}{y}

\newcommand{\Cs}{C}

\newcommand{\Fs}{F}

\newcommand{\Ns}{N}

\newcommand{\Ss}{S}

\newcommand{\xsP}{\xs_{\P}}
\newcommand{\ysP}{\ys_{\P}}
\newcommand{\xsE}{\xs_{\E}}
\newcommand{\ysE}{\ys_{\E}}

\newcommand{\PS}[1]{\mathbbm{P}_{#1}}

\newcommand{\HONE}   {H^1}
\newcommand{\HONEzr} {H^1_0}

\newcommand{\LTWO}  {L^2}

\newcommand{\HS}[1] {H^{#1}}

\newcommand{\CS}[1] {C^{#1}}

\renewcommand{\P} {E}
\newcommand{\E} {e}

\newcommand{\hh}{h}
\newcommand{\Th}{\Omega_{\hh}}

\newcommand{\xvP}{\xv_{\P}}        
\newcommand{\xvE}{\xv_{\E}}        

\newcommand{\hP}{\hh_{\P}}

\newcommand{\hE}{\hh_{\E}}

\newcommand{\mP}{\ABS{\P}}

\newcommand{\mE}{\ABS{\E}}

\newcommand{\Sset}{\mathcal{S}}    

\newcommand{\NMB}{N}
\newcommand{\NS}{\NMB^{\Sset}}   

\newcommand{\dV}{\,dV}
\newcommand{\dS}{\,dS}
\newcommand{\dxv}{\,d\xv}

\newcommand{\fsh}{\fs_{\hh}}

\newcommand{\ush}{\us_{\hh}}

\newcommand{\usI}{\us_{\INTP}}

\newcommand{\vsh}{\vs_{\hh}}

\newcommand{\wsh}{\ws_{\hh}}

\newcommand{\wsI}{\ws_{\INTP}}

\newcommand{\Fsh}{\Fs_{\hh}}

\newcommand{\asP}{\as^{\P}}

\newcommand{\ash}{\as_{\hh}}

\newcommand{\ashP}{\as^{\P}_{\hh}}

\newcommand{\nlen}{\hspace{-0.2mm}}

\newcommand{\snorm}  [2]{|#1|_{#2}}

\newcommand{\norm}   [2]{|\nlen|#1|\nlen|_{#2}}

\newcommand{\Tnorm}  [2]{|\nlen|\nlen|#1|\nlen|\nlen|_{{#2}}}

\newcommand{\ABS}    [1]{\left|#1\right|}

\newcommand{\Vhk} {V^{\hh}_{k}}

\newcommand{\Pin}[1]{\Pi^{\nabla}_{#1}}
\newcommand{\Piz}[1]{\Pi^{0}_{#1}}
\newcommand{\PinP}[1]{\Pi^{\nabla,\P}_{#1}}
\newcommand{\PizP}[1]{\Pi^{0,\P}_{#1}}

\newcommand{\SPh} {S^{\P}_{\hh}}
\newcommand{\SPhdd} {S^{\P,\textrm{dd}}_{\hh}}
\newcommand{\SPhtr} {S^{\P,\textrm{tr}}_{\hh}}

\newcommand{\cbot}{c_*}
\newcommand{\ctop}{c^*}

\newcommand{\restrict}[2]{{#1}_{|{#2}}}

\newcommand{\EOD}{\end{document}}

\newcommand{\roundPrecision}{2}
\newcommand{\bP} {\partial \P}            

\newcommand{\NDOFS}{N^{\textrm{dofs}}}

\newcommand{\matPin}[1]{\boldsymbol{\Pi}^{\nabla}_{#1}}
\newcommand{\matPiz}[1]{\boldsymbol{\Pi}^{0}_{#1}}

\newcommand{\Dtriangle}{\calD_{\rm{Triangle}}}
\newcommand{\Dmaze}{\calD_{\rm{Maze}}}
\newcommand{\Dstar}{\calD_{\rm{Star}}}
\newcommand{\Djenga}{\calD_{\rm{Jenga}}}
\newcommand{\Dslices}{\calD_{\rm{Slices}}}
\newcommand{\Dulike}{\calD_{\rm{Ulike}}}
\newcommand{\DjengaM}{\calD_{\rm{Jenga4}}}
\newcommand{\DslicesM}{\calD_{\rm{Slices4}}}
\newcommand{\DulikeM}{\calD_{\rm{Ulike4}}}

\newcommand{\rP}{r_\P}

\newcommand{\vettore}[1]{\mathbf{#1}}

\newcommand{\stiffness}{\matrice{R}}

\newcommand{\matrice}[1]{\mathbf{#1}}

\begin{document}

\begin{frontmatter}

  \title{VEM and the Mesh}

  \author[IMATIGE]{T.~Sorgente}
  \author[IMATIPV]{D.~Prada}
  \author[IMATIGE]{D.~Cabiddu}
  \author[IMATIGE]{S.~Biasotti}
  \author[IMATIGE]{G.~Patan\`e\\}
  \author[IMATIPV]{M.~Pennacchio}
  \author[IMATIPV]{S.~Bertoluzza}
  \author[IMATIPV]{G.~Manzini}
  \author[IMATIGE]{M.~Spagnuolo}

  \address[IMATIGE]{IMATI, Consiglio Nazionale delle
    Ricerche, Genova, Italy }

  \address[IMATIPV]{IMATI, Consiglio Nazionale delle
    Ricerche, Pavia, Italy }

  \begin{abstract}
    In this work we report some results, obtained within the framework
    of the ERC Project CHANGE, on the impact on the performance of the
    virtual element method of the shape of the polygonal elements of
    the underlying mesh.
    More in detail, after reviewing the state of the art, we present
    a) an experimental analysis of the convergence of the VEM under
    condition violating the standard shape regularity assumptions, b)
    an analysis of the correlation between some mesh quality metrics
    and a set of different performance indexes, and c) a suitably
    designed mesh quality indicator, aimed at predicting the quality
    of the performance of the VEM on a given mesh.
  \end{abstract}


\end{frontmatter}





\section{Introduction}\label{sec:introduction}

Geometrically complex domains are frequently encountered in
mathematical models of real engineering applications.
Their representation in some discrete form is a key aspect in the
numerical approximation of the solutions of the partial differential
equations (PDEs) describing such models and can be extremely
difficult.
The finite element method has been proved to be very successful as
it allows the computational domains to be discretized by using
triangular and quadrilateral meshes in 2D and tetrahedral and
hexahedral meshes in 3D.
Meshes with more complex elements, even admitting curved edges and
faces, can be considered in the finite element formulation through
reference elements and suitable remappings onto the problem space.
To obtain accurate solutions, stringent constraints must be imposed, for
example on the internal angles of the triangles and tetrahedra, thus
requiring in extreme situations meshes with very small sized elements.
To alleviate meshing issues, we can resort to numerical methods that
are designed from the very beginning to provide arbitrary order of
accuracy on more generally shaped elements.
A class of methods with these features is the class of the so called 
 polytopal element
method, or PEM for short, that make it possible to numerically 
solve PDEs using polygonal and polyhedral grids.

\medskip

The PEMs
 allow the user to
incorporate complex geometric features at different scales without
triggering mesh refinement, thus achieving high flexibility in the
treatment of complex geometries.
Moreover, nonconformal meshes can be treated in a straightforward way,
by automatically including hanging nodes (i.e., T-junctions), and the
design of refinement and coarsening algorithms is greatly simplified.

Such polytopal methods normally rely on a special design, as a
straightforward generalization of the FEM is not possible because the
finite element variational formulation requires an explicit knowledge
of the basis functions.
This requirement typically implies that such basis functions are  elements of a subset of scalar and vector
polynomials, and, as a consequence, the FEM is mostly restricted to meshes with elements
having a simple geometrical shape, such as triangles or
quadrilaterals.

The virtual element method (VEM) is a very successful example of PEM.
The VEM formulation and implementation are based on suitable
polynomial projections that we can always compute from the degrees of
freedom of the basis functions.
Since the explicit knowledge of the basis functions for the 
 approximation space is not
required, the method is dubbed as \emph{virtual}.
This fundamental property allows the virtual element to be formulated
on meshes with elements having a very general geometric shapes.

The VEM was originally formulated
in~\cite{BeiraodaVeiga-Brezzi-Cangiani-Manzini-Marini-Russo:2013} as a
conforming FEM for the Poisson problem by rewriting in a variational
setting the \emph{nodal} mimetic finite difference (MFD)
method~\cite{Brezzi-Buffa-Lipnikov:2009,%
  BeiraodaVeiga-Lipnikov-Manzini:2011,%
  BeiraodaVeiga-Manzini-Putti:2015,%
  Manzini-Lipnikov-Moulton-Shashkov:2017}
for solving diffusion problems on unstructured polygonal meshes.
A survey on the MFD method can be found in the review
paper~\cite{Lipnikov-Manzini-Shashkov:2014} and the research
monograph~\cite{BeiraodaVeiga-Lipnikov-Manzini:2014}.
The VEM scheme inherits the flexibility of the MFD method with respect to
the admissible meshes and this feature is well reflected in the many
significant applications that have been developed so far, see, for
example,~\cite{%
  BeiraodaVeiga-Manzini:2014,%
  Benedetto-Berrone-Pieraccini-Scialo:2014,%
  BeiraodaVeiga-Manzini:2015, Berrone-Pieraccini-Scialo-Vicini:2015,%
  Mora-Rivera-Rodriguez:2015,%
  Paulino-Gain:2015,%
  Antonietti-BeiraodaVeiga-Scacchi-Verani:2016,%
  BeiraodaVeiga-Chernov-Mascotto-Russo:2016,%
  BeiraodaVeiga-Brezzi-Marini-Russo:2016a,%
  BeiraodaVeiga-Brezzi-Marini-Russo:2016b,%
  Cangiani-Georgoulis-Pryer-Sutton:2016,%
  Perugia-Pietra-Russo:2016,%
  Wriggers-Rust-Reddy:2016,
  Certik-Gardini-Manzini-Vacca:2018:ApplMath:journal,
  Dassi-Mascotto:2018,%
  Benvenuti-Chiozzi-Manzini-Sukumar:2019:CMAME:journal,%
  Antonietti-Manzini-Verani:2019:CAMWA:journal,
  Certik-Gardini-Manzini-Mascotto-Vacca:2020%
}.
Since the VEM is a reformulation of the MFD method that generalizes
the FEM to polytopal meshes as the other PEMs, it is clearly related
with many other polytopal schemes.
The connection between the VEM and finite elements on
polygonal/polyhedral meshes is thoroughly investigated
in~\cite{Manzini-Russo-Sukumar:2014,
  Cangiani-Manzini-Russo-Sukumar:2015, DiPietro-Droniou-Manzini:2018},
between VEM and discontinuous skeletal gradient discretizations
in~\cite{DiPietro-Droniou-Manzini:2018}, and between the VEM and the
BEM-based FEM method
in~\cite{Cangiani-Gyrya-Manzini-Sutton:2017:GBC:chbook}.
The VEM has been extended to convection-reaction-diffusion problems
with variable coefficients
in~\cite{BeiraodaVeiga-Brezzi-Marini-Russo:2016b}. The issue of preconditioning the VEM has been considered in \cite{Antonietti-Mascotto-Verani:2018,BertoluzzaPennacchioPrada:2017,BertoluzzaPennacchioPrada:2020,Calvo:2019}.

If on one hand, the VEM makes it possible to discretize a PDE
on computational domains partitioned by a polytopal mesh, on the other
hand this flexibility poses the fundamental question of what is a
``good'' polytopal mesh.
All available theoretical and numerical results in the literature
strongly support the fact that the quality of the mesh is crucial to
determine the accuracy of the method and its effectiveness in solving
problems on difficult geometries.
This fact is not surprising as such dependence of the behavior of the
numerical approximation in terms of accuracy, stability, and overall
computational cost on the quality of the underlying mesh has been a
very well known fact for decades in the finite element framework and
has been formalized in concepts like mesh regularity, shape
regularity, etc.
However, the concept of shape regularity of triangular/tetrahedral and
quadrangular/hexahedral meshes is well understood
\cite{ciarlet2002finite,Ja-2008,Shewchuk02whatis}, but the
characterization of a good polytopal mesh is still subject to ongoing
research.
Optimal convergence rates for the virtual element approximations of
the Poisson equation were proved in $\HONE$ and $\LTWO$ norms, see for
instance
\cite{BeiraodaVeiga-Brezzi-Cangiani-Manzini-Marini-Russo:2013,
  Ahmad-Alsaedi-Brezzi-Marini-Russo:2013,CHINOSI2016,
  Beirao-Lovadina-Russo:2016, Brenner:2017:SEV, brenner2018virtual,
  da2020sharper}.
These theoretical results involve an estimate of the approximation
error, which is due to both analytical assumptions (interpolation and
polynomial projections of the virtual element functions) and
geometrical assumptions (the geometrical shape of the mesh elements).

A major point here is that the polytopal framework provides an
enormous freedom to the possible geometric shapes of the mesh
elements.
This freedom makes it difficult to identify which geometric features
may have a negative effect on the performance of the VEM.
Various geometrical (or \textit{regularity}) assumptions have been
proposed to ensure that all elements of any mesh of a given mesh
family in the refinement process are sufficiently regular.
Many papers prove the convergence of the VEM and derive optimal error
estimates under the assumption that the polygonal elements are star-shaped with
respect to all the points in a disc whose diameter is comparable with
the diameter of the elements itself.
This assumption is combined with a second scaling assumption on the
mesh elements in the sequence of refined meshes that is used in the
numerical approximation.
For example, we can assume that either the length of all the edges  or
the distance between any two vertices of a polygonal element scale
comparably with the element diameter, or even weaker conditions.
These assumptions guarantee the VEM convergence and optimal estimates
of the approximation error with respect to different norms.
However, as already observed from the very first papers,
cf.~\cite{Ahmad-Alsaedi-Brezzi-Marini-Russo:2013}, the VEM seems to
maintain its optimal convergence rates also when we use mesh families
that do not satisfy the usual geometrical assumptions.
Since the VEM was proposed in~2013, many more examples have accumulate
in the literature (some also reviewed in this chapter) that show that
a virtual element solver on a simple model problem as the Poisson
equation may still provide a very good behavior, even if all the
theoretical conditions of the analysis are violated.
Good behavior means that the VEM is convergent and the loss of
accuracy is significant only when the degeneracy of the meshes becomes
really extreme.
This suggests that more permissive shape-regularity criteria should be
devised so that the VEM can be considered as effective, and a lot of
work still has to be carried out to identify the specific issues that
may negatively affect its accuracy.
Clearly, these points are also crucial to support the design of better
polygonal meshing algorithm for the tessellation of computational
domains.

{
Understanding the influence of the geometrical characteristics of the elements on the performance of the VEM and, more generally, of the PEM, is one of the goals that the unit based at the \emph{Istituto di Matematica Applicata e Tecnologie Informatiche del CNR} is pursuing within the Advanced Grant Project \emph{New CHallenges in (adaptive) PDE solvers: the interplay of ANalysis and GEometry} (CHANGE), whose final goal is the design of tools embracing geometry and analysis within a multi-level, multi-resolution paradigm.
Indeed, a deeper understanding of such an interrelation can provide,  on the one hand, the geometry processing community with information on the requirements to be incorporated into the meshing tools in order to generate ``good'' polytopal meshes, and, on the other hand, the mathematical community with possible new directions to pursue in the theoretical analysis of the VEM. 
}

In Section~\ref{sec:benchmark} we review the main results of a recent work aimed at  identifying
the correlation between the performance of the virtual element method
and a set of polygonal quality metrics.
To this end, a systematic exploration was carried out to correlate the
performance of the VEM and the geometric properties of the polygonal
elements forming the mesh.
In such study, the performance of the VEM is characterized by
different “performance indexes”, including (but not limited to) the
accuracy of the solution and the conditioning of the associated linear
system.
The quality ``metrics'' measuring the ``goodness'' of a mesh are built
by considering several geometric properties of polygons (see subsection \ref{subsec:benchmark:metrics}) from the
simplest ones such as areas, angles, and edge length, to most complex
ones as kernels, inscribed and circumscribed circles. The individual quality metrics of the polygonal elements of a mesh are combined in a single quality metric for the mesh itself by different aggregation strategies, such as
minima, maxima, averages, worst case scenario and Euclidean norm.
The numerical experiments to collect the results are performed on a
family of parametric elements, which is designed to progressively
stress one or more of the proposed geometric metrics, enriched with
random polygons in order to avoid a bias in the study.

A second critical point developed in Section~\ref{sec:indicators}, is the connection between the performance of
the method and how the regularity of the mesh refinements impacts on the 
approximation process.
Note indeed that the way the geometric objects forming a mesh as edges
and polygons (and faces and polytopal elements in 3D) scale is crucial
in all possible geometric assumptions.
It is a remarkable example shown in Section~\ref{subsec:violating:datasets} that even the star-shaped assumption can be violated on a sequence of rectangular 
meshes (rectangular elements are star-shaped!) if the element aspect
ratio scales badly when the mesh is refined.
To study how the geometrical conditions that are found in the
literature may really impact on the convergence and accuracy of the
VEM, we gradually introduce several pathologies in the mesh datasets
used in the numerical experiments.
These datasets systematically violate all the geometrical assumptions,
and enhance a correlation analysis between such assumptions and the
VEM performance.
As expected from other works in the literature, these numerical
experiments confirm\ the remarkable robustness of the VEM as it fails
only in very few and extreme situations and a good convergence rate is
still visible in most examples.
To quantify this correlation, we build an indicator that measure the
violation of the geometrical assumptions.
This indicator depends uniquely on the geometry of the mesh elements.
A correspondence is visible between this indicator and the performance
of the VEM on a given mesh, or mesh family, in terms of approximation
error and convergence rate.
This correspondence and such an indicator can be used to devise a
strategy to evaluate if a given sequence of meshes is suited to the
VEM, and possibly to predict the behaviour of the numerical
discretization \emph{before} applying the method.

The chapter is organized as follows.
In Section~\ref{sec:VEM}, we present the VEM and the convergence
results for the Poisson equation with Dirichlet boundary conditions.
In Section~\ref{subsec:star:assumptions}, we detail the geometrical
assumptions on the mesh elements that are used in the literature to
guarantee the convergence of the VEM.
In Section~\ref{subsec:star:error} we review the major theoretical results on the
error analysis that are available in the virtual element literature,
reporting the geometrical conditions assumed in each result.
In Section~\ref{subsec:violating:datasets}, we present a number of datasets 
which do not satisfy these assumptions, and 
experimentally investigate the convergence of the VEM over them.
In Section~\ref{sec:benchmark} we present the statistical analysis of the correlation between some notable mesh quality metrics and a selection of quantities measuring differents aspects of the performance of the VEM.
In Section~\ref{sec:indicators}, we propose a mesh quality indicator to
predict the behaviour of the VEM over a given dataset.
In Section~\ref{sec:PEMesh-GUI} we present the open source benchmarking software tools \texttt{PEMesh}~\cite{PAPER-CHANGE-GUI-2020}, developed at IMATI.

\section{Model problem}
\label{sec:VEM}

The elliptic model problem that we focus on in this paper is the Poisson equation with Dirichlet boundary conditions. 
In this section, we briefly review the strong and weak forms of the model equations  and recall the formulation of its virtual element discretization.

\medskip
\PGRAPH{The Poisson equation and its virtual element discretization.}
Let $\Omega$ be an open, bounded, connected subset of $\REAL^2$ with polygonal boundary $\Gamma$.
We consider the Poisson equation with homogeneous Dirichlet boundary
conditions, whose strong form is:
\begin{align}
  -\Delta\us = \fs & \qquad\textrm{in~}\Omega,\label{eq:Poisson:strong:A}\\[0.5em]
  \us        = 0   & \qquad\textrm{on~}\Gamma.\label{eq:Poisson:strong:B}
\end{align}
Remark that, while, for the sake of simplicity, we consider here homogeneous
Dirichlet boundary conditions, the method that we are going to present also applies to the non homogeneous case, the extension being straightforward. The variational formulation of
problem~\eqref{eq:Poisson:strong:A}-\eqref{eq:Poisson:strong:B} takes the form: \textit{Find $\us\in\HONEzr(\Omega)$ such that}
\begin{align}
  \as(\us,\vs)=\Fs(\vs)\qquad\forall\vs\in\HONEzr(\Omega),
  \label{eq:Poisson:weak}
\end{align}
with the bilinear form
$\as(\cdot,\cdot):\HONE(\Omega)\times\HONE(\Omega)\to\REAL$  and the right-hand side linear functional $\Fs:\LTWO(\Omega)\to\REAL$  respectively 
defined as
\begin{align}
  \as(\us,\vs)=\int_{\Omega}\nabla\us\cdot\nabla\vs\dxv
  \label{eq:as:def}
\end{align}
and 
\begin{align}
  \Fs(\vs)=\int_{\Omega}\fs\vs\dxv,
  \label{eq:Fs:def}
\end{align}
where we implicitly assumed that $\fs\in\LTWO(\Omega)$.
The well-posedness of the weak formulation~\eqref{eq:Poisson:weak} can be proven by applying the Lax-Milgram theorem \cite[Section 2.7]{Brenner-Scott:2008}, thanks to the coercivity and
continuity of the bilinear form $\as(\cdot,\cdot)$, and to the continuity
the linear functional $\Fs(\vs)$. 

We consider here the  virtual element approximation of equation~\eqref{eq:Poisson:weak},  mainly based on
References~\cite{Ahmad-Alsaedi-Brezzi-Marini-Russo:2013,BeiraodaVeiga-Brezzi-Cangiani-Manzini-Marini-Russo:2013},
which provides an optimal approximation on polygonal meshes when the
diffusion coefficient is variable in space.

\medskip
The discrete equation 
will take the form: \textit{Find $\ush\in\Vhk$ such that}
\begin{align}
  \ash(\ush,\vsh) = \Fsh(\vsh)
  \qquad\vsh\in\Vhk,
  \label{eq:VEM}
\end{align}
where $\ush$, $\Vhk$, $\ash(\cdot,\cdot)$, $\Fsh(\cdot)$ are the
virtual element approximations of $\us$, $\HONEzr(\Omega)$,
$\as(\cdot,\cdot)$, and $\Fs(\cdot)$.
In the  rest
of this section we recall the construction of these mathematical objects.

\medskip
\PGRAPH{Mesh notation.}
Let $\calT=\{\Th\}_{\hh\in\calH}$ be a family of decompositions $\Th$ of
the computational domain $\Omega$ into a finite set of nonoverlapping
polygonal elements $\P$.
Each  of 
the members $\Th$ of the family $\calT$ will be referred to as the \emph{mesh}.
The  \textit{mesh size} $\hh$, which also serves as subindex, is the maximum of the diameters of the mesh
elements, which is defined by $\hP=\sup_{\xv,\yv\in\P}\ABS{\xv-\yv}$.
We assume that the mesh sizes of the mesh family $\calT$ are in a
countable subset $\calH$ of the real line $(0,+\infty)$ having $0$ as its
unique accumulation point.
We let  $\partial\P$ denote the polygonal boundary of $\P$, which we assume to be  nonintersecting and 
 formed by straight edges $\E$. The center of gravity of $\P$ will be denoted by 
$\xvP=(\xsP,\ysP)$ and its area by $\mP$.
We denote the edge mid-point $\xvE=(\xsE,\ysE)$ and its lenght
$\mE$, and, with a small abuse of notation, we write $\E\in\partial\P$
to indicate that edge $\E$ is running throughout the set of edges
forming the elemental boundary $\partial\P$.
In the next section we will discuss in detail the different assumptions on the mesh family $\calT$, under which he convergence analysis of the VEM and the derivation of the error
estimates in the $\LTWO$ and $\HONE$ are carried out in the literature.

\medskip
\PGRAPH{The virtual element spaces.}
Let $k\geq1$  integer and $\P\in\Th$ a generic mesh
element.
We define the local virtual element space $\Vhk(\P)$, following to 
 the \emph{enhancement strategy} proposed
in~\cite{Ahmad-Alsaedi-Brezzi-Marini-Russo:2013}:
\begin{multline}
\Vhk(\P) = \bigg\{\,
\vsh\in\HONE(\P)\,:\,
\restrict{\vsh}{\partial\P}\in\CS{0}(\partial\P),\,
\restrict{\vsh}{\E}\in\PS{k}(\E)\,\forall\E\in\partial\P,\,\\
\Delta\vsh\in\PS{k}(\P),\,
\textrm{and~}\\[0.25em]
\int_{\P}(\vsh-\PinP{k}\vsh)\,\qs\dV=0
\,\,\forall\qs\in\PS{k}(\P)\backslash\PS{k-2}(\P)
\,\bigg\}.
\label{eq:VhkP:def}
\end{multline}
Here, $\PinP{k}$ is the elliptic projection that will be discussed in
the next section; $\PS{k}(\P)$ and $\PS{k}(\E)$ are the linear spaces
of the polynomials of degree at most $k$, which are respectively
defined over an element $\P$ or an edge $\E$ according to our
notation; and $\PS{k}(\P)\backslash\PS{k-2}(\P)$ is the space of
polynomials of degree equal to $k-1$ and $k$.
By definition, the space $\Vhk(\P)$ contains $\PS{k}(\P)$ and the global
space $\Vhk$ is a conforming subspace of $\HONE(\Omega)$.
The global \emph{conforming virtual element space $\Vhk$ of order $k$ built
on mesh $\Th$} is obtained by gluing together the elemental
approximation spaces, that is
\begin{align}
  \Vhk:=\Big\{\,\vsh\in\HONEzr(\Omega)\,:\,\restrict{\vsh}{\P}\in\Vhk(\P)
  \,\,\,\forall\P\in\Th\,\Big\}.
  \label{eq:Vhk:def}
\end{align}
On every mesh $\Th$, given an integer $k\geq0$, we also define the space of discontinuous piecewise polynomials of degree at most $k$, $\PS{k}(\Th)$,
whose elements are the functions $\qs$ such that $\restrict{\qs}{\P}\in\PS{k}(\P)$ for
every $\P\in\Th$.

\medskip
\PGRAPH{The degrees of freedom.}
For each element $\P$ and each virtual element function
$\vsh\in\Vhk(\P)$, we consider the following set of degrees of freedom~\cite{BeiraodaVeiga-Brezzi-Cangiani-Manzini-Marini-Russo:2013}:
\begin{description}
	\item[]\textbf{(D1)} for $k\geq1$, the values of $\vsh$ at the
	vertices of $\P$;
	
	\medskip
	\item[]\textbf{(D2)} for $k\geq2$, the values of $\vsh$ at the $k-1$
	internal points of the $(k+1)$-point Gauss-Lobatto quadrature rule
	on every edge $\E\in\partial\P$.
	
	\medskip
	\item[]\textbf{(D3)} for $k\geq2$, the cell moments of $\vsh$ of order
	up to $k-2$ on element $\P$:
	\begin{align}
	\frac{1}{\mP}\int_{\P}\vsh\,\qs\dV
	\quad\forall\qs\in\PS{k-2}(\P).
	\label{eq:dofs:D3}
	\end{align}
\end{description}
It is possible to prove that this set of values is unisolvent in $\Vhk(\P)$,
cf.~\cite{BeiraodaVeiga-Brezzi-Cangiani-Manzini-Marini-Russo:2013};
hence, every virtual element function is uniquely identified by it.
The global degrees of freedom of a virtual element function in the 
space $\Vhk$ are given by collecting the elemental degrees of freedom 
\textbf{(D1)}-\textbf{(D3)} for all vertices, edges and elements.
Their unisolvence in $\Vhk$ is an immediate consequence of their
unisolvence in every elemental space $\Vhk(\P)$.

\medskip
\PGRAPH{The elliptic projection operators.}
The \emph{elliptic
projection operator} $\PinP{k}:\HONE(\P)\to\PS{k}(\P)$, whose definition is required in~\eqref{eq:VhkP:def} and which will be instrumental in the definition of the bilinear form $a_h$ in the following, is given,  for any
$\vsh\in\Vhk(\P)$,  by:
\begin{align}  
  \int_{\P}\nabla\PinP{k}\vsh\cdot\nabla\qs\dV  &= \int_{\P}\nabla\vsh\cdot\nabla\qs\dV\quad\forall\qs\in\PS{k}(\P),\label{eq:PiFn}\\[0.5em]
  \int_{\partial\P}\big(\PinP{k}\vsh-\vsh\big)\dS &= 0.                                                              \label{eq:def:Pib_k}
\end{align}
Equation~\eqref{eq:def:Pib_k} allows us to remove the kernel of the
gradient operator from the definition of $\PinP{k}$, so that the
$k$-degree polynomial $\PinP{k}\vsh$ is uniquely defined for every
virtual element function $\vsh\in\Vhk(\P)$.
Furthermore, projector $\PinP{k}$ is a polynomial-preserving operator,
i.e., $\PinP{k}\qs=\qs$ for every $\qs\in\PS{k}(\P)$.
By assembling elementwise contributions we  define a global projection operator
$\Pin{k}:\HONE(\Omega)\to\PS{k}(\Th)$, which is such that
$\restrict{\Pin{k}\vsh}{\P}=\PinP{k}(\restrict{\vsh}{\P}) \ \forall \P \in \Th$.
A key property of the elliptic projection operator is that the
projection $\PinP{k}\vsh$ of any virtual element function
$\vsh\in\Vhk(\P)$ is computable from the degrees of freedom 
\textbf{(D1)}-\textbf{(D3)} 
of $\vsh$
associated with element $\P$.

\medskip
\PGRAPH{Orthogonal projections.}
From the degrees of freedom of a virtual element function
$\vsh\in\Vhk(\P)$ we can also compute the $L^2(\P)$ orthogonal projections
$\PizP{k}\vsh\in\PS{k}(\P)$,
cf.~\cite{Ahmad-Alsaedi-Brezzi-Marini-Russo:2013}.
In fact,  $\PizP{k}\vsh$ of a function
$\vsh\in\Vhk(\P)$ is, by definition, the solution of the variational problem:
\begin{align}
  \int_{\P}\Piz{k}\vsh\,\qs\dV =
  \int_{\P}\vsh\,\qs\dV\qquad\forall\qs\in\PS{k}(\P).
\end{align}
The right-hand side is the integral of $\vsh$ against the polynomial
$\qs$, and,  when $\qs$ is a polynomial of degree up to $k-2$, it is computable from the degrees of freedom \textbf{(D3)} of
$\vsh$.  When $\qs$ is a polynomial of degree $k-1$
or $k$, it is computable from the
moments of $\PinP{k}\vsh$ , cf.~\eqref{eq:VhkP:def}.
Clearly, the orthogonal projection $\PizP{k-1}\vsh$ is also
computable.
Also here we can define a
global projection operator $\Piz{k}:\LTWO(\Omega)\to\PS(\Th)$, which
projects the virtual element functions on the space of discontinuous
polynomials of degree at most $k$ built on mesh $\Th$.
This operator is obtained by assembling the elemental $\LTWO$-orthogonal
projections $\PizP{k}\vsh$ for all mesh elements $\P$, that is
$\restrict{\big(\Piz{k}\vsh\big)}{\P}=\PizP{k}(\restrict{\vsh}{\P})$.

\medskip
\PGRAPH{The virtual element bilinear forms.}
Taking advantage of the elliptic and orthogonal projectors we can now define the
virtual element bilinear form
$\ash(\cdot,\cdot):\Vhk\times\Vhk\to\REAL$.
We start by splitting the discrete bilinear form
$\ash(\cdot,\cdot)$ as the sum of elemental contributions
\begin{align}
  \ash(\ush,\vsh) = \sum_{\P\in\Th}\ashP(\ush,\vsh),\quad
\end{align}
where the contribution of every element is a bilinear form
$\ashP(\cdot,\cdot):\Vhk(\P)\times\Vhk(\P)\to\REAL$, 
approximating the corresponding elemental bilinear form
$\asP(\cdot,\cdot):\HONE(\P)\times\HONE(\P)\to\REAL$,
\begin{align*}
  \asP(\vs,\ws) = \int_{\P}\nabla\vs\cdot\nabla\ws\dV,
  \quad\forall\vs,\ws\in\HONE(\P).
\end{align*}
The bilinear form $\ashP(\cdot,\cdot)$ on each element $\P$ is itself split as the sum of two contributions: 
by
\begin{align}
  \label{eq:ashP:def}
  \ashP(\ush,\vsh) = 
  &\int_{\P}\nabla\PinP{k}\ush\cdot\nabla\PinP{k}\vsh\dV + \\
  &\SPh\Big( \big(I-\PinP{k}\big)\ush, \big(I-\PinP{k}\big)\vsh \Big).\notag
\end{align}
The \textit{stabilization term}  $\SPh(\cdot,\cdot)$ in \eqref{eq:ashP:def} can be any
computable, symmetric, positive definite bilinear form defined on
$\Vhk(\P)$  such that
\begin{align}
  \cbot\asP(\vsh,\vsh)
  \leq\SPh(\vsh,\vsh)
  \leq\ctop\asP(\vsh,\vsh)
  \quad\forall\vsh\in\Vhk(\P)\cap\KER\big(\PinP{k}\big),
  \label{eq:SP:stability}
\end{align}
for two positive constants $\cbot$ and
$\ctop$.
Property \eqref{eq:SP:stability} states that that
$\SPh(\cdot,\cdot)$ scales like $\asP(\cdot,\cdot)$ with respect to
$\hP$.
As $\PinP{k}$ is a polynomial preserving operator, the contribution of the stabilization term in the definition of $\ashP(\cdot,\cdot)$
is zero one (or both) of the two entries is a polynomial of degree
less than or equals to $k$.

\medskip
Condition~\eqref{eq:SP:stability}, is designed so
that $\ashP(\cdot,\cdot)$ satisfies the two following fundamental properties:
\begin{description}
	\item[-] {\emph{$k$-consistency}}: for all $\vsh\in\Vhk(\P)$ and for all
	$\qs\in\PS{k}(\P)$ it holds that
	\begin{align}
	\label{eq:k-consistency}
	\ashP(\vsh,\qs) = \asP(\vsh,\qs);
	\end{align}
	
	\medskip
	\item[-] {\emph{stability}}: there exist two positive constants
	$\alpha_*,\,\alpha^*$, independent of $\hh$ and $\P$, such that
	\begin{align}
	\label{eq:stability}
	\alpha_*\asP(\vsh,\vsh)
	\leq\ashP(\vsh,\vsh)
	\leq\alpha^*\asP(\vsh,\vsh)\quad\forall\vsh\in\Vhk(\P).
	\end{align}
\end{description} 

\medskip
In our implementation of the VEM, we consider the stabilization
proposed in~\cite{Mascotto:2018}, that is we set:
\begin{align}
  \SPh(\vsh,\wsh) =
  \sum_{i=1}^{\NDOFS}\mathcal{A}^{\P}_{ii}\textrm{DOF}_i(\vsh)\textrm{DOF}_i(\wsh),
  \label{eq:stab:D-recipe}
\end{align}
where $\mathcal{A}^{\P}=\big(\mathcal{A}^{\P}_{ij}\big)$ is the matrix
stemming from the implementation of the first term in the bilinear
form $\ashP(\cdot,\cdot)$.
More precisely, let $\phi_i$ be the $i$-th ``canonical'' basis functions generating
the virtual element space, which is the function in $\Vhk(\P)$ whose
$i$-th degree of freedom for $i=1,\ldots,\NDOFS$ (according to a
suitable renumbering of the degrees of freedom in \textbf{(D1)},
\textbf{(D2)}, and \textbf{(D3)}), has value equal to $1$ and all
other degrees of freedom are zero. The $i,j$-th entry of matrix $\mathcal{A}^{\P}$ is
given by
\begin{align}
  \mathcal{A}^{\P}_{ij}:=\asP\big(\PinP{k}\phi_i,\PinP{k}\phi_j\big).
\end{align}
The stabilization in~\eqref{eq:stab:D-recipe} is sometimes called the
``\emph{D-recipe stabilization}'' in the virtual element literature. Observe that, if we instead replace $\mathcal{A}^P$ with the identity matrix, 
we get the  the so called ``\emph{dofi-dofi (dd)
stabilization}'' originally proposed
in~\cite{BeiraodaVeiga-Brezzi-Cangiani-Manzini-Marini-Russo:2013}:
\begin{align}
  \SPhdd(\vsh,\wsh) =
  \sum_{i=1}^{\NDOFS}\textrm{DOF}_i(\vsh)\textrm{DOF}_i(\wsh).
  \label{eq:stab:dofi-dofi}
\end{align}
The stabilization~\eqref{eq:stab:dofi-dofi} underlies 
many of the convergence results available from the literature, which
we briefly review in Section~\ref{subsec:star:error}.

\medskip
\PGRAPH{The virtual element forcing term.}
We approximate the right-hand side of~\eqref{eq:VEM}, by splitting it into
the sum of elemental contributions and approximating every local linear functional by replacing $\vsh$ with $\PizP{k}\vsh$:
\begin{align}
  \Fs(\vsh) = \sum_{\P\in\Th}\big(\fs,\PizP{k}\vsh\big)_{\P}.
  \textrm{~~where~~}\big(\fs,\PizP{k}\vsh\big)_{\P} = \int_{\P}\fs\,\PizP{k}\vsh\dV.
  \label{eq:extensionc}
\end{align}

\medskip
\PGRAPH{Approximation properties in the virtual element space.}
Under a suitable regularity assumption on the mesh family $\mathcal{T}$ (see assumption \textbf{G\ref{g1}} in
the next section), we can prove the following estimates on the projection and interpolation 
operators:
\begin{enumerate}
	\item for every $s$ with $1\le\ss\le\ks+1$ and for every $\ws\in\HS{s}(E)$ there exists a
	$\ws_{\pi}\in\PS{k}(\P)$ such that
	\begin{equation}
	\ABS{\ws-\ws_{\pi}}_{0,\P} +
	\hP\ABS{\ws-\ws_{\pi}}_{1,E}\leq\Cs\hh^s_{\P}\ABS{\ws}_{s,\P};
	\end{equation}
	\item for every $s$ with $2\le\ss\leq\ks+1$, for every $\hh$, for all $\P\in\Th$ and for
	every $\ws\in\HS{s}(E)$ there exists a $\wsI\in\Vhk(\P)$ such that
	\begin{equation}
	\ABS{\ws-\wsI}_{0,\P} +
	\hP\ABS{\ws-\wsI}_{1,\P}\leq\Cs\hP^s\ABS{\ws}_{s,\P}.
	\end{equation}  
\end{enumerate}
Here, $C$ is a  positive constant depending only on the polynomial degree 
$k$ and on some mesh regularity constants that we will present and discuss in the next
section.

\medskip
\PGRAPH{Main convergence properties.}
Thanks to the coercivity and continuity of the bilinear form $\ash(\cdot,\cdot)$, and to the
continuity of the right-hand side linear functional
$\big(\fs,\Piz{k}\cdot\big)$, by applying  the Lax-Milgram
theorem~\cite[Section 2.7]{Brenner-Scott:2008}, we obtain the well-posedness of the discrete formulation~\eqref{eq:VEM}.

\medskip
Let then $\us\in\HS{k+1}(\Omega)$ be the solution to the variational
problem~\eqref{eq:Poisson:weak} on a convex domain $\Omega$ with
$\fs\in\HS{k}(\Omega)$.
Let $\ush\in\Vhk$ be the solution of the virtual element method
\eqref{eq:VEM} on every mesh of a mesh family $\calT=\{\Th\}$
satisfying a suitable set of mesh geometrical assumptions.
Under suitable assumptions on the mesh familty $\mathcal{T}$, it is possible to prove that the following error estimates hold:
\begin{itemize}
\item the $\HONE$-error estimate holds:
  \begin{align}
    \label{eq:source:problem:H1:error:bound}
    \norm{\us-\ush }{1,\Omega}\leq
    \Cs\hh^{k}\left(
    \norm{\us}{k+1,\Omega} 
    + \snorm{\fs}{k,\Omega}
    \right);
  \end{align}

  \medskip
\item the $\LTWO$-error estimate holds:
  \begin{align}
    \label{eq:source:problem:L2:error:bound}
    \norm{\us-\ush}{0,\Omega}\leq 
    \Cs\hh^{k+1}\left(
    \norm{\us}{k+1,\Omega}
    + \snorm{\fs}{k,\Omega}
    \right).
  \end{align}
\end{itemize}
Constant $\Cs$ in \eqref{eq:source:problem:H1:error:bound} and in \eqref{eq:source:problem:L2:error:bound} may depend on the stability constants $\alpha_*$ and
$\alpha^*$, on mesh regularity constants which we will
introduce in the next section, on the size of the computational domain
$\ABS{\Omega}$, and on the approximation degree $k$.
Constant $\Cs$ is normally independent of $\hh$, but for the most
extreme meshes it may depend on the ratio between the longest and
shortest edge lenghts, cf. Section~\ref{subsec:star:error}.

\medskip
Finally, we note that the approximate solution $\ush$ is not
explicitly known inside the elements.
Consequently, in the numerical experiments of
Section~\ref{subsec:violating:performance}, we approximate the error
in the $\LTWO$-norm as follows:
\begin{align*}
  \norm{\us-\ush}{0,\Omega}\approx
  \norm{\us-\Piz{k}\ush}{0,\Omega},
\end{align*}
where $\Piz{k}\ush$ is the global $\LTWO$-orthogonal projection of the
virtual element approximation $\ush$ to $\us$.
On its turn, we approximate the error in the energy norm as follows:
\begin{align*}
  \snorm{ \us-\ush }{1,\Omega}\approx\ash(\usI-\ush,\usI-\ush)
\end{align*}
where $\usI$ is the virtual element interpolant of the exact solution
$\us$.

{In this work, we are interested in checking whether 
	optimal convergence rates  put forward by these estimates are maintained on different mesh families that may display
	some pathological situations.
	From a theoretical viewpoint, the convergence estimates hold under
	some constraints on the shapes of the elements forming the
	mesh, called \emph{mesh geometrical (or regularity) assumptions}.
	We summarize the major findings from the literature in
	Section~\ref{subsec:star:error} and in the next sections we will
	investigate how breaking such constraints may affect these results.}\section{State of the art}
\label{sec:star}
Various geometrical (or \textit{regularity}) assumptions have been
proposed in the literature to ensure that all elements of all meshes
in a given mesh family in the refinement process are sufficiently
regular.
These assumptions guarantee the convergence of the VEM and optimal
estimates of the approximation error with respect to different
norms.\\
In this section, we overview the geometrical assumptions introduced in
the VEM literature to guarantee the convergence of the method, and we
provide a list of the main convergence results based on such
assumptions.


\subsection{Geometrical assumptions}
\label{subsec:star:assumptions}
We start by reviewing the geometrical assumptions appeared in the VEM
literature since their definition in
\cite{BeiraodaVeiga-Brezzi-Cangiani-Manzini-Marini-Russo:2013}.
Note that these assumptions are defined for a single mesh $\Th$, but
the conditions contained in them are required to hold independently of
the mesh size $h$.
As a consequence, when an assumption is imposed to a mesh family
$\calT=\{\Th\}_h$, it has to be verified simultaneously by every $\Th
\in \calT$.\\

It is well-known from the FEM literature that the approximation
properties of a method depend on specific assumptions on the geometry
of the mesh elements.
Classical examples of geometrical assumptions for a family of
triangulations $\{\Th\}_{\hh\to 0}$, are the ones introduced in
\cite{ciarlet2002finite} and \cite{zlamal1968finite}, respectively:

\smallskip
\begin{description}
\item[$(a)$] \textit{Shape regularity condition:} there exists a real
  number $\gamma\in(0,1)$, independent of $h$, such that for any
  triangle $\P \in \Omega_h$ we have
  \begin{align} \label{eq:shape:regularity}
    \rP\ge\gamma h_\P,
  \end{align}
  being $h_\P$ the longest edge in $\P$ and $\rP$ its inradius;
\item[$(b)$] \textit{Minimum angle condition:} there exists an angle
  $\alpha_0>0$, independent of $h$, such that for any triangle
  $\P\in\Omega_h$ we have
  \begin{align} \label{eq:min:angle}
    \alpha_\P\ge\alpha_0,
  \end{align}
  being $\alpha_\P$ the minimal angle of $\P$.
\end{description}

\smallskip
When we turn our focus on polygonal meshes, a preliminar consideration
is needed on the definition of the polygonal elements.
It is commonly accepted, even if not always explicitly specified, that
a mesh $\Th$ has to be made of a finite number of \textit{simple
  polygons}, i.e. open simply connected sets whose boundary is a
non-intersecting line made of a finite number of straight line
segments.\\
The other regularity assumptions proposed in the VEM literature to
ensure approximation properties have been deduced in analogy to the
similar conditions developed for the FEMs.
The main assumption, which systematically recurs in every VEM paper,
is the so-called \textit{star-shapedness} of the mesh elements.

\smallskip
\begin{gass} \label{g1}
  there exists a real number $\rho\in(0,1)$, independent of $h$, such
  that every polygon $\P\in\Th$ is star-shaped with respect to a disc
  with radius
  \[\rP\ge\rho\hP.\]
\end{gass}
We denote $\rP$ the radius of the greatest possible inscribed disk in
$\P$ and \textit{star center} the center of such disk, where it
exists.
We stress the fact that \textbf{G1} does not accept polygons
star-shaped with respect to a single point, as both $\rho$ and $\hP$
are greater than zero, and we conventionally say that $\rP=0$ if $\P$
is not star-shaped.
\begin{figure}[htbp] 
  \centering
  \begin{tabular}{c c}
    \includegraphics[width=.2\columnwidth]{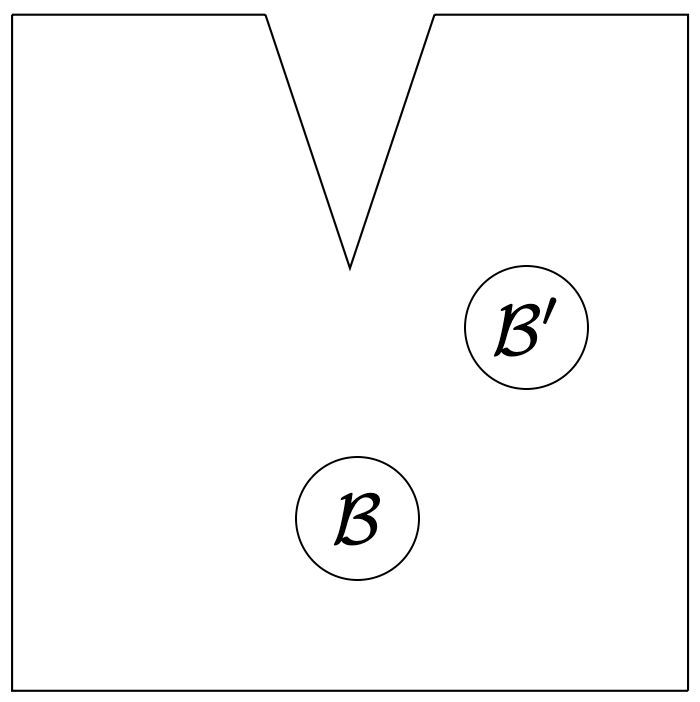}
    &\qquad
    \includegraphics[width=.2\columnwidth]{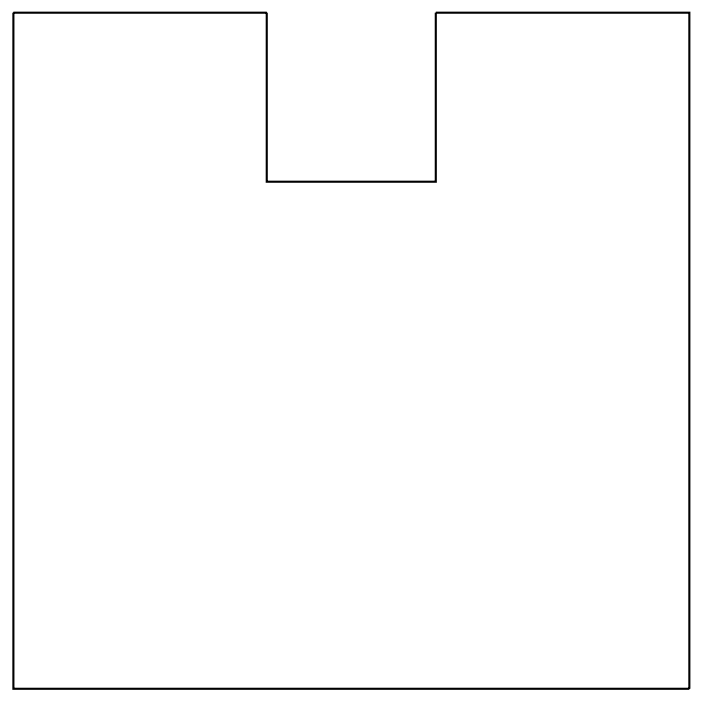}
    \\
    (a) &\qquad (b)
  \end{tabular}
  \caption{The element in (a) is star-shaped with respect to the disk
    $\calB$ but not with respect to the disk $\calB'$, while the
    element in (b) is not star-shaped with respect to any disk.}
  \label{fig:star}
\end{figure}\\
Assumption~\textbf{G\ref{g1}} is nothing but the polygonal extension
of the classical \textit{shape regularity condition} for triangular
meshes.
In fact, any triangular element $\P$ is star-shaped with respect to
its maximum inscribed disk (the one with radius $\rP$) and the
diameter $\hP$ coincides with its longest edge.
Moreover, \textbf{G\ref{g1}} can also be stated in the following weak
form, as specified in
\cite{BeiraodaVeiga-Brezzi-Cangiani-Manzini-Marini-Russo:2013} and
more accurately in \cite{Brenner:2017:SEV}:

\medskip
\noindent
\textbf{Assumption G1-weak} \textit{there exists a
  real number $\rho\in(0,1)$, independent of $h$, such that 
  every polygon $\P\in\Th$ can be split into a finite number $\Ns$ of disjoint
  polygonal subcells $\P_1,\dots,\P_{\Ns}$ where, for
  $j=1,\dots,\Ns$,
  \begin{description}
  \item[$(a)$] element $\P_j$ is star-shaped with respect to a disc with radius $r_{\P_j}\ge\rho\hh_{\P_j}$;
  \item[$(b)$] elements $\P_j$ and $\P_{j+1}$ share a common edge.
\end{description}}
\noindent
From a practical point of view, assumptions \textbf{G\ref{g1}} and
\textbf{G1-weak} are almost equivalent, and they are treated
equivalently in all papers reviewed in
Section~\ref{subsec:star:error}.

\smallskip
\noindent
Assumption~\textbf{G\ref{g1}} plays a key role in most of the
theoretical results regarding polygonal methods.
It is needed by the Bramble-Hilbert lemma \cite{Brenner-Scott:2008},
an important result on polynomial approximation that is often used for
building approximation estimates, and also by the following lemma.
\begin{lemma}\label{lemma:g1}
  If a mesh $\Th$ satisfies Assumption~\textbf{G\ref{g1}} then for all
  polygons $\P\in\Th$ there exists a mapping $F:\calB_1\to\P$, with
  the Jacobian $J$ of $F$ satisfying
  \begin{align}
    \norm{J}{2}\lesssim\hh, \qquad
    \snorm{\det(J)}{}\lesssim\hh^2 \qquad \textit{and} \qquad
      \norm{J^{-1}}{2}\lesssim\hh^{-1},
  \end{align}
  and, for a sufficiently regular $\us$, the following relations hold
  \begin{align*}
    &\norm{\us}{0,\P} \simeq \hP\norm{\us\circ F}{0,\calB_1},
    &\norm{\us}{0,\bP} \simeq \hP^{1/2}\norm{\us\circ F}{0,\partial\calB_1},\\
    &\snorm{\us}{1,\P} \simeq \snorm{\us\circ F}{1,\calB_1},
    &\snorm{\us}{1/2,\bP} \simeq \snorm{\us\circ F}{1/2,\partial\calB_1},
  \end{align*}
  where all the implicit constants only depend on the constant $\rho$
  from \textbf{G\ref{g1}}.
\end{lemma}
Thanks to the relations in Lemma~\ref{lemma:g1}, inequalities that we
have on the unit circle $\calB_1$, such as the Poincaré inequality or
the trace inequalities, may be transferred to the polygon $\P$ by a
``scaling'' argument.

\medskip
In the very first VEM paper
\cite{BeiraodaVeiga-Brezzi-Cangiani-Manzini-Marini-Russo:2013}, where
the method was introduced, Assumption~\textbf{G\ref{g1}} was followed
by another condition on the maximum point-to-point distance.

\smallskip
\noindent
\textbf{Assumption G2-strong} \textit{there exists a real number
  $\rho\in(0,1)$, independent of $h$, such that for every polygon
  $\P\in\Th$, the the distance $d_{i,j}$ between any two vertices
  $v_i, v_j$ of $\P$ satisfies}
  \[d_{i,j}\ge\rho\hP.\]

\smallskip
\noindent
In fact, Assumption~\textbf{G2-strong} was soon replaced in the
following works \cite{Ahmad-Alsaedi-Brezzi-Marini-Russo:2013},
\cite{Beirao-Lovadina-Russo:2016} and \cite{Brenner:2017:SEV} by a
weaker condition on the length of the elemental edges.
This new version allows, for example, the existence of four-sided
polygons with equal edges but one diagonal much smaller than the
other.

\smallskip
\begin{gass} \label{g2}
  there exists a real number $\rho\in(0,1)$, independent of $h$, such
  that for every polygon $\P\in\Th$, the length $\mE$ of every edge
  $\E\in\bP$ satisfies
  \[\mE\ge\rho\hP.\]
\end{gass}
The Authors consider a single constant $\rho$ for both
Assumption~\textbf{G\ref{g1}} and \textbf{G\ref{g2}} and refer to it
as the \textit{mesh regularity constant} or \textit{parameter}.\\
Under Assumption~\textbf{G1-weak} and \textbf{G\ref{g2}}, it can be
proved \cite{Brenner:2017:SEV} that the simplicial triangulation of
$\P$ determined by the star-centers of $\P_1,\dots,\P_{\Ns}$ satisfies
the \textit{shape regularity} and the \textit{minimum angle}
conditions.
The same holds under Assumptions \textbf{G\ref{g1}} and
\textbf{G\ref{g2}}, as a special case of the previous statement.
Moreover, Assumption~\textbf{G\ref{g2}} implies that for $1\le j,k\le
N$ it holds ${h_{\P_{j}}}\slash{h_{\P_{k}}}\le\rho^{-|j-k|}$.\\
As already mentioned in the very first papers, these assumptions are
more restrictive than necessary, but at the same time they allow the
VEM to work on very general meshes.
For example, Ahmad et
al. in~\cite{Ahmad-Alsaedi-Brezzi-Marini-Russo:2013} state that:
\begin{quote}
  \textit{``Actually, we could get away with even more general
    assumptions, but then it would be long and boring to make precise
    (among many possible crazy decompositions that nobody will ever
    use) the ones that are allowed and the ones that are not.''}
\end{quote} 

\medskip
In the subsequent papers \cite{Beirao-Lovadina-Russo:2016} and
\cite{brenner2018virtual} assumption~\textbf{G\ref{g1}} is preserved,
but assumption~\textbf{G\ref{g2}} is substituted by the alternative
version:

\smallskip
\begin{gass}  \label{g3}
  There exists a positive integer $\Ns$, independent of $\hh$, such
  that the number of edges of every polygon $\P\in\Th$ is (uniformly)
  bounded by $\Ns$.
\end{gass}

\smallskip
\noindent
The Authors show how assumption~\textbf{G\ref{g2}} implies
assumption~\textbf{G\ref{g3}}, but assumption~\textbf{G\ref{g3}} is
weaker than assumption~\textbf{G\ref{g2}}, as it allows for edges
arbitrarily small with respect to $\hP$.
Both combinations~\textbf{G\ref{g1}}+\textbf{G\ref{g2}}
and~\textbf{G\ref{g1}}+\textbf{G\ref{g3}} imply that the number of
vertices of $\P$ and the minimum angle of the simplicial triangulation
of $\P$ obtained by connecting all the vertices of $\P$ to its
star-center, are controlled by the constant $\rho$.

\medskip
Another step forward in the refinement of the geometrical assumptions
was made by Beir\~ao Da Veiga et al. in \cite{da2020sharper}.
Besides assuming \textbf{G\ref{g1}}, the Authors imagine to "unwrap"
the boundary $\bP$ of each element $\P \in \Omega_h$ onto an interval
of the real line, obtaining a one-dimensional mesh $\calI_\P$.
The mesh $\calI_\P$ can be subdivided into a number of disjoint
sub-meshes $\calI_\P^1, \ldots, \calI_\P^N$, corresponding to the
edges of $\P$.
Then, the following condition is assumed.

\smallskip
\begin{gass} \label{g4}
  There exists a real number $\delta>0$, independent of $\hh$, such that for every polygon $\P\in\Th$:
  \begin{description}
  \item[$(a)$] the one-dimensional mesh $\calI_\P$ can be subdivided into a finite number $N$ of disjoint sub-meshes $\calI_\P^1, \ldots, \calI_\P^N$;
  \item[$(b)$] for each sub-mesh $\calI_\P^i$, $i=1, \ldots, N$, it holds that
  \[
  \frac{\max_{e\in\calI_\P^i}|e|}{\min_{e\in\calI_\P^i}|e|} 
  \le \delta.
  \]
  \end{description}
\end{gass}

\begin{figure}[htbp]  
  \centering
  \includegraphics[width=.8\columnwidth]{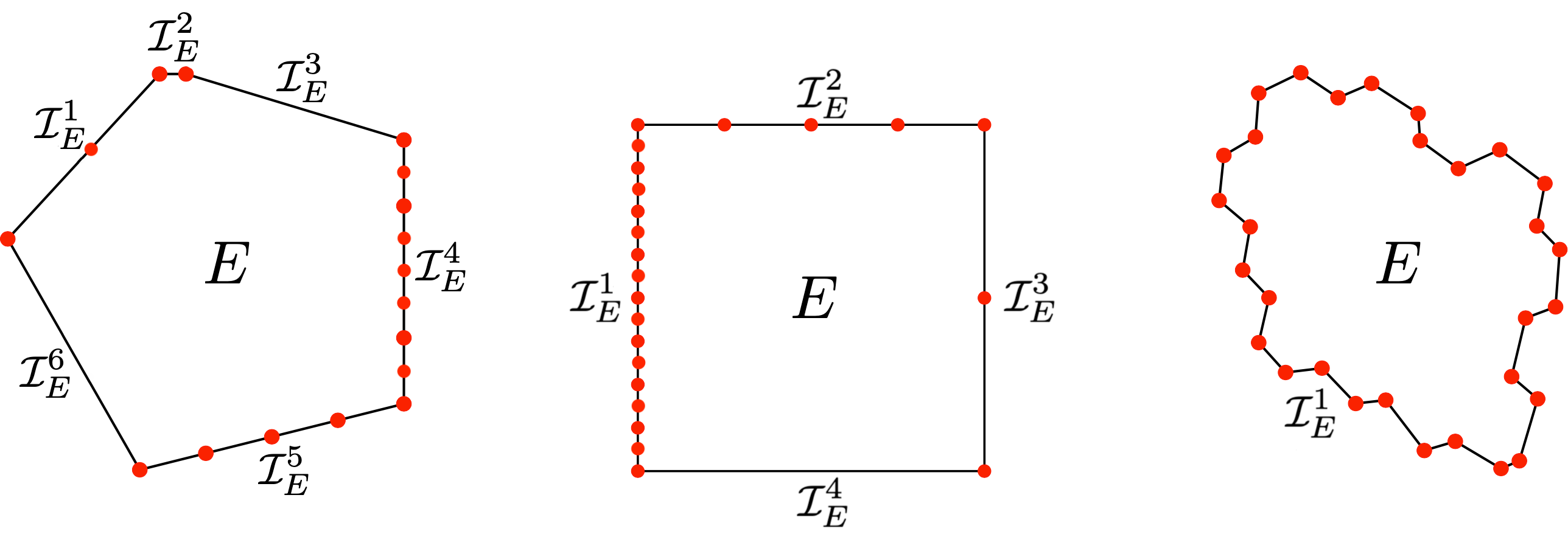}
  \caption{Examples of admissible elements according to
    assumption~\textbf{G\ref{g4}}. Red dots indicate the vertices of
    the element.}
    \label{fig:g4}
\end{figure}

\smallskip
Each polygon $\P$ corresponds to a
one-dimensional mesh $\calI_\P$, but a sub-mesh
$\calI_\P^i\subset\calI_\P$ might contain more than one edge of $\P$, cf. Fig.~\ref{fig:g4}.
Therefore assumption~\textbf{G\ref{g4}} does not require a
uniform bound on the number of edges in each element and does not
exclude the presence of small edges.
Mesh families created by agglomeration, cracking,
gluing, etc.. of existing meshes are admissible according
to~\textbf{G\ref{g4}}.
	
\smallskip
\noindent
As we will see in the next section, possible assumption pairs
requested in the literature to guarantee the convergence of the VEM
are given by combining \textbf{G\ref{g1}} (or, equivalently,
\textbf{G1-weak}) with either \textbf{G2-strong}, \textbf{G\ref{g2}},
\textbf{G\ref{g3}} or \textbf{G\ref{g4}}.

\medskip
Last, we report a condition (together with its weaker version) that
appears in the literature related to the \textit{non-conforming}
version of the VEM, which we do not cover in this work.
It is defined for instance in \cite{bertoluzzamanzini2021}, but it can
also be found under the name of \textit{local quasi-uniformity}.

\smallskip
\begin{gass}\label{g5} 
  there exists a real number $\gamma\in(0,1)$, independent of $h$,
  such that for every polygon $\P\in\Th$, for all $\E, \E'$ adjacent
  edges of $\P$ it holds that
  \begin{align*}
  \frac{1}{\gamma} \leq \frac{|\E|}{|\E'|} \leq\gamma.
  \end{align*}
\end{gass}
In this case, the polygonal elements of the mesh are allowed to have a
very large number of very small edges, provided that every two
consecutive edges scale proportionally.
This assumption also has a weak version, in which we essentially ask that, for every $\P\in\Th$, a part of $\bP$ satisfies \textbf{G\ref{g5}} and the remaining part satisfies \textbf{G\ref{g3}}.

\smallskip
\noindent
\textbf{Assumption G5-weak} \textit{there exist a real number $\gamma\in(0,1)$ and a positive integer $N$, both independent of $\hh$, such that for every polygon $\P\in\Th$, the set $\calE_\P$ of the edges of $\P$ can be split as $\calE_\P=\calE_\P^1\cup\calE_\P^2$, where $\calE_\P^1$ and $\calE_\P^2$ are such that
\begin{description}
\item[-] for any pair of adjacent edges $\E,\E'\in\calE_\P$ with $\E\in\calE_\P^1$, the inequality holds
  \[\frac{1}{\gamma} \leq \frac{|\E|}{|\E'|} \leq \gamma;\]
\item[-] $\calE_\P^2$ contains at most $N$ edges.
\end{description}}
Assumption \textbf{G5-weak} allows for situations where a large number of small edges coexists with some large edges.
We can think of families of meshes for which such an assumption is not satisfied, but they would be extremely pathological.
Assumptions \textbf{G5} and \textbf{G5-weak} are here reported for the sake of completeness but are not considered in the analysis of the following sections, as they refer to a different context.


\subsection{Convergence results in the VEM literature}
\label{subsec:star:error}
In this section, we briefly overview the literature on the main results of the convergence analysis of the VEM method.
%
For each article, we explicitly report (if available) the theoretical results and highlight the geometric assumptions used, reporting the abstract energy error, the $\HONE$ error estimate, and the $\LTWO$ error.
For a greater uniformity of the presentation with the rest of the chapter, we have slightly modified some notations and introduced minimal variations to some statements of the theorems.

\medskip
\PGRAPH{``Basic Principles of Virtual Elements Methods''
\cite{BeiraodaVeiga-Brezzi-Cangiani-Manzini-Marini-Russo:2013}}\\
This is the paper in which the VEM method was introduced and defined.
In the original formulation, the paper introduced the regular conforming virtual element space. For simplicity and with a small abuse of notation, the regular conforming virtual element space is denoted by $\Vhk$ as in \eqref{eq:Vhk:def} and
\eqref{eq:VhkP:def}:
\begin{subequations}
  \label{eq:Vhk:regular:conforming}
  \begin{align}
    \Vhk:= 
    &\{\vsh\in\HONE(\Omega):\ \restrict{\vsh}{\P}\in\Vhk(\P) \ \forall\P\in\Th\},
  \end{align}
  where
  \begin{align}
    \Vhk(\P):= 
    \{\vsh\in\HONE(\P):\ 
    &\restrict{\vsh}{\bP}\in C^0(\bP), \ 
    \restrict{\vsh}{\E}\in\PS{k}(\E) \ \forall\E\in\bP, \\
    &\text{and } \ \Delta\vsh\in\PS{k-2}(\P)\},\notag
  \end{align}
\end{subequations}
and the \textit{dofi-dofi} formulation $\SPhdd$ defined in
\eqref{eq:stab:dofi-dofi} is used for the stabilization bilinear
form.\\
There, the authors introduced the concept of \emph{simple polygon} and the geometric regularity assumptions \textbf{G\ref{g1}} and \textbf{G\ref{g2} - strong}. The following statement on the convergence of VEM in the energy norm is general and largely used in the VEM literature, even if not explicitly stated in \cite{BeiraodaVeiga-Brezzi-Cangiani-Manzini-Marini-Russo:2013}.
To this purpose, for functions
$\vs\in\HONE(\Th)$ the broken $\HONE$-seminorm is defined as follows:
\begin{align}    \label{eq:h1-broken}
    \snorm{\vs}{\hh,1}:=\left(\sum_{\P\in\Th} \snorm{\nabla \vs}{0,\P}^2\right)^{1/2}.
\end{align}
\begin{theo}[abstract energy error]  \label{theo:basic:abstract}
Under the \textit{k-consistency} and \textit{stability} assumptions
  defined in Subsection~\ref{sec:VEM}, cf. \eqref{eq:k-consistency} and
  \eqref{eq:stability}, the discrete problem has a unique solution
  $\ush$.
  Moreover, for every approximation $\usI\in \Vhk$ of $\us$ and every
  approximation $\us_\pi$ of $\us$ that is piecewise in $\PS{k}(\Th)$,
  we have
  \begin{equation}
    \snorm{\us-\ush}{1,\Omega} \leq
    \Cs(\snorm{\us-\usI}{1,\Omega} +
    \snorm{\us-\us_{\pi}}{\hh,1} +
    \mathfrak{F}_{\hh}),
  \end{equation}
  where $\Cs$ is a constant depending only on $\alpha_*$ and
  $\alpha^*$ (the constants in \eqref{eq:stability}), and, for any
  $\hh$, $\mathfrak{F}_{\hh}=\snorm{\fs-\fsh}{{\Vhk}^\prime}$ is the
  smallest constant such that
  \begin{equation*}
    (\fs,\vs)-\langle\fsh-\fs,\vs\rangle \leq
    \mathfrak{F}_{\hh}\snorm{\fs}{1}
    \qquad\forall\vs\in\Vhk.
  \end{equation*}
\end{theo}
In \cite{BeiraodaVeiga-Brezzi-Cangiani-Manzini-Marini-Russo:2013} it was claimed that it is possible to estimate the convergence rate in terms of the $\LTWO$ error using duality argument techniques.

\medskip
\PGRAPH{``Equivalent projectors for virtual element methods'' \cite{Ahmad-Alsaedi-Brezzi-Marini-Russo:2013}}\\
In \cite{Ahmad-Alsaedi-Brezzi-Marini-Russo:2013}, the enhanced VEM space
\eqref{eq:VhkP:def} adopted in this chapter replaces $\Vhk(\P)$ (in \cite{Ahmad-Alsaedi-Brezzi-Marini-Russo:2013} the enhanced VEM space is named ``modified VEM space''). This paper uses the \textit{dofi-dofi} stabilization.
Under the geometrical assumptions \textbf{G\ref{g1}} and
\textbf{G\ref{g2}}, the paper provides an explicit estimation of $\HONE$ and $\LTWO$ errors and introduces the Theorem \ref{theo:basic:abstract} for the abstract energy error.
\begin{theo}[$\HONE$ error estimate]
  Assuming~\textbf{G\ref{g1}},~\textbf{G\ref{g2}}, let the right-hand
  side $\fs$ belong to $H^{k-1}(\Omega)$, and that the exact solution
  $\us$ belong to $H^{k+1}(\Omega)$. Then
  \begin{equation}
    \norm{\us-\ush}{1,\Omega} \le \Cs|\hh|^k \snorm{\us}{k+1,\Omega}
  \end{equation}
  with $\Cs$ a positive constant independent of $\hh$.
\end{theo}
\begin{theo}[$\LTWO$ error estimate]
  Assuming~\textbf{G\ref{g1}},~\textbf{G\ref{g2}} and with $\Omega$
  convex, let the right-hand side $\fs$ belong to $H^k(\Omega)$, and
  that the exact solution $\us$ belong to $H^{k+1}(\Omega)$. Then
  \begin{equation}
    \norm{\us-\ush}{0,\Omega} + 
    |\hh|\norm{\us-\ush}{1,\Omega} \le
    \Cs|\hh|^{k+1}\snorm{u}{k+1,\Omega},
  \end{equation}
  with $\Cs$ a constant independent of $\hh$.
\end{theo}

\medskip
\PGRAPH{``Stability analysis for the virtual element method'' \cite{Beirao-Lovadina-Russo:2016}}\\
This contribution deals with the regular conforming VEM space
\eqref{eq:Vhk:regular:conforming} defined in
\cite{BeiraodaVeiga-Brezzi-Cangiani-Manzini-Marini-Russo:2013}.
The paper \cite{Beirao-Lovadina-Russo:2016} provides a new estimation of the abstract energy error and analyses the $\HONE$ error with respect to two different stabilization techniques.
Moreover, new analytical assumptions on the bilinear
form $\ash(\cdot,\cdot)$ to replace \eqref{eq:stability} are introduced:
\begin{subequations} 
\begin{align}
\label{eq:stability2} 
    \ashP(\vsh,\vsh) 
    &\le \Cs_1(\P)\Tnorm{\vsh}{\P}^2, 
    \text{ for all } \vsh\in\Vhk(\P);\\
    \Tnorm{\qs}{\P}^2 
    &\le \Cs_2(\P)\ashP(\qs,\qs), 
    \text{ for all } \qs\in\PS{k}(\P),
\end{align}
\end{subequations}
being $\Tnorm{\cdot}{\P}$ a discrete semi-norm induced by the
stability term and $\Cs_1(\P), \Cs_2(\P)$ positive constants which
depend on the shape and possibly on the size of $\P$.
In this paper, the estimate is necessary for the
polynomials $\qs\in\PS{k}(\P)$ only, while in the standard analysis in
\cite{BeiraodaVeiga-Brezzi-Cangiani-Manzini-Marini-Russo:2013} a
kind of bound \eqref{eq:stability2}(b) was required for every
$\vsh\in\Vhk(\P)$.
Thus, even when $\Cs_1(\P)$ and $\Cs_2(\P)$ can be chosen independent
of $\P$, on $\Vhk(\P)$ the semi-norm induced by the stabilization term
may be stronger than the energy $\ashP(\cdot,\cdot)^{1/2}$.\\
\begin{theo}[abstract energy error]
  Under the \textit{stability} assumptions \eqref{eq:stability2}, let
  the continuous solution $\us$ of the problem satisfy
  $\restrict{\us}{\P}\in \calV_\P$ for all $\P\in\Th$, where $\calV_\P
  \subseteq \Vhk(\P)$ is a subspace of sufficiently regular functions.
  Then, for every $\usI\in \Vhk$ and for every $\us_\pi$ such that
  $\us_{\pi|\P}\in\PS{k}(\P)$, the discrete solution $\ush$ satisfies
  \begin{equation}
    \snorm{\us-\ush}{1,\Omega} \lesssim
    \Cs_{\text{err}}(\hh)\ (\mathfrak{Fs}_{\hh}
    + \Tnorm{\us-\usI}{}
    + \Tnorm{\us-\us_{\pi}}{}
    + \snorm{\us-\usI}{1,\Omega}
    + \snorm{\us-\us_{\pi}}{\hh,1}),
  \end{equation}
  where the constant $\Cs_{\text{err}}(\hh)$ is given by
  \begin{equation*}
    \Cs_{\text{err}}(\hh) = \max\left\{
    1,\,\tilde{\Cs}(\hh)\Cs_1(\hh),\tilde{\Cs}(\hh)^{3/2}\sqrt{\Cs^*(\hh)\Cs_1(\hh)}
    \right\}.
  \end{equation*}
\end{theo}
From the previous theorem, the following quantities are derived
from the constants in \eqref{eq:stability2}:
\begin{align*}
    \tilde{\Cs}(\hh) = \max_{\P\in\Th}\{1,\,\Cs_2(\P)\},\quad
    \Cs_1(\hh) = \max_{\P\in\Th}\{\Cs_1(\P)\},\\
    \Cs^*(\hh) = \frac{1}{2}\max_{\P\in\Th}\{\min\{1,\Cs_2(\P)^{-1}\}\}.
\end{align*}
In \cite{Beirao-Lovadina-Russo:2016} the stability term $\SPh(\cdot,\cdot)$ is considered as the combination of two contributions: the first, $\Ss_h^{\bP}$, related to the
boundary degrees of freedom; the second, $\Ss_h^{o\P}$, related to the
internal degrees of freedom.
In the following statements, we restrict the
analysis to $\Ss_h^{\bP}$ without losing generality. In this case,$\Ss_h^{\bP}$ is expressed in the
\textit{dofi-dofi} form $S_{\hh}^{\bP,\textrm{dd}}$ 
defined in \eqref{eq:stab:dofi-dofi}, or in the \textit{trace} form
introduced in \cite{Wriggers-Rust-Reddy:2016}:
\begin{align}  \label{eq:stab:trace}
  S_{\hh}^{\bP,\textrm{tr}}(\vsh,\wsh) =
  \hh_{\P} \int_{\bP} \partial_s\vsh \partial_s\wsh ds,
\end{align}
where $\partial_s\vsh$ denotes the tangential derivative of $\vsh$
along $\bP$.
\begin{theo}[$\HONE$ error estimate with \textit{dofi-dofi} stabilization]
  Assuming~\textbf{G\ref{g1},~G\ref{g3}}, let $\us\in\HS{s}(\Omega)$,
  $s>1$, be the solution of the problem with
  $\SPh=\Ss_h^{\bP,\textrm{dd}}$.
  Let $\ush$ be the solution of the discrete problem, then it holds
  \begin{align}
    \norm{\us-\ush}{1,\Omega} \lesssim \Cs(\hh)\hh^{s-1}\snorm{\us}{s,\Omega}\quad 1<\ss\le\ks+1,
  \end{align}
  with
  \[\Cs(\hh)=\max_{\P\in\Th}(\log(1+\hP/\hh_m(\P))),\]
  where $\hh_m(\P)$ denotes the length of the smallest edge of $\P$.
\end{theo}
\begin{corollary}
  Assuming~\textbf{G\ref{g1}} and \textbf{G\ref{g2}} instead, then
  $\cs(\hh)\lesssim 1$ and therefore
  \begin{equation*}
    \norm{\us-\ush}{1,\Omega}\lesssim\hh^{s-1}\snorm{\us}{s,\Omega}\qquad 1<s\le\ks+1.
  \end{equation*}
\end{corollary}        
\begin{theo}[$\HONE$ error estimate with \textit{trace} stabilization]
  Under Assumption~\textbf{G\ref{g1}}, let $\us\in\HS{s}(\Omega)$,
  $s>3/2$ be the solution of the problem with
  $\SPh=\Ss_h^{\bP,\textrm{tr}}$.
  Let $\ush$ be the solution of the discrete problem, then it holds
  \begin{equation}
    \norm{\us-\ush}{1,\Omega}\lesssim\hh^{s-1}\snorm{\us}{s,\Omega}\quad 3/2<s\le k+1.
  \end{equation}
\end{theo}

\medskip
\PGRAPH{``Some Estimates for Virtual Element Methods'' \cite{Brenner:2017:SEV}}\\
In \cite{Brenner:2017:SEV}, the enhanced VEM space is defined as in the following:
\begin{equation}
    \label{eq:Vhk:enhanced2}
    \begin{split}
    \Vhk(\P) := \big\{ \vsh\in\HONE(\P):\
    &\restrict{\vsh}{\bP}\in\PS{k}(\bP),\\
    &\exists\,\qs_{\vsh}(=-\Delta\vsh)\in\PS{k}(\P)\text{~such~that~} \\
    &\int_{\P}\nabla\vsh\cdot\nabla\wsh\dxv = \int_{\P}\qs_{\vsh}\wsh\dxv\quad\forall\wsh\in\HONEzr(\P), \\
    &\text{~and~}
    \Pi_k^{0,\P}\vsh - \Pi_k^{\nabla, \P}\vsh\in\PS{k-2}(\P) \big\}.
    \end{split}
\end{equation}
i.e., in a slightly different but equivalent way than \eqref{eq:VhkP:def}. In this work, the Authors consider different types of stabilization, but the convergence results are independent of $\SPh$.
The geometrical assumptions used in \cite{Brenner:2017:SEV} are
\textbf{G\ref{g1}} and \textbf{G\ref{g2}}.
\begin{theo}[abstract energy error]
  Assuming~\textbf{G\ref{g1},~G\ref{g2}}, if $\fs\in\HS{s-1}(\Omega)$
  for $1\le\ss\le\ks$, then there exists a positive constant $\Cs$
  depending only on $k$ and $\rho$ from \textbf{G\ref{g1}} such that
  \begin{equation}
    \snorm{\us-\ush}{1,\Omega}\le\Cs(\inf_{\vs\in\Vhk}\snorm{\us-\vs}{1,\Omega}
    + \inf_{\ws\in\PS{k}(\Th)}\snorm{\us-\ws}{\hh,1}
    + \hh^s\snorm{\fs}{s-1,\Omega}).
  \end{equation}
\end{theo}
\begin{theo}[$\HONE$ error estimate]
  Assuming~\textbf{G\ref{g1},~G\ref{g2}}, if $\us\in\HS{s+1}(\Omega)$
  for $1\le\ss\le\ks$, then there exists positive constants $\Cs_1$,
  $\Cs_2$ depending only on $k$ and $\rho$ from \textbf{G\ref{g1}}
  such that
  \begin{equation}
    \snorm{\us-\ush}{1,\Omega} + \snorm{\us-\Pin{k}\ush}{\hh,1}
    \leq\Cs_1\hh^s\snorm{\us}{s+1,\Omega}.
  \end{equation}
\end{theo} 
\begin{theo}[$\LTWO$ error estimate]
  Assuming~\textbf{G\ref{g1},~G\ref{g2}}, with $\Omega$
  convex, if $\us\in\HS{s+1}(\Omega)$ for for $1\le\ss\le\ks$, then
  there exists a positive constant $\Cs$ depending only on $\Omega$,
  $k$ and $\rho$ from \textbf{G\ref{g1}} such that
  \begin{equation}
    \norm{\us-\ush}{0,\Omega}\le\Cs\hh^{s+1}\snorm{\us}{s+1,\Omega}.
  \end{equation}
\end{theo}

\medskip
\PGRAPH{``Virtual element methods on meshes with small edges or faces'' \cite{brenner2018virtual}}\\
In this paper, error estimates of the VEM are yield for polygonal or polyhedral meshes possibly equipped with small edges $(d = 2)$ or faces $(d = 3)$.
In this case, the VEM space is formulated as \eqref{eq:Vhk:enhanced2} (i.e., the so-called enhanced space). The local stabilizing bilinear form is
considered in both the \textit{dofi-dofi} ($\SPhdd$) and in the
\textit{trace} ($\SPhtr$of \eqref{eq:stab:trace}) formulations. The Authors introduce also the constants:
\begin{align}
  \label{eq:error:small}
  \calH := \sup_{\P\in\Th}\left(
  \frac{\max_{\E\in\bP}\hE}{\min_{\E\in\bP}\hE}
  \right), \qquad
  \alpha_{\hh}:= \begin{cases}
    \ln\left( 1+\calH \right) & \mbox{~with } \SPhdd\\
    1 & \mbox{~with } \SPhtr
  \end{cases}
\end{align}
The geometrical assumptions considered in this work are
\textbf{G\ref{g1}} and \textbf{G\ref{g3}}.
In particular, a mesh-dependent energy norm $\norm{\cdot}{\hh}
:= \sqrt{\ash(\cdot,\cdot)}$ and a functional
$\Xi_{\hh}:\Vhk\to\PS{k}(\Th)$ are introduced. The function $\Xi_{\hh}$ is defined as:
\begin{align}
  \Xi_{\hh} =
  \begin{cases}
    \Piz{1} & \mbox{ if } k=1,2\\
    \Piz{k-1} & \mbox{ if } k\ge3.
  \end{cases}
\end{align}
\begin{theo}[abstract energy error]
  Assuming~\textbf{G\ref{g1},~G\ref{g3}}, let $\us$ and $\ush$ be the
  solutions of the continuous and discrete problems. We have:
  \begin{align}
    \norm{\us-\ush}{\hh}\lesssim
    &\inf_{\ws\in\Vhk} \norm{\us-\ws}{h} +
    \norm{\us-\Pin{k}\us}{\hh} +\\
    &\sqrt{\alpha_{\hh}}\left(\norm{\us-\Pin{k}\us}{\hh,1} +
    \sup_{\ws\in\Vhk}\frac{(\fs,\ws-\Xi_{\hh}\ws)}{\snorm{w}{1,\Omega}}
    \right).
  \end{align}
\end{theo}
\begin{theo}[$\HONE$ error estimate]
  Assuming~\textbf{G\ref{g1},~G\ref{g3}}, if the solution $\us$
  belongs to $\HS{s+1}(\Omega)$ for some $1\le\ss\le\ks$, we have:
  \begin{align}
    &\norm{\us-\ush}{\hh} \lesssim \sqrt{\alpha_{\hh}}\hh^s\snorm{u}{s+1,\Omega}, \ \textit{and}\\
    &\snorm{\us-\ush}{1,\Omega}
    + \sqrt{\alpha_{\hh}}\left[
      \snorm{\us-\Pin{k}\ush}{\hh,1}
      + \snorm{\us-\Piz{k}\us}{\hh,1}
      \right]
    \lesssim \alpha_{\hh}\hh^s\snorm{\us}{s+1,\Omega}.
  \end{align}
\end{theo}
\begin{theo}[$\LTWO$ error estimate]
  Assuming~\textbf{G\ref{g1},~G\ref{g3}}, if the solution $\us$ belongs to \\
  $\HS{s+1}(\Omega)$ for some $1\le\ss\leq\ks$, we have:
  \begin{equation}
    \norm{\us-\ush}{0,\Omega}\leq
    \Cs\ \alpha_{\hh}\hh^{s+1}
    \snorm{\us}{s+1,\Omega}.
  \end{equation}
\end{theo} 
The notation $A\lesssim B$ denotes that $A\le\Cs B$, with a positive
constant $\Cs$ depending on: i) the mesh regularity parameter $\rho$ of
\textbf{G\ref{g1}}, ii) the degree $k$ in the case of $\SPhtr$, and iii)
the maximum number of edges $N$ of \textbf{G\ref{g3}} in the
case of $\SPhdd$.

\medskip
\PGRAPH{``Sharper error estimates for Virtual Elements and a bubble-enriched version'' \cite{da2020sharper}}\\
This paper shows that the $\HONE$ interpolation error
$\snorm{\us-\usI}{1,\P}$ on each element $\P$ can be split into two parts: a
boundary and a bulk contribution.
The intuition behind this work is that it is possible to decouple the polynomial order on the boundary and in the bulk of the element.
Let $\ks_o$ and $\ks_{\partial}$ be two positive integers with
$\ks_o\ge\ks_{\partial}$ and let $\kv=(\ks_o,\ks_{\partial})$.
For any $\P\in\Th$, the \textit{generalized virtual
  element space} is defined as follows:
\begin{subequations}
\begin{align}
    V^h_{\kv} := 
    &\{\vsh\in\HONEzr(\Omega):
    \restrict{\vsh}{\P}\in) V^h_{\kv}(\P), \ \forall\P\in\Th\},
\end{align}
where
\begin{align}
    V^h_{\kv}(\P) := 
    \{\vsh\in\HONEzr(\P):
    &\restrict{\vsh}{\bP}\in\CS{0}(\bP), \
    \restrict{\vsh}{\E}\in\PS{\ks_\partial}(\E) \ \forall\E\in\bP, \\
    &\text{and } \ \Delta\vsh\in\PS{\ks_o-2}(\P)\}.\notag
\end{align}
\end{subequations}
For $\ks_o=\ks_{\partial}$, the space $V^h_{\kv}(\P)$ coincides with
the regular virtual element space in
\eqref{eq:Vhk:regular:conforming}.
In addition, given a function $v \in \HONEzr \cap H^s(\Th)$, on each
element $\P \in \Th$ the interpolant function
$\calI_hv$ is defined as the solution of an elliptic problem as follows:
\begin{align*}
  \begin{cases} \Delta\calI_hv = \Pi^{0,\P}_{k_o-2} \Delta v& \mbox{in } \P \\
    \calI_hv = v_b & \mbox{on } \bP,
  \end{cases}
\end{align*}
where $v_b$ is the standard 1D piecewise polynomial interpolation of
$v|_{\bP}$.
\begin{theo}[abstract energy error]
  Under Assumption~\textbf{G\ref{g1}}, let $u \in \HONEzr(\Th) \cap
  H^s(\Th)$ with $s > 1$ be the solution of the continuous problem and
  $\ush\in V^h_{\kv}$ be the solution of the discrete
  problem. Consider the functions
  \begin{equation*}
    e_h=\ush-\calI_hu, \ e_{\calI}=u-\calI_hu, \ e_{\pi}=u-u_{\pi}, \ e_u=u_{\pi}-\calI_hu,
  \end{equation*}
  where $u_{\pi} \in \PS{k_o}(\Th)$ is the piecewise polynomial
  approximation of $\us$ defined in Bramble-Hilbert Lemma. Then it
  holds that
  \begin{align}
    \snorm{u-\ush}{1,\Omega}^2 + \alpha \ a_h(e_h,e_h) \lesssim 
    &\alpha^2 \sum_{\P\in \Th} h^2_\P \norm{f-f_h}{0,\P}^2 + 
    \alpha^2 \snorm{e_\pi}{1,\Th}^2 + \\
    &\alpha \snorm{e_\calI}{1,\Omega}^2 + 
    \alpha \sum_{\P \in \Th}\sigma^\P
  \end{align}
  where $\alpha$ is the coercivity constant and $\sigma^\P := \SPh
  ((I-\Pi^{\nabla,\P}_{k_0})e_u, (I-\Pi^{\nabla,\P}_{k_0})e_u)$.
\end{theo}
\begin{theo}[$\HONE$ error estimate with \textit{dofi-dofi} stabilization] 
  Assuming~\textbf{G\ref{g1},~G\ref{g4}}, let $\us\in \HONEzr(\Th)$ be
  the solution of the continuous problem and $\ush\in V^h_{\kv}$ be
  the solution of the discrete problem obtained with the
  \textit{dofi-dofi} stabilization.
  Assume moreover that $\us\in\HS{\bar{k}}(\Th)$ with
  $\bar{\ks}=\max\{\ks_o+1,\,\ks_{\partial}+2\}$ and
  $\fs\in\HS{\ks_o-1}$.
  Then it holds that
  \begin{equation}
    \snorm{\us-\ush}{1,\Omega}^2
    \lesssim\alpha\sum_{\P\in\Th}\left(
    (\alpha + \calN_{\P})^{1/2}\hP^{k_o}+\hh_{\bP}^{\ks_\partial}
    \right)^2,
  \end{equation}
  where $\hh_{\bP}$ denotes the maximum edge length, $\alpha$
  is the constant defined in \eqref{eq:error:small}, and $\calN_{\P}$
  is the number of edges in $\P$.
\end{theo}
\begin{theo}[$\HONE$ error estimate with \textit{trace} stabilization] 
  Under Assumption~\textbf{G\ref{g1}}, let $\us\in \HONEzr(\Th)$ be the
  solution of the continuous problem and $\ush\in V^h_{\kv}$ be the
  solution of the discrete problem obtained with the \textit{trace}
  stabilization.
  Assume moreover that $\us\in\HS{\bar{\ks}}(\Th)$ with
  $\bar{\ks}=\max\{\ks_o+1,\,\ks_{\partial}+2\}$ and
  $\fs\in\HS{\ks_o-1}$.
  Then it holds that
  \begin{equation}
    \snorm{\us-\ush}{1,\Omega}^2
    \lesssim\sum_{\P\in\Th}\left( \hP^{\ks_o} + \hh_{\bP}^{\ks_{\partial}}\right)^2.
  \end{equation}
\end{theo}
\section{Violating the geometrical assumptions}
\label{sec:violating}
We are here interested in testing the behaviour of the VEM on a number
of mesh ``datasets'' which systematically stress or violate the
geometrical assumptions from Section~\ref{subsec:star:assumptions}.
This enhances a correlation analysis between such assumptions and the
VEM performance, and we experimentally show how the VEM presents a
good convergence rate on most examples and only fails in very few
situations.


\subsection{Datasets definition}
\label{subsec:violating:datasets}
We start with defining the concept of dataset over a domain $\Omega$,
that for us will be the unit square $[0,1]^2$.
\begin{definition}
  We call a \emph{dataset} a collection
  $\calD:=\{\Omega_n\}_{n=0,\ldots,N}$ of discretizations of the
  domain $\Omega$ such that
  \begin{enumerate}
  \item[(i)] the mesh size of $\Omega_{n+1}$ is smaller than the mesh
    size of $\Omega_{n}$ for every $n=0,\ldots,N-1$;
  \item[(ii)] meshes from the same dataset follow a common refinement
    pattern, so that they contain similar polygons organized in
    similar configurations.
  \end{enumerate}
\end{definition}
Each mesh $\Omega_n$ can be uniquely identified via its size as $\Th$,
therefore every dataset $\calD$ can be considered as a subset of a
mesh family: $\calD=\{\Th\}_{h\in\calH'}\subset\calT$, being $\calH'$
a finite subset of $\calH$.

\medskip
Together with the violation of the geometrical assumptions, we are
also interested in measuring the behaviour of the VEM when the sizes
of mesh elements and edges scale in a nonuniform way during the
refinement process.
To this purpose, for each mesh $\Omega_n\in\calD$ we define the
following \textit{scaling indicators} and study their trend for $n\to
N$:
\begin{align}  \label{eq:scaling:indicators}
  A_n=\frac{\max_{\P\in \Omega_n} |\P|}{\min_{\P\in \Omega_n} |\P|}
  \quad\textrm{and}\quad 
  e_n=\frac{\max_{e\in \Omega_n} |e|}{\min_{e\in \Omega_n} |e|}.
\end{align}

We designed six particular datasets in order to cover a wide range of
combinations of geometrical assumptions and scaling indicators, as
shown in Table~\ref{table:ratios}.
For each of them we describe how it is built, how $A_n$ and $e_n$
scale in the limit for $n\to\Ns$ and which assumptions it fulfills or
violates.\\
Note that none of the considered datasets (exception made for the
reference dataset $\Dtriangle$) fulfills any of the sets of
assumptions required by the convergence analysis found in the
literature, c.f. Section~\ref{subsec:star:error}.\\
Each dataset is built around (and often takes its name from) a
particular polygonal element, or elements configurations, contained in
it, which is meant to stress one or more assumptions or indicators.
The detailed construction algorithms, together with the explicit
computations of $A_n$ and $e_n$ for all datasets, can be found in
\cite[Appendix B]{sorgente2021role}, while the complete collection of
the dataset can be downloaded at
(\footnote{https://github.com/TommasoSorgente/vem-quality-dataset}).

\medskip
\PGRAPH{Reference dataset.}
The first dataset, $\Dtriangle$, serves as a reference to evaluate the
other datasets by comparing the respective performance of the VEM over
each one.
It contains only triangular meshes, built by inserting a number of
vertices in the domain through the \textit{Poisson Disk Sampling}
algorithm \cite{bridson2007fast}, and connecting these vertices in a
Delaunay triangulation through the \textit{Triangle} library
\cite{Shewchuk-Triangle}.
The refinement is obtained by increasing the number of vertices
generated by the Poisson algorithm and computing a new Delaunay
triangulation.
The meshes in this dataset perfectly satisfy all the geometrical
assumptions and the indicators $A_n$ and $e_n$ are almost constant.

\medskip
\PGRAPH{Hybrid datasets.}
Next, we define some \textit{hybrid} datasets, which owe their name to
the presence of both triangular and polygonal elements (meaning
elements with more than three edges).
A number of identical polygonal elements (called the \emph{initial
polygons}) is inserted in $\Omega$, and the rest of the domain is
tessellated by triangles with area smaller than the one of the initial
polygons.
These triangles are created through the \textit{Triangle} library,
with the possibility to add Steiner points~\cite{Shewchuk-Triangle}
and to split the edges of the initial polygons, when necessary, with
the insertion of new vertices.
The refinement process is iterative, with parameters regulating the
size, the shape and the number of initial polygons.

\begin{figure}[htbp]  
  \centering
  \includegraphics[width=.8\linewidth]{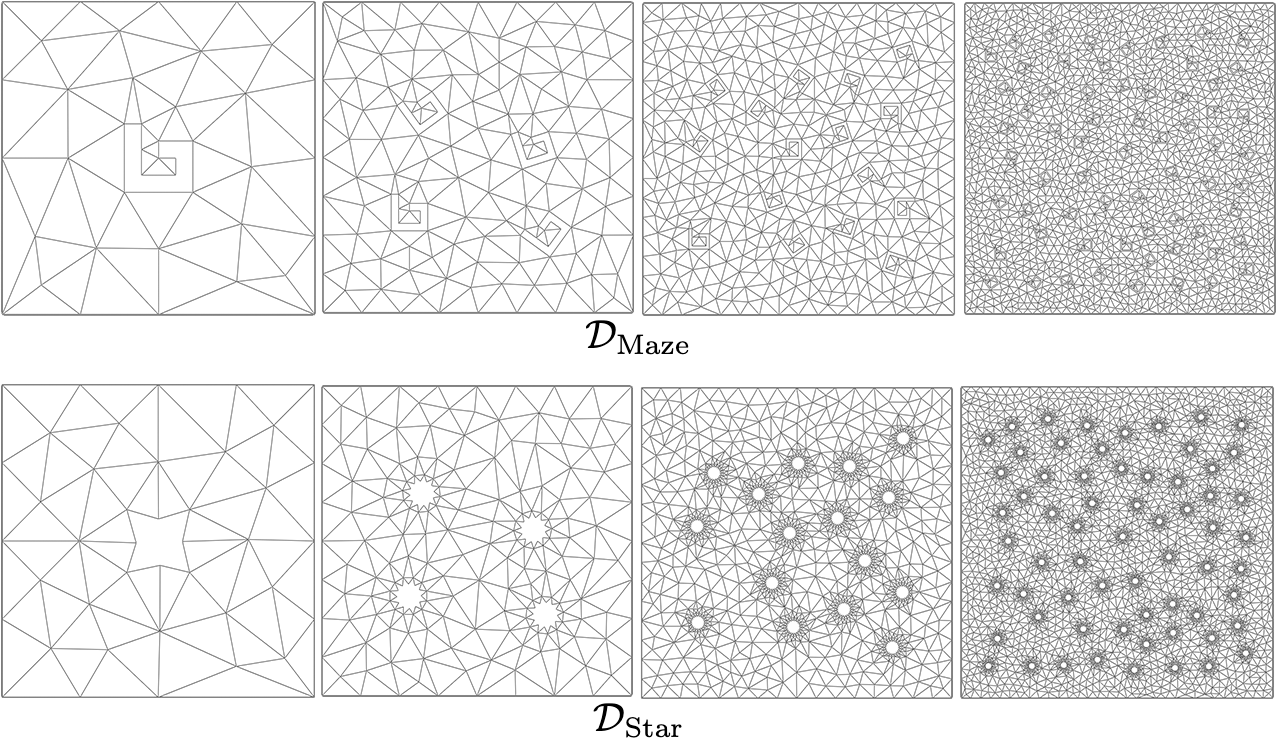}
  \caption{Meshes $\Omega_0, \Omega_2, \Omega_4, \Omega_6$ from datasets $\Dmaze$ and $\Dstar$.}
  \label{fig:hybrid-mesh}
\end{figure}

\medskip
We defined two hybrid datasets, $\Dmaze$ and $\Dstar$ shown in
Fig.~\ref{fig:hybrid-mesh}, as they violate different sets of
geometrical assumptions.
Other choices for the initial polygons are possible, for instance
considering the ones in Section~\ref{subsec:benchmark:results}.

\medskip
Dataset $\Dmaze$ is named after a $10$-sided polygonal element $\P$,
called ``maze'', spiralling around an external point.
Progressively, as $n\to N$, each mesh $\Omega_n$ in $\Dmaze$ contains
an increasing number of mazes $\P$ with decreasing thickness.
Concerning the scaling indicators, we have $A_n\sim a^n$ for a
constant $e<a<3$ and $e_n\sim n\log(n)$, and the geometrical
assumptions violated by this dataset are:

\begin{description}
\item[\textbf{G1}] because every $\P$ is obviously \textit{not} star-shaped;
\item[\textbf{G1-weak}] because of the decreasing thickness of the
  mazes.
  Indeed, it would be trivial to split $\P$ into a collection of
  star-shaped elements $\P_i$ (rectangles, for instance), but as $n\to
  N$ each radius $r_{\P_i}$ would decrease faster than the respective
  diameter $\hh_{\P_i}$, unless considering an infinite number of
  them.
  Therefore, it would not be possible to bound the ratio
  $r_{\P_i}/\hh_{\P_i}$ from below with a global $\rho$ independent of
  $h$;
\item[\textbf{G2, G2-strong}] because the length of the shortest edge
  $\E$ of $\P$ decreases faster than the diameter $\hP$.
\end{description}

\begin{figure}[htbp]  
  \centering
  \includegraphics[width=.5\linewidth]{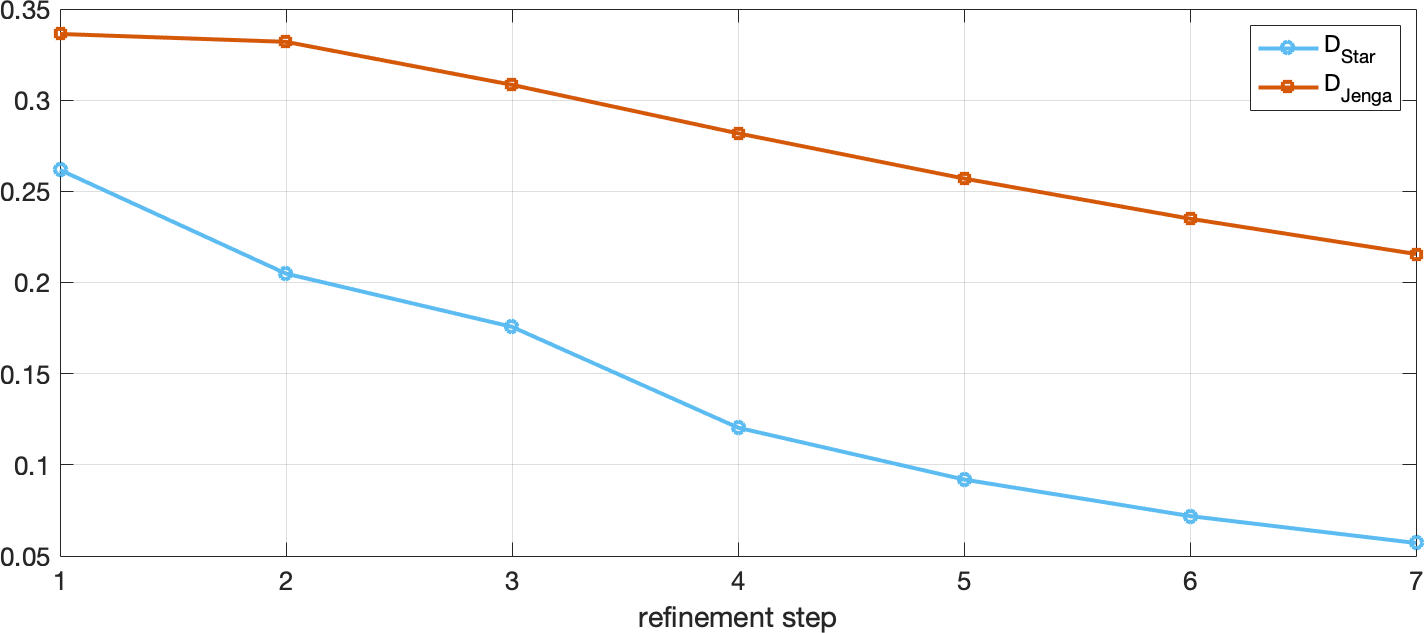}
  \caption{Evolution of the ratio $\rP/\hP$ for the meshes in datasets $\Dstar$ and $\Djenga$.}
  \label{fig:ratio}
\end{figure}

\medskip
Dataset $\Dstar$ is built by inserting star-like polygonal elements,
still denoted by $\P$.
As $n\to N$, the number of spikes in each $\P$ increases and the inner
vertices are moved towards the barycenter of the element.
Both indicators $A_n$ and $e_n$ scale linearly, and the dataset
violates the assumptions:
\begin{description}
\item[\textbf{G1}] each star $\P$ is star-shaped with respect to the
  maximum circle inscribed in it, but as shown in
  Fig.~\ref{fig:ratio}, the radius $\rP$ of such circle decreases
  faster than the elemental diameter $\hP$, therefore it is not
  possible uniformly bound from below the quantity $\rP/\hP$ with a
  global $\rho$;
\item[\textbf{G1-weak}] in order to satisfy it, we should split each
  $\P$ into a number of sub-polygons, each of them fulfilling
  \textbf{G\ref{g1}}.
  Independently of the way we partition $\P$, the number of sub-polygons
  would always be bigger than or equal to the number of spikes in $\P$,
  which is constantly increasing, hence the number of sub-polygons would
  tend to infinity;
\item[\textbf{G2-strong}] because the distance between the inner
  vertices of $\P$ decreases faster than $\hP$ (but \textbf{G\ref{g2}}
  holds, because the edges scale proportionally to $\hP$);
\item[\textbf{G3}] because the number of spikes in each $\P$ increases
  from mesh to mesh, therefore the total number of vertices and edges
  in a single element cannot be bounded uniformly.
\end{description}

\medskip
\PGRAPH{Mirroring datasets.}
As an alternative strategy to build a sequence of meshes whose
elements are progressively smaller, we adopt an iterative
\textit{mirroring} technique: given a mesh $\calM$ defined on
$\Omega$, we generate a new mesh $\calM'$ containing four adjacent
copies of $\calM$, opportunely scaled to fit $\Omega$.\\
The starting point for the construction of a dataset is the first base
mesh $\widehat{\Omega}_0$, which coincides with the first
computational mesh $\Omega_0$.
At every step $1\le n\le N$, a new base mesh $\widehat{\Omega}_n$ is
built from the previous base mesh $\widehat{\Omega}_{n-1}$.
The computational mesh $\Omega_{n}$ is then obtained by applying the
mirroring technique $4^n$ times to the base mesh
$\widehat{\Omega}_{n}$.
This construction allows us to obtain a number of vertices and degrees
of freedom in each mesh that is comparable to that of the meshes at
the same refinement level in datasets $\Dmaze$ and $\Dstar$.

\begin{figure}[htbp]  
  \centering
  \includegraphics[width=.8\linewidth]{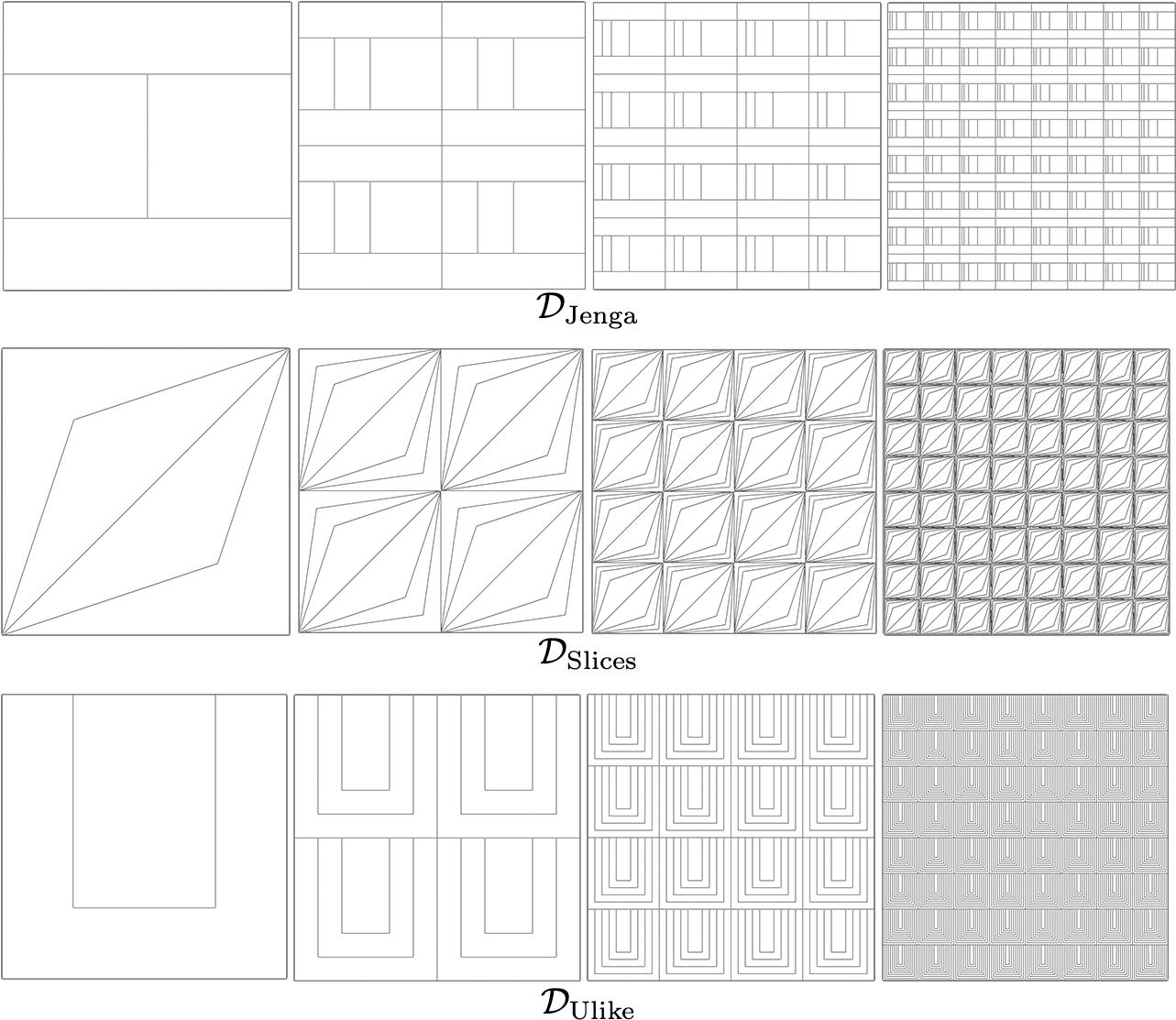}
  \caption{Meshes $\Omega_0, \Omega_1, \Omega_2, \Omega_3$ from datasets $\Djenga$, $\Dslices$ and $\Dulike$.}
  \label{fig:nonuniform-mesh}
\end{figure}

\medskip
The $n$-th base mesh $\widehat{\Omega}_{n}$ of dataset $\Djenga$
(Fig.~\ref{fig:nonuniform-mesh}, top) is built as follows.
We start by subdividing the domain $(0,1)^2$ into three horizontal
rectangles with area equal to $1/4$, $1/2$ and $1/4$ respectively.
Then, we split the rectangle with area $1/2$ vertically, into two
identical rectangles with area $1/4$.
This concludes the construction of the base mesh $\widehat{\Omega}_0$,
which coincides with mesh $\Omega_0$.
At each next refinement step $n\geq1$, we split the left-most
rectangle in the middle of the base mesh $\widehat{\Omega}_{n-1}$ by
adding a new vertical edge, and apply the mirroring technique to
obtain $\Omega_{n}$.
Despite being made entirely by simple rectangular elements, this mesh
family is the most complex one: both $A_n$ and $e_n$ scale like
$2^{n}$ and it breaks all the assumptions.
\begin{description}
\item[\textbf{G2, G2-strong}] because the ratio $|\E|/h_\P$ decreases
  unboundedly in the left rectangle $\P$, as shown in
  Fig.~\ref{fig:ratio}.
  This implies that a lower bound with a uniform constant $\rho$
  independent of $\hh$ cannot exist;
\item[\textbf{G1, G1-weak}] since the length of the radius $\rP$ of
  the maximum disc inscribed into a rectangle is equal to $1/2$ of its
  shortest edge $\E$, the ratio $\rP/\hP$ also decreases;
\item[\textbf{G3}] because the number of edges of the top (or the
  bottom) rectangular element grows unboundedly;
\item[\textbf{G4}] because the one-dimensional mesh built on the
  boundary of the top rectangular element cannot be subdivided into a
  finite number of quasi uniform sub-meshes.
  In fact, either we have infinite sub-meshes or an infinite edge ratio.
\end{description}

\medskip
For the dataset $\Dslices$ (Fig.~\ref{fig:nonuniform-mesh}, middle),
the $n$-th base mesh $\widehat{\Omega}_{n}$ is built as follows.
First, we sample a set of points along the diagonal (the one
connecting the vertices $(0,1)$ and $(1,0)$) of the reference square
$[0,1]^2$, and connect them to the vertices $(0,0)$ and $(1,1)$.
In particular, at each step $n\geq0$, the base mesh
$\widehat{\Omega}_{n}$ contains the vertices $(0,0)$ and $(1,1)$, plus
the vertices with coordinates $(2^{-i},1-2^{-i})$ and
$(1-2^{-i},2^{-i})$ for $i=1,\ldots,n+2$.
Then, we apply the mirroring technique.
Since no edge is ever split, we find that $e_n\sim c$, while $A_n\sim
2^{n}$.
The dataset $\Dslices$ violates assumptions
\begin{description}
\item[\textbf{G1, G1-weak}] because, up to a multiplicative scaling
  factor depending on $\hh$, the radius $\rP$ is decreasing faster
  than the diameter $\hP$, which is constantly equal to $\sqrt{2}$
  times the same scaling factor.
  Moreover, any finite subdivisions of the mesh elements would suffer
  the same issue;
\item[\textbf{G2-strong}] the vertices sampled along the diagonal have
  accumulation points at $(0,1)$ and $(1,0)$, therefore the distance
  between two vertices decreases progressively.
\end{description}

\medskip
In $\Dulike$ (Fig.~\ref{fig:nonuniform-mesh}, bottom), we build the
mesh $\widehat{\Omega}_{n}$ at each step $n\geq0$ by inserting $2^n$
equispaced $U$-shaped continuous polylines inside the domain, creating
as many $U$-like polygons.
Then, we apply the mirroring technique as usual.
Edge lengths scale exponentially and areas scale uniformly, i.e., $e_n
\sim 2^{n}$, $A_n\sim c$, and the violated assumptions are:
\begin{description}
\item[\textbf{G1, G1-weak, G2, G2-strong}] for arguments similar to
  the ones seen for $\Dmaze$;
\item[\textbf{G3}] in order to preserve the connectivity, the lower
  edge of the more external $U$-shaped polygon in every base mesh must
  be split into smaller edges before applying the mirroring technique.
  Hence the number of edges of such elements cannot be bounded from
  above.
\end{description}

\medskip
\PGRAPH{Multiple mirroring datasets.}
As a final test, we modified datasets $\Djenga$, $\Dslices$ and
$\Dulike$ in order to particularly stress the scaling indicators $A_n$
and $e_n$.\\
We obtained this by simply inserting four new elements at each step
instead of one.
The resulting datasets, $\DjengaM$, $\DslicesM$ and $\DulikeM$, are
qualitatively similar to the correspondent mirroring datasets above.
Each of these datasets respects the same assumptions of its original
version, but the number of elements at each refinement step now
increases four times faster.
Consequently, the indicators $A_n$ and $e_n$ change from $2^n$ to
$2^{4n}$, but $A_n$ remains constant for $\DulikeM$, and $e_n$ remains
constant for $\DslicesM$.
\begin{table}[t]  
  \centering
  \caption{Summary of the geometrical conditions violated and the
    asymptotic trend of the indicators $A_n$ and $e_n$ for all
    datasets.  \textbf{G1w} and \textbf{G2s} stand for
    \textbf{G1-weak} and \textbf{G2-strong}, respectively, and $a$ is
    a constant such that $e<a<3$).}
  \vspace{0.5cm}
  \begin{tabular}{ccccccc}
    \hline\noalign{\smallskip}
    \textbf{dataset \ } & $\Dtriangle$ & $\Dmaze$ & $\Dstar$ & $\Djenga$ & $\Dslices$ & $\Dulike$\\
    \noalign{\smallskip}\hline\noalign{\smallskip}
    \textbf{G1, G1w, G2s} & 
    & $\times$ & $\times$ & $\times$ & $\times$ & $\times$ \\
    \textbf{G2} & 
    & $\times$ & & $\times$ & & $\times$ \\
    \textbf{G3} & 
    & & $\times$ & $\times$ & & $\times$ \\
    \textbf{G4} & 
    & & & $\times$ & & \\
    $A_n$ & 
    $c$ & $a^n$ & $n$ & $2^n$ & $2^n$ & $c$ \\
    $e_n$ & 
    $c$ & $n\log(n)$ & $n$ & $2^n$ & $c$ & $2^n$\\
    \noalign{\smallskip}\hline
  \end{tabular}
  \label{table:ratios}
\end{table}


\subsection{VEM performance over the datasets}
\label{subsec:violating:performance}
We solved the discrete Poisson problem \eqref{eq:Poisson:weak} with
the VEM \eqref{eq:VEM} described in Section~\ref{sec:VEM} for
$k=1,2,3$ over each mesh of each of the datasets defined in
Section~\ref{subsec:violating:datasets}, using as groundtruth the
function
\begin{equation}
  u(x,y)=\frac{\sin(\pi x)\sin(\pi y)}{2\pi^2}, \hspace{0.3cm} (x,y)\in\Omega=(0,1)^2.
\end{equation}
This function has homogeneous Dirichlet boundary conditions, and this
choice was appositely made to prevent the boundary treatment from
having an influence on the approximation error.

\medskip
\PGRAPH{Performance indicators}
After computing the solution of the problem on a particular dataset,
we want to measure the performance of the method over that dataset,
that is, the accuracy and the convergence rate.
We selected, among many possible alternatives, some quantities which
can indicate if the error between the continuous solution $\us$ and
the computed solution $\ush$ is small and if the VEM worked properly
in the computation of $\ush$.
A more complete and accurate analysis of the possible performance
indicators can be found in Section~\ref{subsec:benchmark:performance}.

\begin{figure}[htbp] 
  \centering
  \includegraphics[width=\linewidth]{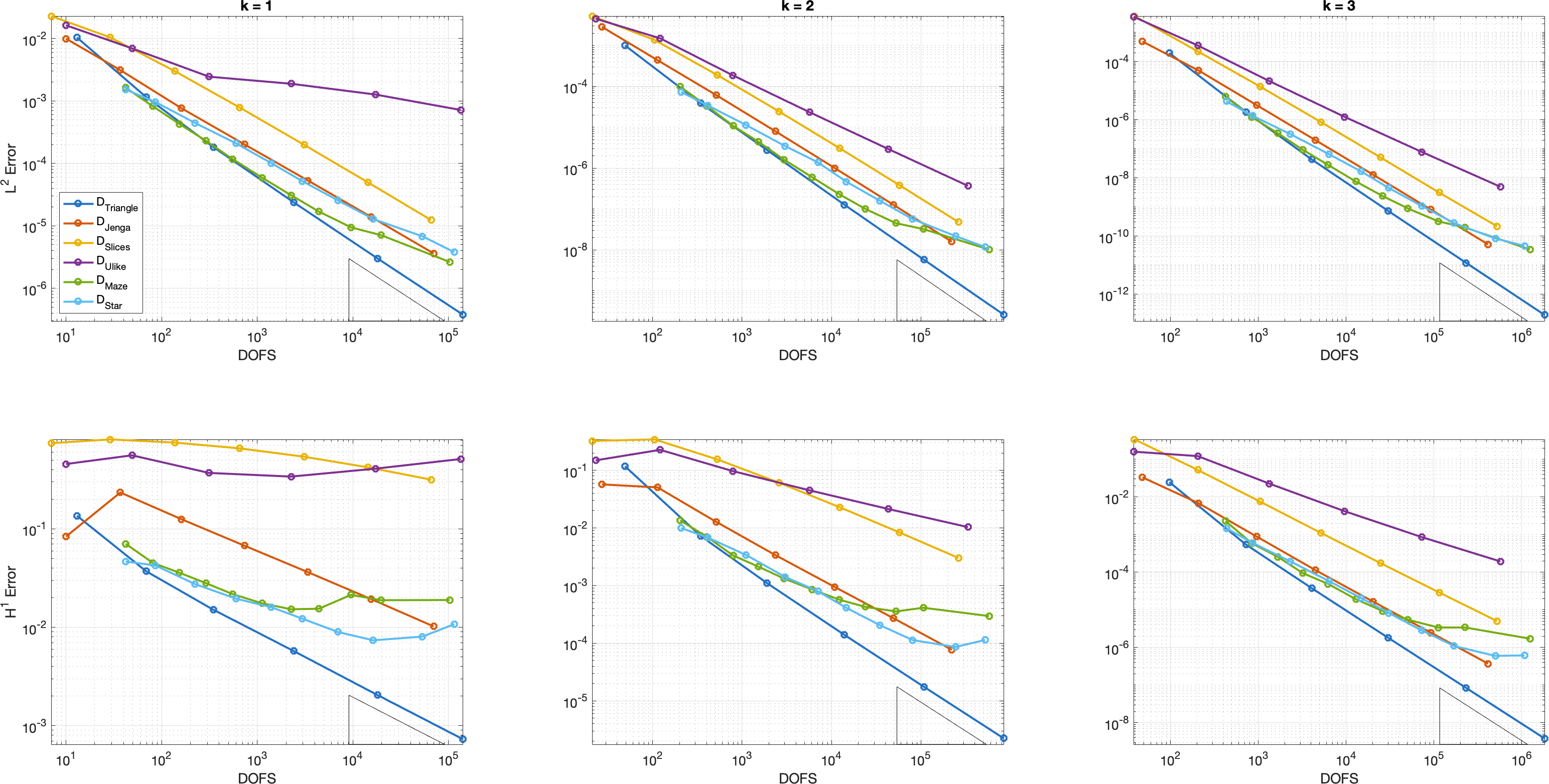}
  \caption{$L^2$-norm and $H^1$-seminorm of the approximation errors
    of the reference, hybrid and mirroring datasets for $k=1,2,3$.}
  \label{fig:performance1}
\end{figure}

\begin{figure}[htbp]  
  \centering
  \includegraphics[width=\linewidth]{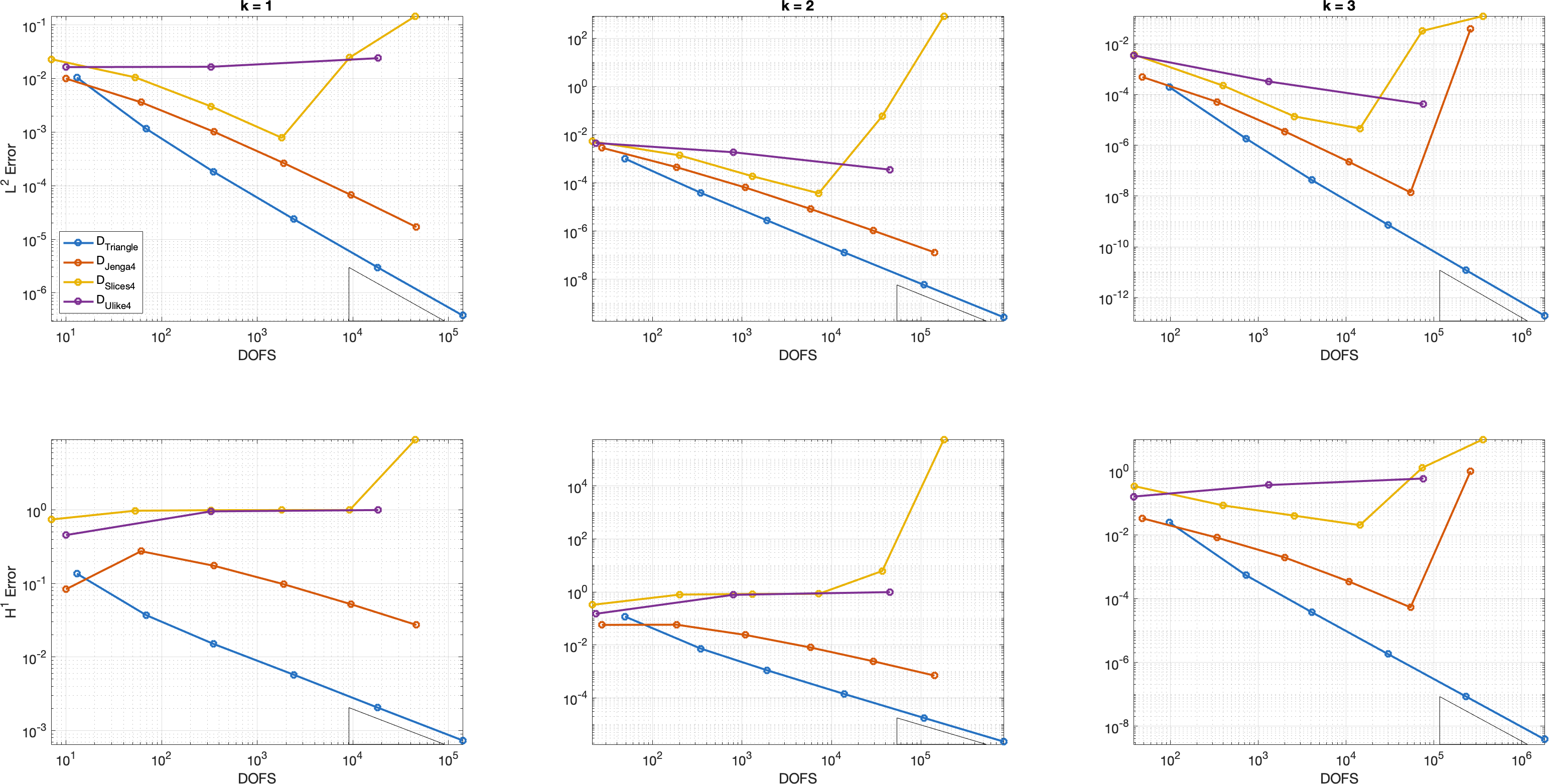}
  \caption{$L^2$-norm and $H^1$-seminorm of the approximation errors
    of the reference and multiple mirroring datasets for $k=1,2,3$. }
  \label{fig:performance4}
\end{figure}

\smallskip
The approximation error might be measured in different norms: the most
widely used in our framework are the relative $\HONE$-seminorm and the
relative $\LTWO$-norm (in the following we use the generic term
\textit{norms} to indicate both of them):
\begin{align*}
  \frac{\snorm{\us-\ush}{1,\Omega}}{\snorm{\us}{1,\Omega}}
  \quad\text{and}\quad
  \frac{\norm{\us-\ush}{0,\Omega}}{\norm{\us}{0,\Omega}}.
\end{align*}
In Fig.~\ref{fig:performance1} and Fig.~\ref{fig:performance4} we plot
the two norms of the error generated by the VEM on each dataset as the
number of DOFs increases (that is, as $n\to N$).

\smallskip
Another quantity, not directly related to the error, which can be of
interest is the condition number of the matrices \textbf{G} and
\textbf{H} (with the notation adopted in
\cite{BeiraodaVeiga-Brezzi-Marini-Russo:2014}):
\begin{align*}
  \text{cond}(\textbf{G})=\|\textbf{G}\|\|\textbf{G}^{-1}\|,\qquad 
  \text{cond}(\textbf{H})=\|\textbf{H}\|\|\textbf{H}^{-1}\|
\end{align*}
Matrix \textbf{G} is involved in the computation of the \textit{local
  stiffness matrix} and the projector $\Pin{k}$, while matrix
\textbf{H} is involved in the computation of the \textit{local mass
  matrix} and the projector $\Piz{k}$.

\smallskip
Last, as an estimate of the error produced by projectors $\Pin{k}$ and
$\Piz{k}$, represented by matrices $\matPin{k}$ and $\matPiz{k}$, we
check the identities
\begin{align}\label{eq:projector:identities}
  |\matPin{k}\textbf{D}-\textbf{I}|=0 \quad\text{and}\quad |\matPiz{k}\textbf{D}-\textbf{I}|=0,
\end{align}
where for $k<3$ we have $\matPiz{k}:=\matPin{k}$.
The computation of the projectors is obviously affected by the
condition numbers of \textbf{G} and \textbf{H}, but the two indicators
are not necessarily related.\\
Condition numbers and identity values for $k=1,2,3$ are reported in
Table~\ref{table:numerical_performance}.
All of these quantities are computed element-wise and the maximum
value among all elements of the mesh is selected.

\medskip
\PGRAPH{Performance}
The VEM performs perfectly over the reference dataset $\Dtriangle$,
and this guarantees for the correctness of the method.
The approximation error evolves in accordance to the theoretical
results (the slopes being indicated by the triangles) both in $L^2$
and in $H^1$ norms for all $k$ values, and condition numbers and
optimal errors on the projectors $\Piz{k}$ and $\Pin{k}$ remain
optimal.

\smallskip
For the hybrid datasets $\Dstar$ and $\Dmaze$, errors decrease at the
correct rate for most of the meshes, and only start deflecting for
very complicated meshes with very high numbers of DOFs.
These deflections are probably due to the extreme geometry of the star
and maze polygons and not to numerical problems, as in both datasets
we have cond(\textbf{G}) $<10^6$ and cond(\textbf{H}) $<10^9$, which
are still reasonable values.
Projectors seem to work properly: $|\matPin{k}\textbf{D}-\textbf{I}|$
remains below $10^{-8}$ and $|\matPiz{k}\textbf{D}-\textbf{I}|$ below
$10^{-7}$.
In a preliminary stage of this work, we obtained similar plots (not
reported here) using other hybrid datasets built in the same way, with
polygons surrounded by triangles.
In particular, we did not notice big differences when constructing
hybrid datasets as in Section~\ref{subsec:violating:datasets} with any
of the initial polygons of Section~\ref{subsec:benchmark:results}.

\smallskip
On the meshes from mirroring datasets $A_n$ or $e_n$ may scale
exponentially instead of uniformly, as reported in
Table~\ref{table:ratios}.
This reflects to cond(\textbf{G}) and cond(\textbf{H}), which grow up
to, respectively, $10^{10}$ and $10^{14}$ for $\Djenga$ in the case
$k=3$.
Nonetheless, the discrepancy of the projectors identities
\eqref{eq:projector:identities} remains below $10^{-5}$, which is not
far from the results obtained with datasets $\Dmaze$ and $\Dstar$.
The method exhibits an almost perfect convergence rate on dataset
$\Djenga$, even though $L^2$ and $H^1$ errors are bigger in magnitude
than the ones measured for hybrid datasets.
$\Dslices$ produces even bigger errors and a non-optimal convergence
rate, and $\Dulike$ is the dataset with the poorest performance, but
the VEM still converges at a decent rate for $k>1$.\\
This may be due to the fact that for $k=1$ the DOFs correspond to the
vertices of the mesh, which are disposed in a particular configuration
that generates horizontal bands in the domain completely free of
vertices, and therefore of data.
For $k>1$ instead, we have DOFs also on the edges and inside the
elements, hence the information is more uniformly distributed.

\smallskip
In the case of multiple mirroring datasets the method diverges badly
on all datasets (see Fig.~\ref{fig:performance4}), and this is
principally due to the very poor condition numbers of the matrices
involved in the calculations (see
Table~\ref{table:numerical_performance}).
The plots relative to datasets $\DjengaM$ and $\DslicesM$ maintain a
similar trend to those of $\Djenga$ and $\Dslices$ in
Fig.~\ref{fig:performance1}, until numerical problems cause
cond(\textbf{G}) and cond(\textbf{H}) to explode up to over $10^{30}$
for $\DjengaM$ and $10^{18}$ for $\DslicesM$.
In such conditions, the projection matrices $\matPin{k}$ and
$\matPiz{k}$ become meaningless and the method diverges.
The situation slightly improves for $\DulikeM$: cond(\textbf{H}) is
still $10^{16}$, but the discrepancy of $\matPin{k}$ and $\matPiz{k}$
remain acceptable.
As a result, the method on $\DulikeM$ does not properly explode, but
the approximation error and the convergence rate are much worse than
those seen for $\Dulike$ in Fig.~\ref{fig:performance1}.
\begin{table}[H]  
  \centering
  \caption{Summary of numerical performance for all datasets.
  We report the $\log_{10}$ of the original values for the condition
  number of \textbf{G} and \textbf{H} and the discrepancy of
  projection matrices $\matPin{k}$ and $\matPiz{k}$.
  Note that for $k<3$ we have $\matPiz{k}:=\matPin{k}$.}
  \resizebox{\columnwidth}{!}{
  \begin{tabular}{c|ccc|ccc|ccc|ccc|ccc|ccc|ccc|ccc|ccc}
    \hline
    \multicolumn{1}{c}{\textbf{dataset}} &
    \multicolumn{3}{c}{$\Dtriangle$} & \multicolumn{3}{c}{$\Dmaze$} &  \multicolumn{3}{c}{$\Dstar$} & \multicolumn{3}{c}{$\Djenga$} &  \multicolumn{3}{c}{$\Dslices$} & \multicolumn{3}{c}{$\Dulike$} & \multicolumn{3}{c}{$\DjengaM$} &  \multicolumn{3}{c}{$\DslicesM$} &  \multicolumn{3}{c}{$\DulikeM$}\\
    \hline\noalign{\smallskip}
    $k$ & 1 & 2 & 3 & 1 & 2 & 3 & 1 & 2 & 3 & 1 & 2 & 3 & 1 & 2 & 3 & 1 & 2 & 3 & 1 & 2 & 3 & 1 & 2 & 3 & 1 & 2 & 3\\
    \noalign{\smallskip}\hline\noalign{\smallskip}
    cond(\textbf{G}) 
    & 0 & 2 & 5 
    & 2 & 3 & 6 
    & 1 & 3 & 6 
    & 1 & 5 & 10 
    & 2 & 4 & 6 
    & 1 & 4 & 7 
    & 6 & 18 & 31 
    & 6 & 8 & 10 
    & 2 & 6 & 11\\ 
    cond(\textbf{H})
    & 2 & 5 & 7 
    & 2 & 5 & 8 
    & 3 & 6 & 9 
    & 4 & 9 & 14 
    & 2 & 8 & 10 
    & 3 & 7 & 10 
    & 13 & 26 & 39 
    & 2 & 15 & 18 
    & 5 & 10 & 16\\ 
    $|\matPin{k}\textbf{D}-\textbf{I}|$
    & -13 & -11 & -9 
    & -12 & -10 & -8 
    & -12 & -10 & -8 
    & -12 & -8 & -5 
    & -12 & -10 & -9 
    & -13 & -10 & -8 
    & -9 & 3 & 13 
    & -8 & -6 & -5 
    & -13 & -8 & -5\\ 
    $|\matPiz{k}\textbf{D}-\textbf{I}|$
    &  &  & -10 
    &  &  & -8 
    &  &  & -7 
    &  &  & -5 
    &  &  & -5 
    &  &  & -7 
    &  &  & 20 
    &  &  & 8 
    &  &  & -4\\ 
    \noalign{\smallskip}\hline
  \end{tabular}
  }
  \label{table:numerical_performance}
\end{table}

\smallskip
As a preliminary conclusion, by simply looking at the previous plots
we observe that the relationship between the geometrical assumptions
respected by a certain dataset and the performance of the VEM on it is
not particularly strong.
In fact, we obtained reasonable errors and convergence results on
datasets violating several assumptions (all of them, in the case of
$\Djenga$).
\section{Mesh quality metrics}
\label{sec:benchmark}

The aim of this section is to introduce some geometrical parameters,
that we will refer to as {\em quality metrics}, which are potentially
well suited to measure the shape regularity of a polygon, and study,
statistically, the behavior of a VEM solver, as such measures
degrade. In the following we present a list of polygon quality metrics
and different strategies to combine them to form a quality metric for
a polygonal tessellation. We will also introduce a list of parameters
providing us with different ways of measuring the performance of the
VEM solver at hand.


\subsection{Polygon quality metrics}
\label{subsec:benchmark:metrics}

Different parameters provide us with some information on how much an
element is far from being ``nice'', and the different assumptions
presented in Section \ref{subsec:star:assumptions} give us a first
list of possibilities. In the following we present the list of polygon
quality metrics that we singled out for our study.

\begin{itemize}

\item \emph{Circumscribed circle radius (\CC)}: this is defined as the
  radius of the smallest circle fully containing~$P$.  The parameter
  \CC~ is computed by treating the vertices of~$P$ as a point cloud
  and running the Welzl's algorithm to solve the minimum covering
  circle problem~\cite{welzl1991smallest}. We point out that this
  choice does not yield the classical definition of circumscribed
  circle, requiring that all the vertices of the polygon lie on the
  circle, that does not necessarily exist for all polygons
  (Fig.~\ref{fig:cc}).
	
  \begin{figure}[htbp]
    \centering
    \includegraphics[width=.5\columnwidth]{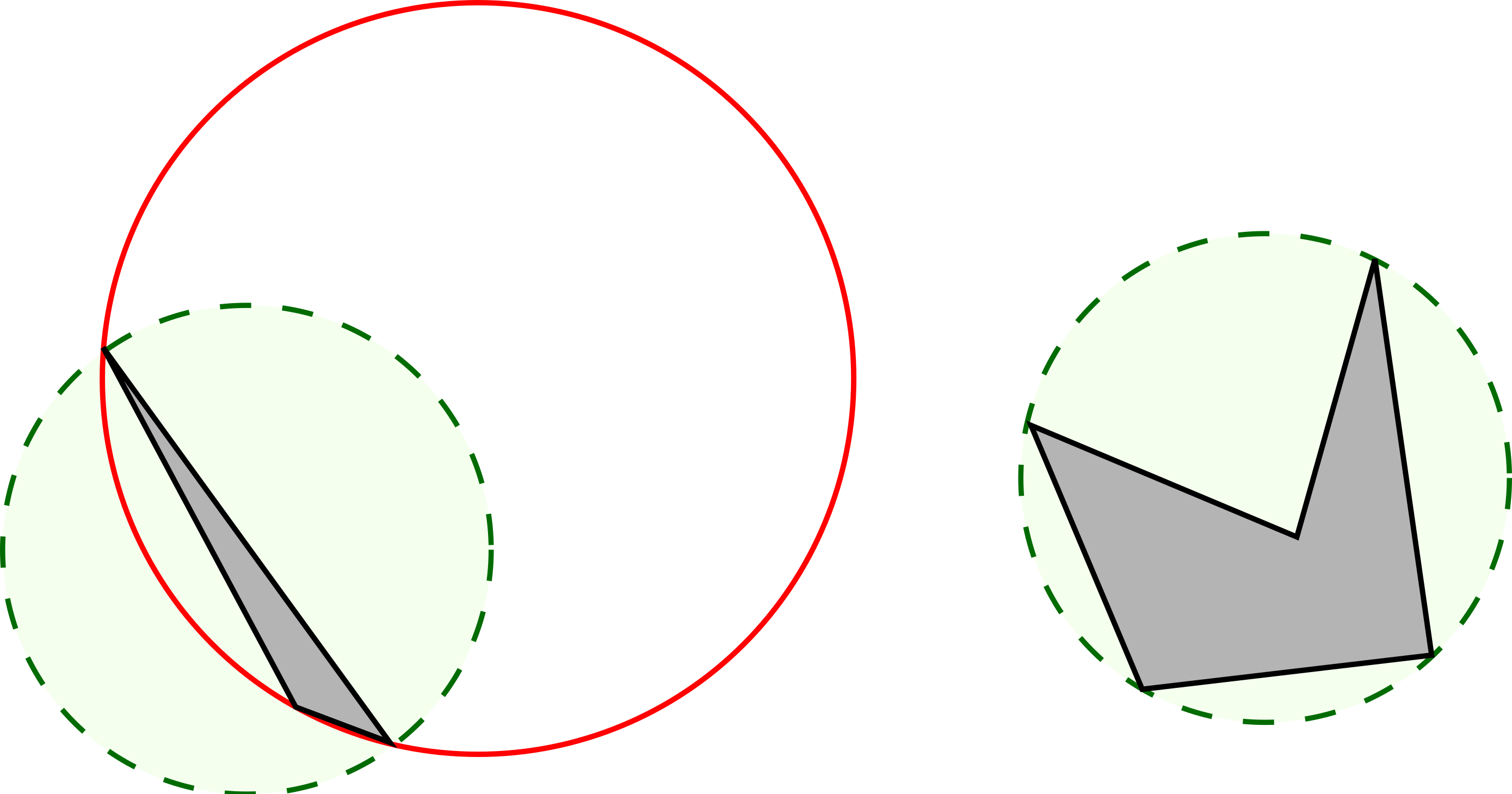}
    \caption{To be able to scale on general polygons, we define the
      radius of the circumscribed circle (CC) as the radius of the
      smallest circle containing the polygon itself. (Left) For a
      skinny triangle, CC is the radius of the smallest circle that
      passes through the endpoints of its longest edge (green), and
      not of the circle passing through all its three vertices
      (red). (Right) A general polygon and its CC.}
     \label{fig:cc}
  \end{figure}
  
\item \emph{Inscribed circle radius (\IC)}: this is defined as the
  radius of the biggest circle fully contained in~$P$.  For the
  computation of \IC, starting from a Voronoi diagram of the edges
  of~$P$, we select the corner in the diagram that is furthest from
  all edges as center of the circle. The radius of the IC is the
  minimum distance between such point and any of the edges of~$P$.
  For the computation of the diagram of the set of edges, which,
  differently from the point case, has curved boundaries between
  cells, we apply the Boost Polygon
  Library~\cite{simonson2009geometry}.
	
\item \emph{Circle ratio (\CR)}: the ratio between \IC~ and
  \CC. Differently from the previous two, this measure does not depend
  on the scale of the polygon, and is always defined in the
  range~$(0,1]$.
	
  \item \emph{Area (\AR)}: the area of the polygon~$P$.
	
  \item \emph{Kernel area (\KE)}: the area of the kernel of the
    polygon, defined as the set of points~$p \in P$ from which the
    whole polygon is visible. If the polygon is convex, then the area
    of the polygon and the area of the kernel are equal. If the
    polygon is star-shaped, then the area of the kernel is a positive
    number. If the kernel is not star-shaped, then the kernel of the
    polygon is empty and KE will be zero.
	
  \item \emph{Kernel-area/area ratio (\KA)}: the ratio between the
    area of the kernel of~$P$ and its whole area. For convex polygons,
    this ratio is always 1. For concave star-shaped polygons, \KA~ is
    strictly defined in between 0 and 1. For non star-shaped polygons,
    \KA~ is always zero.
	
  \item \emph{Area/perimeter ratio (\AP)}, or \emph{compactness}:
    defined as
    \mbox{$\frac{2\pi*\textsf{area}(P)}{\textsf{perimeter}(P)^2}$}. This
    measure reaches its maximum for the most compact 2D shape (the
    circle), and becomes smaller for less compact polygons.
	
  \item \emph{Shortest edge (\SE)}: the length of the shortest edge
    of~$P$.
	  
  \item \emph{Edge ratio (\ER)}: the ratio between the length of the
    shortest and the longest edge of~$P$.
	
  \item \emph{Minimum vertex to vertex distance (\MP)}, which is the
    minimum distance between any two vertexes in~$P$. \MP~is always
    less then or equal to \SE. In case the two vertexes realizing the
    minimum distance are also extrema of a common edge, \MP~ and \SE~
    are equals.
	
	
  \item \emph{Minimum angle (\MA)}: the minimum inner angle of the
    polygon~$P$.
	
  \item \emph{Maximum angle (\MX)}: the maximum inner angle of the
    polygon~$P$.
	
  \item \emph{Number of edges ($\NS$)}: the number of edges of the
   polygon.
	
  \item \emph{Shape regularity (\SR)}: the ratio between the radius of
    the circle inscribed in the kernel of the element and the radius
    \CC~of the circumscribed circle
    (Fig.~\ref{fig:shape_regularity}). This assumes the value $0$ for
    polygons which are not star shaped.
	
    \begin{figure}[htbp]
      \centering
      \begin{tabular}{c}
        \includegraphics[height=110pt]{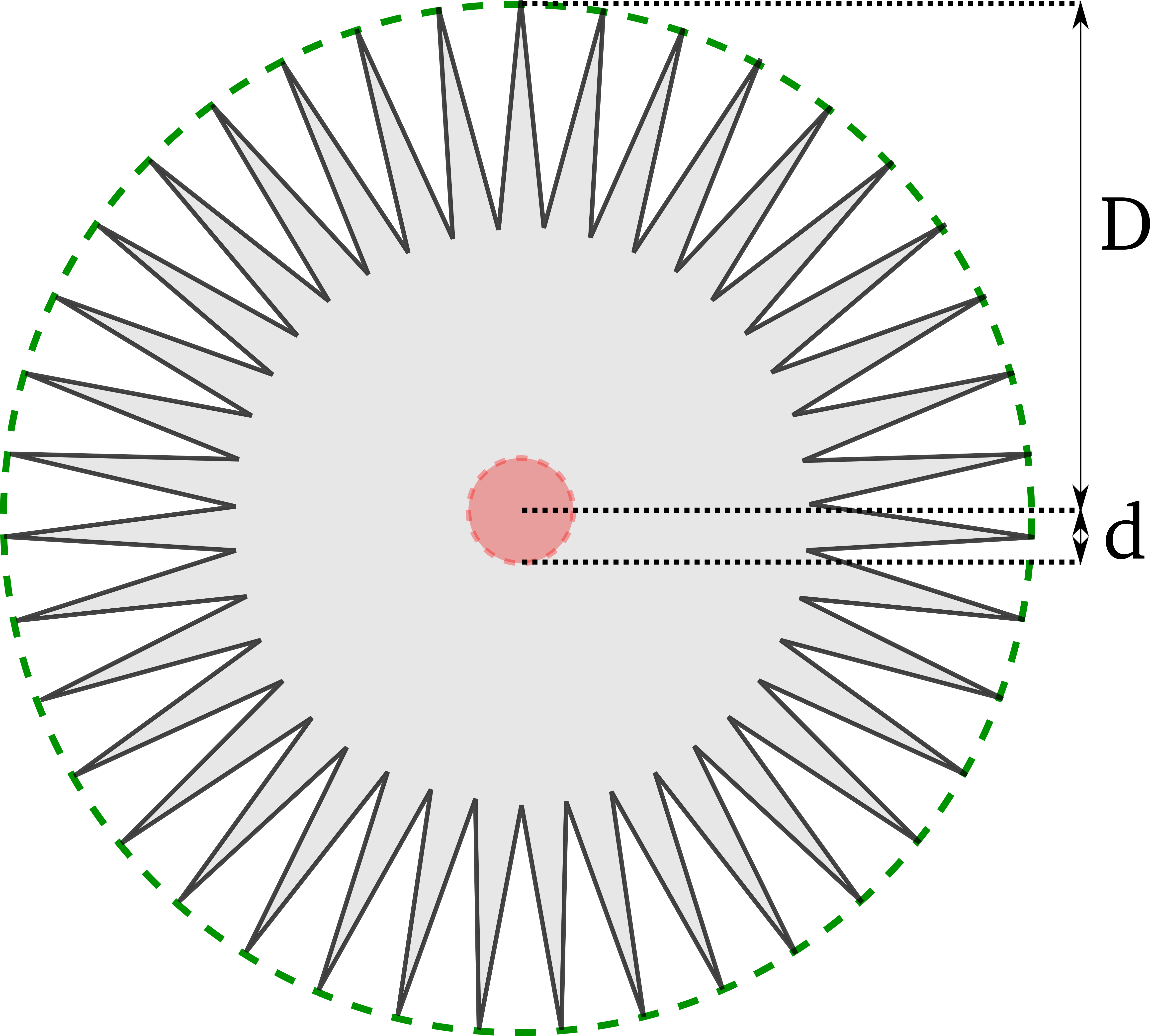}
      \end{tabular}
\caption{Polygon shape regularity expressed in terms of the ratio
  between the radius maximal ball inscribed in the kernel of the
  polygon ($d$), and the radius of the maximal ball inscribing the
  element ($D$).}
  \label{fig:shape_regularity}
    \end{figure}
    
\end{itemize}

Remark that, while in triangular elements all the previous metrics are
strictly bound to each other, for general polygons this is not the
case, and, given a couple of such metrics, it is in general possible
to find sequences of polygons for which one of the two degenerates,
while the other stays constant. Remark also that some of the metrics,
namely \CR, \KA, \AP, \ER, \MA, \MX, $\NS$~and \SR, are scale
invariant, and only depend on the shape of the polygon and not on its
size.

\newcommand{\metric}{{\mathfrak m}}

\PGRAPH{Aggregating polygon quality metrics into mesh quality metrics}
Let now consider a tessellation $\Th=\{P_i, i=1,\cdots,N\}$ with $N$
elements. For $\metric$ being any of the above polygon quality metrics
we can consider the vector
\[
\boldsymbol{\metric}(\Th) = (\metric(P_1),\cdots,\metric{(P_N)}).
\] 
Starting from such a vector we can define a measure of the quality of
the whole mesh, by following different strategies. More precisely, we
considered the following possibilities.
\begin{itemize}
\item \emph{Average}
  \[
  \metric_{\text{av}}(\Th) = \text{average}(\boldsymbol{\metric}(\Th)) = \frac 1 N \sum_{i=1}^N \metric(P_i).
  \]
\item \emph{Euclidean norm}
  \[
  \metric_{\text{l}_2}(\Th) = \| \boldsymbol{\metric}(\Th)\|_2 = \left( \sum_{i=1}^N \metric(P_i)^2 \right)^{1/2}
  \]
\item \emph{Maximum}
  \[
  \metric_{\text{max}}(\Th) =  \max \{ \metric(P_i), i = 1,\cdot,N \}
  \]
\item \emph{Minimum}
  \[
  \metric_{\text{min}}(\Th) =  \min \{ \metric(P_i), i = 1,\cdot,N \}
  \]
\end{itemize}
To these four strategies, we add a fifth strategy, namely we define
$\metric_{\text{worst}}$ as the one, between $\metric_{\text{min}}$
and $\metric_{\text{max}}$, that singles out the worst polygon. More
precisely, for $\metric \in \{\IC,\CR,\KE,\KA,\AP,\SE,\ER,\MP,\MA,\SR
\}$ we set
\[
\metric_{\text{worst}}(\Th) =	\metric_{\text{min}}(\Th),
\]
while for $\metric \in \{\CC,\MX,\NS \}$ we set
\[
\metric_{\text{worst}}(\Th) =	\metric_{\text{max}}(\Th).
\]



\subsection{Performance indicators}
\label{subsec:benchmark:performance}

Let us now consider how we can measure the performance of a PDE solver. Depending on the context, saying that a  solver is ``good''  can have different meanings. Typically, the first thing that comes to mind, is that a solver is good when the error between the computed and the true solution is small. However the error might be measured in different norms: for instance, while the most natural norm in our framework is the $H^1$ norm (which is spectrally equivalent to the energy norm), if one is interested in the point value of the solution, the right information on the accuracy of the method is provided by the $L^\infty$ norm of the error. The $L^2$ norm of the error is also frequently used to measure the accuracy. On the other hand, other quantities, not directly related to the error, can be of interest. For instance, if the problem considered involves inexact data,  in order to limit the effect on the  computed solution of the error on the data, the condition number of the stiffness matrix in the linear system of equations stemming from the discrete problem \eqref{eq:VEM} should be kept as small as possible, compatibly with the known fact that such a quantity increases as the mesh size decreases. The condition number of the stiffness matrix is also of interested if one aims at solving a large number of different problems,  therefore needing a good computational performance of the numerical method. 
Here we introduce the different {\em performance indicators}  that we chose to consider in our statistical analysis of the the Virtual Element Method.

\

\noindent\emph{Energy norm error}. The first indicator that we consider is the error between the computed solution $u_h$ and the continuous solution $u$, measured in the energy norm.  We recall that the quantity on the left hand side of inequality \eqref{eq:source:problem:H1:error:bound}
 is not computable. What is then usually done is to evaluate instead some computable quantity, such as the energy norm relative to the discrete equation, or the broken $H^1(\Th)$ norm of $u - \Pinabla u_h$. For our statistical analysis we chose, as first performance indicator
\[
\perf{1} = \frac{\sqrt{a_h(u_h-u_I,u_h -  u_I)}}{\sqrt{a_h(u_I,u_I)}}.
\]
Observe that this indicator depends on the data $f$ and $g$ of the model problem considered, so that for different values of such data we have different values of the indicator.

\

\noindent\emph{$L^\infty$ error}. In the finite element framework, under a shape regularity condition for the underlying triangular/quadrangular mesh, it is possible to bound, a priori, the maximum pointwise error. More precisely, under suitable assumptions on the discretization space which are satisfied by the most commonly used finite elements, if $u_h$ is the finite element solution  to Problem \eqref{eq:Poisson:strong:A}--\eqref{eq:Poisson:strong:B} defined on a quasi uniform triangular/quadrangular grid of mesh size $h$, it is possible to prove (see \cite{Brenner-Scott:2008}) that 
\[
\| u - u_h \|_{0,\infty,\Omega} \leq C \hmax^k \| u \|_{k,\infty,\Omega}.
\]While to our knowledge the problem of giving an  $L^\infty$ a priori bound on the error for the virtual element method has not yet been addressed in the literature, measuring the error in such a norm is relevant to many applications. It is then interesting to measure such an error and see how it is affected by the shape of the elements of the tessellation. However, also in this case, as we do not have access to the point values of the discrete function $u_h$, we will instead compute, as a performance indicator, the quantity
\[
\perf{2} = \frac{ \max_{i} | \dofiglob{u_h-u} |}{\max_{i} | \dofiglob{u} }.
\]
Also this indicator depends on the data $f$ and $g$ of the model problem considered.

\

\noindent\emph{$L^2$ error}. The third quantity that we consider, is the $L^2$ norm of the error. 
As usual, the proof of an  inequality  of the form \eqref{eq:source:problem:L2:error:bound} involves a duality argument and relies on the same a priori interpolation estimates used for proving \eqref{eq:source:problem:H1:error:bound}. It requires, therefore, the same shape regularity assumption on the tessellation. As for the $H^1$ and the $L^\infty$ norm errors, the $L^2$ norm of the error is not computable, and we replace it, in our experiments with the quantity
\[
\perf{3} = \frac{\| u_h - \Pi^0_k u \|_{0,\Omega}}{\| \Pi^0_k u\|_{0,\Omega}},
\]
Also this indicator depends on the data $f$ and $g$ of the model problem considered.

\

\noindent\emph{Condition number}. The condition number of the global stiffness matrix has a twofold effect:
\begin{itemize}
	\item it provides reliable information on the efficiency of iterative solvers for the linear system arising from the discretization, ad on the need (or lack thereof) of resorting to some kind of preconditioning;
	\item more importantly, it provides a bound on how errors are propagated and, possibly, amplified, in the solution process. We recall, in fact, that the condition number is defined as 
	\[
	\kappa(\stiffness) = \| \stiffness \|\,\|\stiffness^{-1} \|,
	\] 
	where $\stiffness$ is the stiffness matrix stemming from equation \eqref{eq:VEM}.
	A large condition number for $\stiffness$ might signify that the norm of $\stiffness^{-1}$ is large. In turn, this implies that possible errors on the evaluation of right hand side of \eqref{eq:VEM}, which might derive not only from round off, but also from error on the data, are amplified by the solver resulting in a possibly much larger error in the computed discrete solution. 
\end{itemize}

Estimating, a priori, the condition number of the stiffness matrix $\stiffness$ is not difficult and relies on the use of an inverse inequality of the form
\begin{equation}\label{inverse}
\| v_h \|_{1,\Omega} \leq C_{\text{inv}} \hmin^{-1} \| v_h \|_{0,\Omega}.\end{equation}

Under assumptions \textbf{G1} and \textbf{G2}, such an inequality can be proven to hold with a constant $C_{\text{inv}}$ independent of the polygon.   If \eqref{inverse} holds,  and if the chosen scaling for  degrees of freedom is such that, letting $v_h \in \Vh$ and $\vettore{v}$ the vector of its degrees of freedom, we have
\begin{equation}\label{Riesz}
\| v_h \|_{0,\Omega}^2 \simeq \vettore{v}^T \vettore v
\end{equation}
(if assumptions \textbf{G1} and \textbf{G2} hold, this is always possible \cite{ChenHuang:2018,BertoluzzaPennacchioPrada:2020}), then, for $\lambda$ eigenvalue of $\stiffness$ and $\vettore{v}$ correponding eigenvector, with $v_h$ denoting the corresponding  function in $\Vh$, we can write
\[
\vettore v^T \vettore v \lesssim \| v_h \|^2_{0,\Omega} \leq \| v_h \|_{1,\Omega}^2 \lesssim a_h(v_h,v_h) = \vettore v^T \stiffness \vettore v = \lambda \vettore v^T \vettore v,
\]
as well as
\[
\lambda \vettore v^T \vettore v =  \vettore v^T \stiffness \vettore v  = a_h(v_h,v_h) \lesssim \| v_h \|_{1,\Omega}^2 \leq  C_{\text{inv}} \hmin^{-2} \| v_h \|_{0,\Omega} \lesssim C_{\text{inv}} h^{-2} \vettore v^T \vettore v 
\]
finally yielding
\[
\vettore v^T \vettore v \lesssim \lambda \vettore v^T \vettore v  \lesssim C_{\text{inv}} \hmin^{-2} \vettore v^T \vettore v,
\]
which implies 
\[
\kappa_2 (\stiffness) \lesssim C_{\text{inv}} \hmin^{-2}.
\]
While the geometry of the elements also affects the implicit constant in the above inequality in different way (in particular through the continuity and coercivity constants and through the constants implicit in the equivalence relation \eqref{Riesz}), we underlined here the effect of the constant appearing in the inverse inequality \eqref{inverse}. Observe that if Assumption \textbf{G2} is violated, then $ C_{\text{inv}}$ is known to explode as the ratio between the diameter and the length of the smaller edge.

\

As computing the condition number exactly can be, for large matrices, extremely expensive, we used, as the performance indicator $\perf{4}$, a lower estimate of the $1$-norm condition number, computed by the Lanczos method, according to \cite{Hager:1984,Higham:2000}.

\

\noindent\emph{Condition number after preconditioning.} As already observed, for ill conditioned stiffness matrices, it is customary to resort to some form of preconditioning, and, consequently, a fairer estimation of the efficiency attained by a given solver requires taking into account the effect of preconditioning. Several approaches are available to precondition the stiffness matrix arising from the VEM discretization. We recall, among others, preconditioners based on domain decomposition, such as Schwarz \cite{Calvo:2019}, and FETI-DP and  BDDC \cite{BertoluzzaPennacchioPrada:2017, BertoluzzaPennacchioPrada:2020}, and multigrid \cite{Antonietti-Mascotto-Verani:2018}. Here we rather consider a simpler algebraic preconditioner, namely, the Incomplete Choleski factorization preconditioner \cite{GolubVanLoanBook}. Also here we use, as the preformance indicator $\perf{5}$, the lower  estimate of the $1$-norm condition number of the preconditioned stiffenss marix, computed according to \cite{Hager:1984,Higham:2000}.

\

\noindent\emph{Constant in the error estimate.} While the theoretical error bound \eqref{eq:source:problem:H1:error:bound} puts forward the dependence of the  error, measured in the $H^1(\Omega)$ norm, on the diameter of the largest element in the tessellation, the correlation analysis which, we will present in Section \ref{subsec:benchmark:results}, shows a higher correlation of the error with the average diameter size, suggesting that, in practice, an estimate such as
\begin{equation}\label{eq:error:H1:fixedtess}
\| u - u_h \|_{1,\Omega} \leq C(u) \hav^k
\end{equation}
might hold. Of course, for fixed tessellation $\Th$ and solution $u$ an equality will hold in \eqref{eq:error:H1:fixedtess}, for a suitable constant $C(u)$ depending on the tessellation $\Th$.
 We then use such a constant as a performance indicator $\perf{6}$:
\[
\perf{6} = \frac{\perf{1}}{\hav^k}.
\]

\

\noindent\emph{Aubin-Nitsche trick constant}. In the finite element framework, the Aubin-Nitsche duality trick allows to bound the $L^2$-norm of the error as
\[
\| u - u_h \|_{0,\Omega} \leq C_{\text{AN}} \hmax \| u - u_h \|_{1,\Omega},
\]
from which a bound of the form \eqref{eq:source:problem:L2:error:bound} immediately results.
While an analogous estimate cannot be proven for the VEM method, for which \eqref{eq:source:problem:L2:error:bound} is proven directly, as, once again, the two quantities $\| u - u_h \|_{0,\Omega}$ and $\hmax \| u - u_h \|_{1,\Omega}$ are both positive and finite, the above inequality can be replaced by an equality for a given constant $C_{\text{AN}}$, depending on the data and on the tessellation. Also in this case our correlation suggests to replace, in an estimate of this kind, the diameter $\hmax$ of the largest element with the average diameter, and we propose, as a performance indicator, the quantity
\[
\perf{7} = \frac{\perf{2}}{\hav \perf{1}}.
\]

\

\noindent\emph{Effectiveness of the preconditioner}. The last performance indicator aims at evaluating if and how the geometry of the tessellation affects the effectiveness of the preconditioner by comparing the condition number of the preconditioned stiffness matrix with the one of the same matrix without any preconditioning. More precisely, the last performance indicator is defined as 
\[
\perf{8} = \frac{\perf{5}}{\perf{4}}.
\]
Remark that we expect $\perf{8}$ to be less than $1$. $\perf{8}$ close to $1$ indicates an ineffective preconditioner, while $\perf{8}\ll 1$ indicates that the preconditioner is performing its role, while $\perf{8} \geq 1$ indicates the failure of the preconditioning algorithm.


\subsection{Results}
\label{subsec:benchmark:results}
\begin{figure}[htbp]
\centering
  \includegraphics[width=.9\columnwidth]{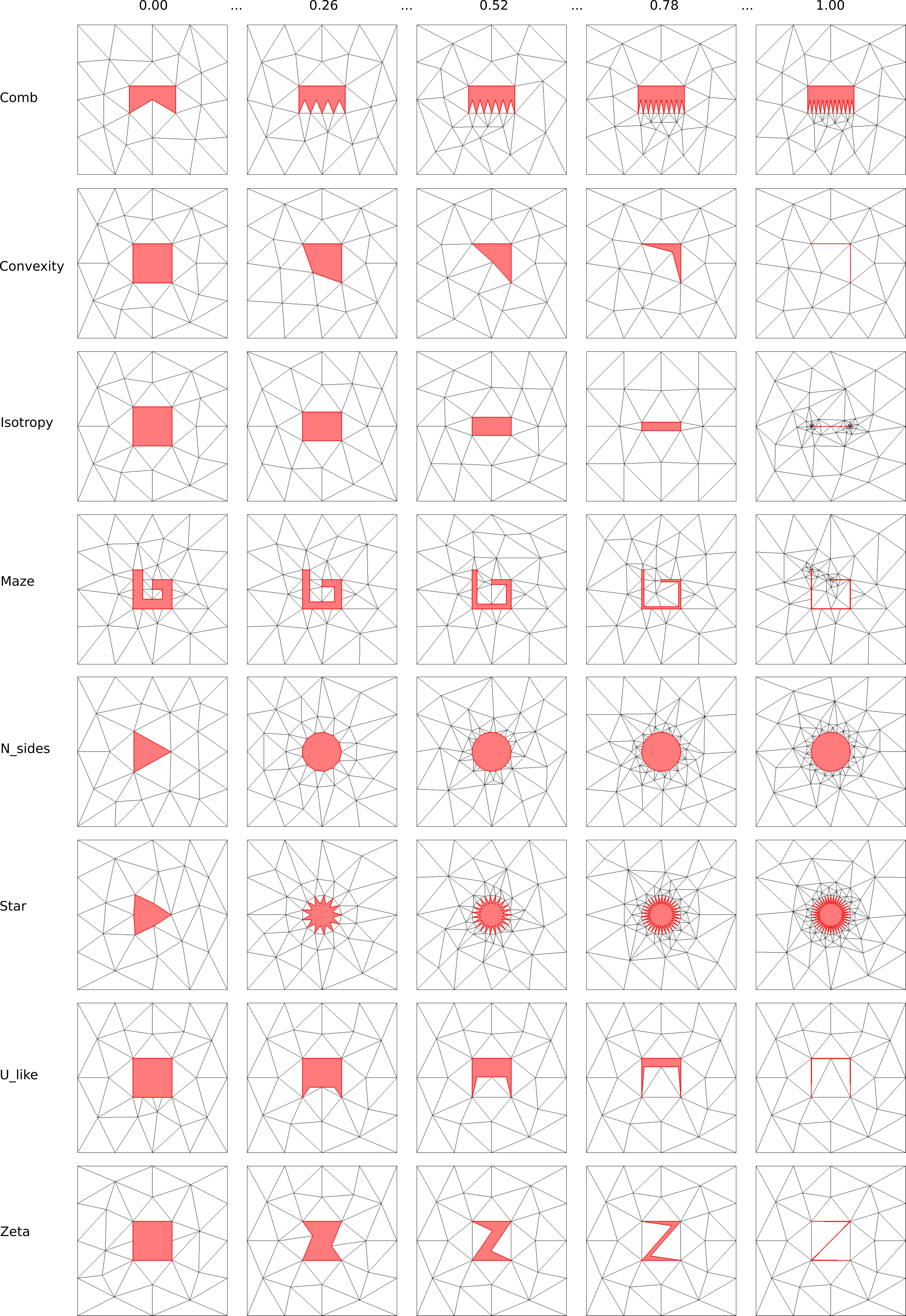}
  \caption{A subset of meshes in our dataset. On the left column, the
    list of names associated to parametric polygons (in red). On the
    top row, the values of the deformation parameter. }
  \label{fig-parametric-polygons}
\end{figure}

We considered two different test problem corresponding to two
different solutions to the Poisson equation
\eqref{eq:Poisson:strong:A}--\eqref{eq:Poisson:strong:B}, both with
$\Omega = (0,1)^2$.
\begin{description}
	\item[Test 1] For the first test problem, the groundtruth is, once again,
	\begin{equation*}
	u_1(x,y)=\frac{\sin(\pi x)\sin(\pi y)}{2\pi^2}.
	\end{equation*} 
	\item[Test 2] For the second test problem the groundtruth is the Franke function, namely
	\begin{multline*}
	u_2(x,y):=\frac 3 4 e^{-\left(
		(9x-2)^2+(9y-2)^2
		\right)/4} + \frac 3 4 e^{-
		\left(
		(9x+1)^2/49 + (9y+1)/10\right)}\\
	+ \frac 1 2 e^{-
		\left(
		(9x-7)^2 + (9y-3)^2	
		\right)/4
	}
	+\frac 1 5 e^{-
		\left(
		(9x-4)^2 + (9y-7)^2
		\right)
	}.
	\end{multline*}
\end{description}
We solved each problem with the VEM solver of order $k=1,2,3$ on 260
hybrid tessellations explicitly designed to progressively stress the
polygon quality metrics described in
Sect.~\ref{subsec:benchmark:metrics}. Specifically, a family of
parametric polygons have been designed: 
the baseline configuration of each polygon (P(0)) does not present
critical geometric features, and 20 versions of the same polygon are
generated by progressively deforming the baseline by a parameter
\textit{t}. Each polygon (its baseline version and its deformations)
is placed at the center of a canvas representing the squared domain
$[0,1]^2$. the space of the domain complementary to the polygon is
filled with triangles \cite{Shewchuk-Triangle}
\cite{livesu2019cinolib}. Figure~\ref{fig-parametric-polygons} shows
the list of our parametric polygons and how they have been used to
generate the dataset.
Note that the parametric polygons "maze" and "star" are similar to the
initial polygons of the hybrid datasets defined in
Section~\ref{subsec:violating:datasets}, meaning that they contain the
same pathologies even if edges and areas scale differently.  The
resulting meshes, however, have little in common, as meshes from
$\Dmaze$ and $\Dstar$ contain several of such polygons with different
sizes.

\begin{figure}[htbp]
  \centering
  \foreach \i in {11,...,16} { 
    \includegraphics[width=4cm]{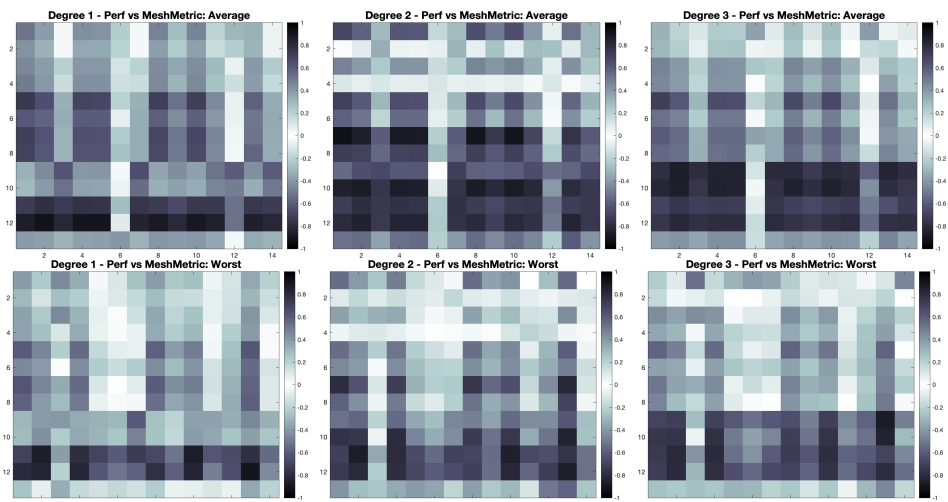}
  }
  \caption{Mesh metrics (average agglomeration method, top and worst
    agglomeration method, bottom) vs Performance metrics for $k$
    equals $1$ through $3$.}\label{fig:5av}
\end{figure}

For each test case we compute the 5 performance indexes $\perf{1}$,
$\perf{2}$, $\perf{3}$, $\perf{6}$, $\perf{7}$ corresponding to each
tessellation. Moreover, for each tessellation, we computed the three
performance indexes $\perf{4}$, $\perf{5}$ and $\perf{8}$, for a grand
total of 13 performance indexes for each tessellation. We then
computed the Spearman correlation index
\cite{Spearman:1904,MeyersWellBook} between the 14 quality metrics,
with the five aggregation methods, and the 13 performance indexes. The
results are displayed in figure \ref{fig:5av}.  As the Euclidean norm
aggregation method gives results that are substantially similar to the
average aggregation method and as the worst aggregation method
summarized the most significant of the minimum and maximum aggregation
methods, we only show the results for average and worst aggregation
method.

\begin{figure}[htbp]
  \centering
  \foreach \i in {26,...,34} { 
    \includegraphics[width=4cm]{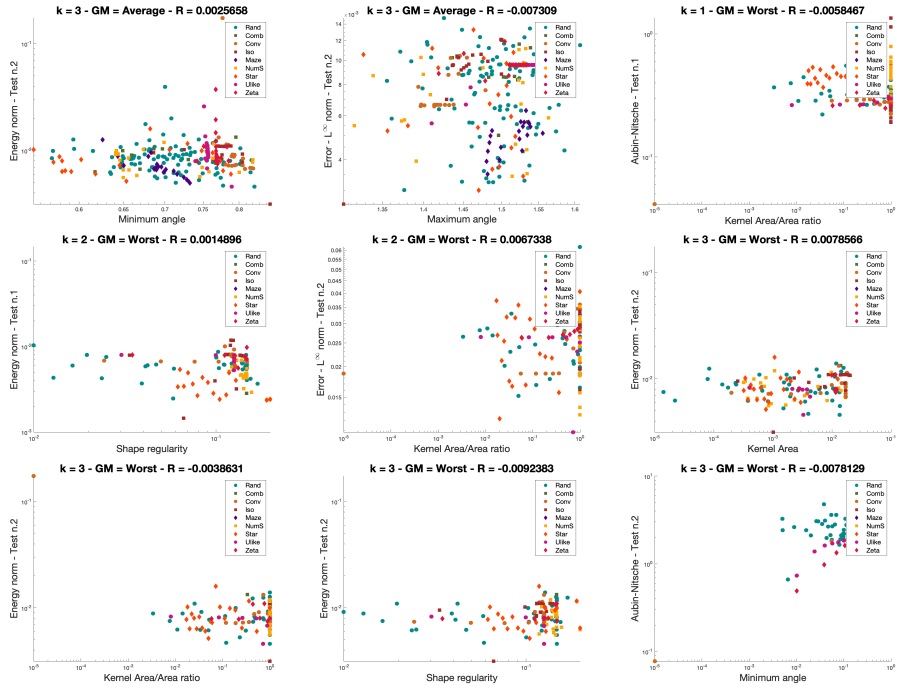}
  }
  \caption{Scatterplots corresponding to an absolute value of the
    Spearman correlation lower than .03. On the x-axis the geometric
    metric, on the $y$-axis the performance metric.
  }\label{fig:scatterlow}
\end{figure}

\begin{figure}[htbp]
  \centering
  \foreach \i in {17,...,25} {
    \includegraphics[width=4cm]{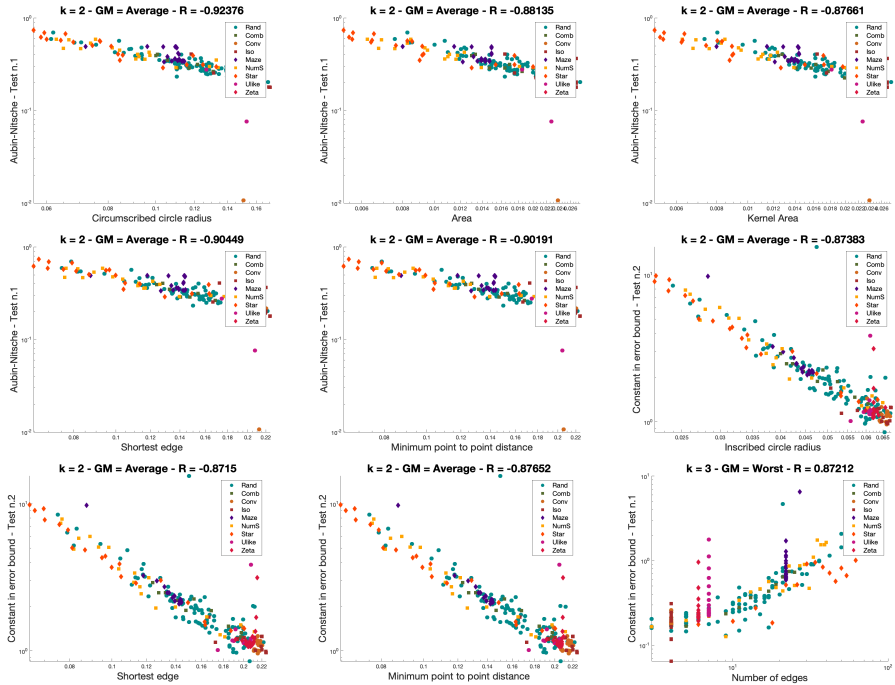}  
  }
  \caption{Scatterplots corresponding to an absolute value of the
    Spearman correlation greater than .9. On the x-axis the geometric
    metric, on the $y$-axis the performance
    metric.}\label{fig:scatterhigh}
\end{figure}

We next consider the cases that result in the overall lowest and
highest correlation index. For such cases, we present the related
scatterplots (displayed either in semi-logarithmic or in logarithmic
scale). More precisely the scatterplot corresponding to low
correlation are shown in Figure \ref{fig:scatterlow}, while in Figure
\ref{fig:scatterhigh} we can look at the scatterplots corresponding to
the high correlation cases. The first thing that jumps to the eye is
that all high correlation case except one correspond to the average
aggregation method, while the converse holds for the low correlation
(all cases but two correspond to the worst aggregation
method). Looking at Figure \ref{fig:scatterlow} we also see that the
SR metric has low correlation with the error in the energy norm.  This
confirms the experimental observations, already pointed out in the
literature, and is particularly interesting, as it shows that
violating Assumption~\textbf{G1}, which is the commonly used
assumption in the analysis of VEM methods, does not necessarily result
in a deterioration of the accuracy of the method. This suggests to try
to carry out the convergence analysis in the absence of such a bound
(and, in particular, for meshes of non star shaped polygons).  The
minimum and maximum angle also seem to have little correlation with
the error.  Conversely, two quantities that are clearly highly
correlated with several of the metrics are the constants in the error
estimate $\perf{6}^{(1)}$ and $\perf{6}^{(2)}$ and the Aubin Nitsche
trick constants $\perf{7}^{(1)}$ and $\perf{7}^{(2)}$.  We also
observe that the constant in the error bound is highly correlated with
the highest number of edges. This is also very interesting, as it
suggests that even a single polygon with a high number of edges can
affect the performance of the method. It is then worth devoting some
effort to the design of variants of the VEM, aimed at attaining
robustness with respect to the number of edges.

  


\begin{figure}[htbp]
  \centering
  \begin{tabular}{c c}
  \includegraphics[width=.47\textwidth]{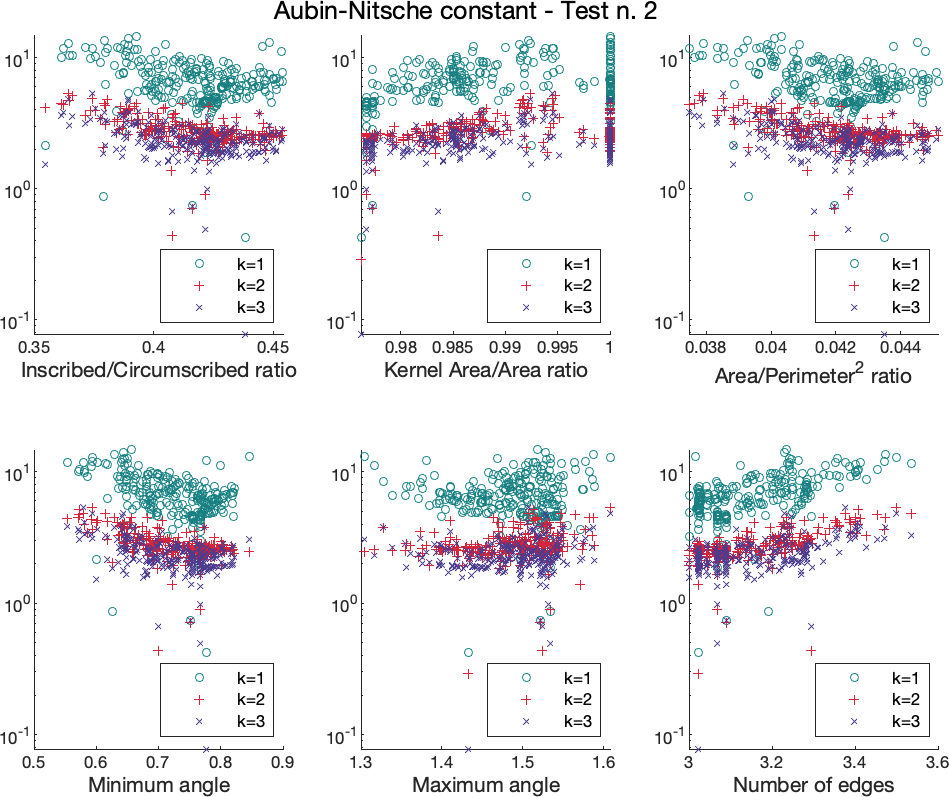} &
  \includegraphics[width=.47\textwidth]{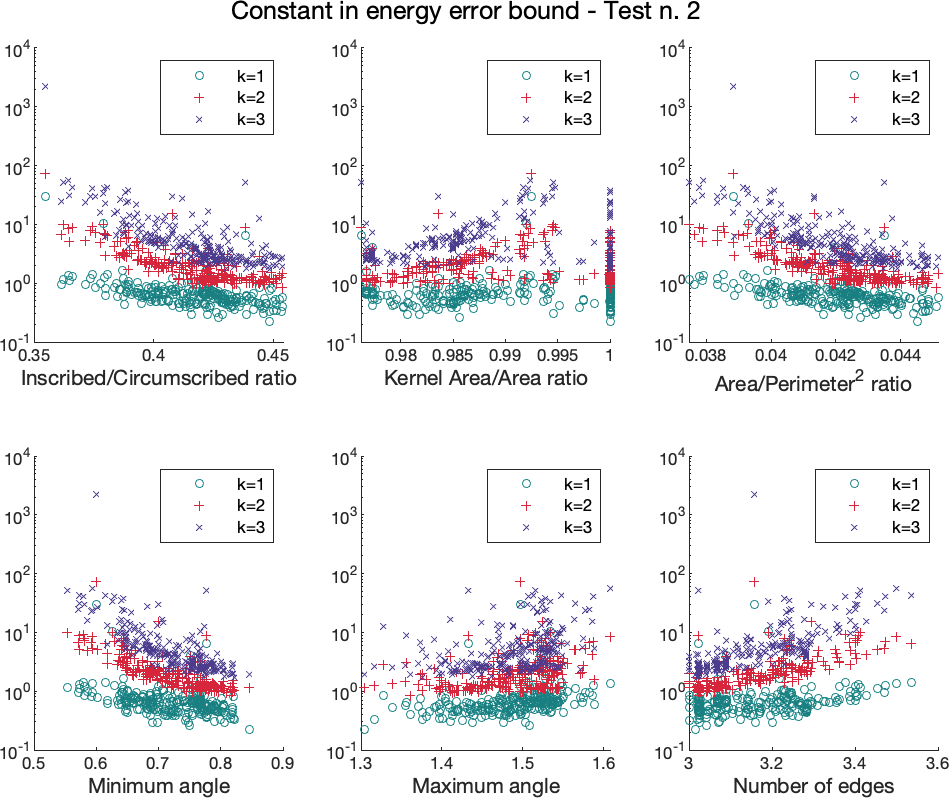}\\
  (a) & (b)
  \end{tabular}
  \caption{Relation between the Aubin-Nitsche constant (a) or the constant 
  in the classical a priori error bound (b) and the scale
  invariant metrics for $k=1$, $2$ and $3$.} \label{fig:combined1}
\end{figure}

\begin{figure}[htbp]
  \centering
  \includegraphics[width=.47\linewidth]{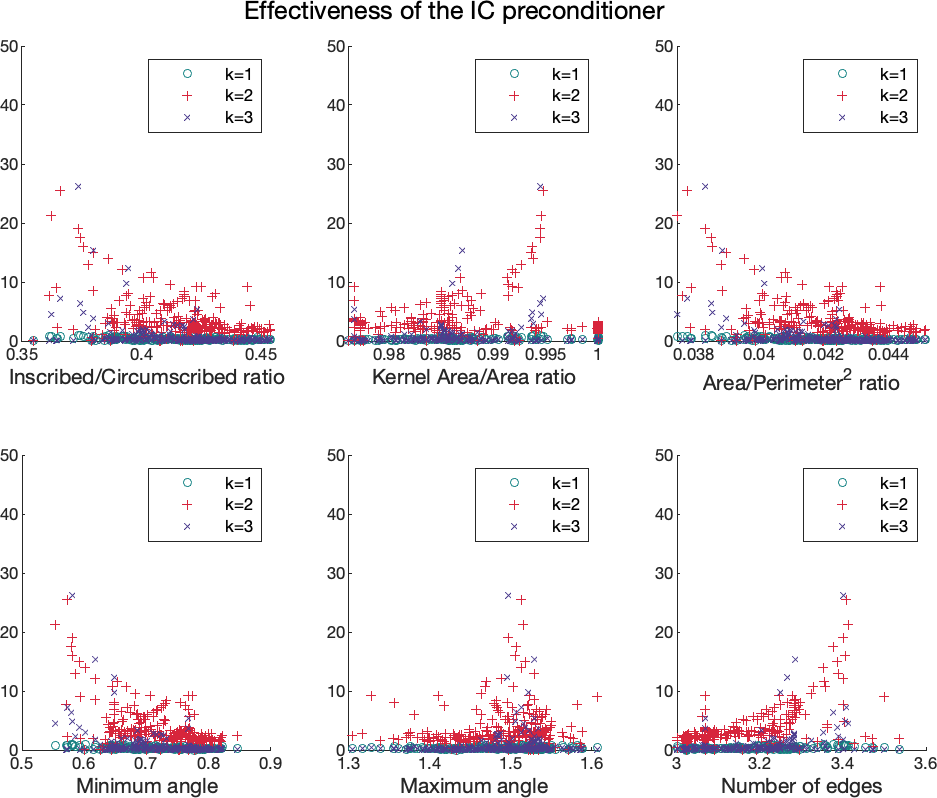}
  \caption{Relation between the effectiveness of the incomplete
    Choleski preconditioner and the scale invariant metrics for $k=1$,
    $2$ and $3$.}\label{fig:combined3}
\end{figure}

\PGRAPH{Comparison between methods of different orders} The previously
presented set of tests is also aimed at gaining some insight on the
combined influence of geometry and order of the method on the
performance of the method itself.  Theoretical results on the
performance of the VEM rarely track the dependence of the constant in
the inequality on the order $k$ of the method. Notable exceptions are
\cite{BeiraodaVeiga-Chernov-Mascotto-Russo:2016,BeiraodaVeiga-Chernov-Mascotto-Russo:2018},
that however, in order to give an a priori error bound bound, require
the polygonal meshes to satisfy the very restrictive assumptions
\textbf{G1} and \textbf{G2}. Focusing on three of the performance
indexes, which we expect to be scale independent, and on the six scale
independent geometric quality measures, we present the superposed
scatter plots for $k=1$, $2$ and $3$. The resulting comparison yields
some interesting observations.  In Figure \ref{fig:combined1}(a), we see
that the Aubin-Nitsche constant for $k=2$ and $3$ is lower and less
spread out than for $k=1$.  The opposite happens for the constant in
the classical error bound of the form
\eqref{eq:source:problem:H1:error:bound}. Such constant increases as
$k$ increases (as it is to be expected), but it also appear that
higher values of $k$ amplify the negative influence of the bad quality
of the mesh. Figure \ref{fig:combined3} is particularly interesting:
for $k=2$, the incomplete Choleski preconditioner is less effective
than it is for $k = 1$ and $3$, and it even appears to be failing in a
number of cases, giving an effectiveness index $> 1$. This might be
the consequence of a difference between odd and even order method,
that is sometimes encountered in non conforming discretizations.

\section{Mesh quality indicators}
\label{sec:indicators}

In this section we define a \textit{mesh quality indicator}, that is,
a scalar function capable of providing insights on the behaviour of
the VEM over a particular mesh $\Th$, before actually computing the
approximated solution $\ush$.
This indicator depends uniquely on the geometry of the mesh elements,
but has interesting correspondences with the performance of the VEM,
in terms of approximation error and convergence rate.


\subsection{Definition}
\label{subsec:indicators:definition}
We start from the geometrical assumptions defined in
Section~\ref{subsec:star:assumptions}.
The driving idea is that instead of imposing an absolute condition
that a mesh can only satisfy or violate, we want to measure
\textit{how much} the mesh satisfies or violate that condition.
This approach is more accurate, as it captures small quality
differences between meshes, and it does not exclude a priori all the
particular cases of meshes only slightly outside the geometrical
assumptions.\\
From each geometrical assumption \textbf{Gi}, \textbf{i
}$=1,\ldots,4$, we derive a scalar function $\varrho_i: \{\P \subset
\Omega_h\}\to [0,1]$ defined element-wise, which measures how well the
element $\P\in \Omega_h$ meets the requirements of \textbf{Gi} from 0
($\P$ does not respect \textbf{Gi}) to 1 ($\P$ fully respects
\textbf{Gi}).

\smallskip
From assumption \textbf{G1} we derive the indicator $\varrho_1$, which
can be interpreted as an estimate of the value of the constant $\rho$
from \textbf{G1} on the polygon $\P$:
\[
  \varrho_1(\P) = \frac{k(\P)}{|\P|},
\]
being $k(\P)$ the area of the \textit{kernel} of a polygon $\P$,
defined as the set of points in $\P$ from which the whole polygon is
visible (cf. \textit{(KE)} of Section~\ref{subsec:benchmark:metrics}).
Therefore, $\varrho_1(\P)=1$ if $\P$ is convex, $\varrho_1(\P)=0$ if
$\P$ is not star-shaped and $\varrho_1(\P)\in(0,1)$ if $\P$ is concave
but star-shaped.

\smallskip
Similarly, the function $\varrho_2$ returns an estimate of the
constant $\rho$ introduced in \textbf{G2}, expressed trough the ratio
$|\E|/h_\P$:
\[
  \varrho_2(\P) = \frac{\min(\sqrt{|\P|}, \ \min_{e \in \bP}|e|)}
  {\max(\sqrt{|\P|}, \ h_\P)}.
\]
The insertion of the quantity $\sqrt{|\P|}$ is needed in order to
scale the indicator in the range $[0,1]$ and to avoid pathological
situations.

\smallskip
Function $\varrho_3$ is a simple counter of the number of edges of the
polygon, which penalizes elements with numerous edges as required by
\textbf{G3}:
\[
  \varrho_3(\P) = \frac{3}{\#\left\{\E\in\bP\right\}}.
\]
It returns 1 if $\P$ is a triangle, and it decreases as the number of edges increases.

\begin{figure}[htbp]
  \centering
  \includegraphics[width=.3\linewidth]{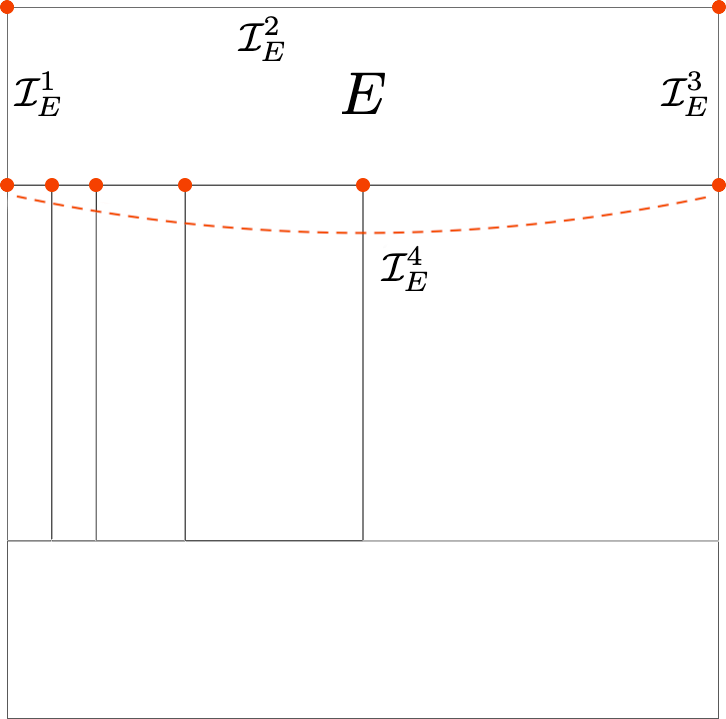}
  \caption{One-dimensional mesh relative to the top bar element $\P$ of a $\Djenga$ base mesh.  The boundary of $\P$ is represented by a mesh $\calI_\P=\{\calI_\P^1, \calI_\P^2, \calI_\P^3, \calI_\P^4\}$, where the sub-meshes $\calI_\P^1, \calI_\P^2$ and $\calI_\P^3$ contain,
respectively, the left, top and right edge of $\P$, while $\calI_\P^4$
contains all the aligned edges in the bottom of $\P$.}
  \label{fig:rho4}
\end{figure}

\smallskip
Last, we recall from Section~\ref{subsec:star:assumptions} that the
boundary of a polygon $\P$ can be considered as a 1-dimensional mesh
$\calI_\P$, which can be subdivided into disjoint sub-meshes
$\calI_\P^1, \ldots, \calI_\P^N$, each one containing possibly more
than one edge of $\P$.
In practice, we consider as a sub-mesh the collection of all edges
whose vertices lie on the same line: an example is shown in Fig.~\ref{fig:rho4}.

The indicator $\varrho_4$ returns the minimum ratio between the
smallest and the largest element in every $\calI_\P$, that is a
measure of the quasi-uniformity of $\calI_\P$ imposed by \textbf{G4}:
\[
\varrho_4(\P) = \min_i\frac{\min_{\E\in\calI_\P^i}|\E|}
{\max_{\E\in\calI_\P^i}|\E|}.
\]

\medskip
Combining together $\varrho_1, \varrho_2$, $\varrho_3$ and
$\varrho_4$, we define a global function $\varrho:\{\Th\}_h\to [0,1]$
which measures the overall quality of a mesh $\Omega_h$.
Given a dataset $\calD$, we can study the behaviour of $\varrho(\Th)$
for $\hh\to0$ and determine the quality of the meshes through the
refinement process.
In particular, we combine the indicators with the formula
$\varrho_1\varrho_2 + \varrho_1\varrho_3 + \varrho_1\varrho_4$ as it
reflects the way in which these assumptions are typically imposed:
\textbf{G1} and \textbf{G2}, \textbf{G1} and \textbf{G3} or
\textbf{G1} and \textbf{G4} (but not, for instance, \textbf{G2} and
\textbf{G3} simultaneously):
\begin{equation}  \label{eq:indicator}
  \varrho(\Omega_h) = \sqrt{\frac{1}{\#\left\{ \P \in \Omega_h \right\}} \ \sum_{\P\in\Omega} \frac{\varrho_1(\P)\varrho_2(\P) + \varrho_1(\P)\varrho_3(\P)+ \varrho_1(\P)\varrho_4(\P)}{3}}.
\end{equation}
We list some observations on the newly defined mesh quality indicator
$\varrho$:
\begin{itemize}
\item we have $\varrho(\Omega_h)=1$ if and only if $\Omega_h$ contains
  only equilateral triangles, $\varrho(\Omega_h)=0$ if and only if
  $\Omega_h$ contains only non star-shaped polygons, and
  $0<\varrho(\Omega_h)<1$ otherwise;
\item the indicator $\varrho$ only depends on the geometrical
  properties of the mesh elements (because the same holds for
  $\varrho_1, \varrho_2$, $\varrho_3$ and $\varrho_4$), therefore it
  can be computed before applying the numerical scheme;
\item the construction of $\varrho$ is easily upgradeable to future
  developments: whenever new assumptions on the features of a mesh
  should come up, one simply needs to introduce a new function
  $\varrho_i$ that measures the violation of the new assumption and
  opportunely insert it into equation~(\ref{eq:indicator});
\item similarly, this indicator is easily extendable to other
  numerical schemes by substituting the assumptions designed for the
  VEM with the assumptions on the new scheme, and defining the
  relative indicators $\varrho_i$.
\end{itemize}


\subsection{Results}
\label{subsec:indicators:results}
\begin{figure}[htbp]
  \centering
  \begin{tabular}{cc}
    \includegraphics[width=.45\linewidth]{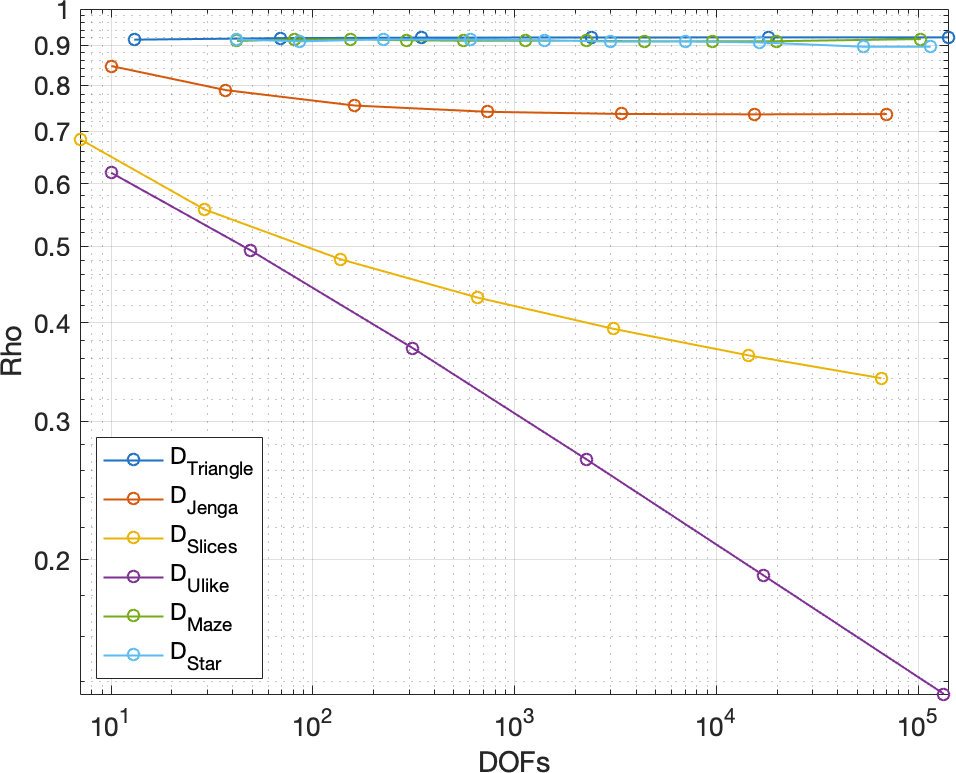}
    & \qquad
    \includegraphics[width=.45\linewidth]{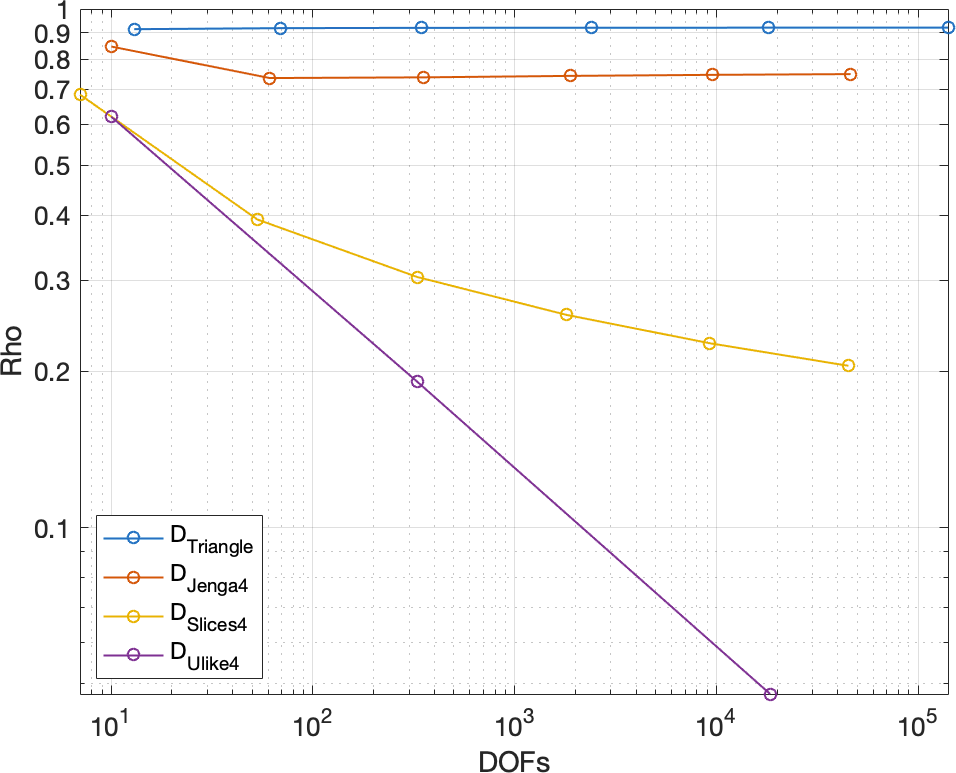}
    \\
    (a) & (b)\\
  \end{tabular}
  \caption{Indicator $\varrho$ for all datasets from
    Section~\ref{subsec:violating:datasets}.}
\label{fig:indicator}
\end{figure}

We evaluated the indicator $\varrho$ over the datasets defined in
Section~\ref{subsec:violating:datasets}; results are shown in
Fig.~\ref{fig:indicator}.
We want to investigate if the $\varrho$ values of the meshes of a
certain dataset in Fig.~\ref{fig:indicator}, computed before solving
the problem, are somehow related to the performance of the VEM over
the meshes of the same dataset, shown in Fig.~\ref{fig:performance1}
and \ref{fig:performance4}, in terms of error approximation and
convergence rate.

Since $\varrho$ does not depend on the polynomial degree $k$ nor on
the type of norm used, we have to consider an hypothetical average of
the plots for the different $k$ values and for the different norms
($L^2$ and $H^1$) from Fig.~\ref{fig:performance1} and
\ref{fig:performance4}.\\
As already seen, for an ideal dataset made by meshes containing only
equilateral triangles, $\varrho$ would be constantly equal to 1.
We assume this value as a reference for the other datasets: the closer
$\varrho$ is on a dataset to the line $y=1$, the smaller is the
approximation error that we expect that dataset to produce.
Similarly, the more negative is the $\varrho$ slope, the worse is the
convergence rate that we expect the method to obtain over that
dataset.

\smallskip
For meshes belonging to $\Dtriangle$, $\varrho$ is almost constant and
very close to 1, thus foreseeing the excellent convergence rates and
the low errors relative to $\Dtriangle$ in every sub-figure of
Fig.~\ref{fig:performance1}.

\smallskip
The plots for $\Dmaze$ and $\Dstar$ in Fig.~\ref{fig:indicator}(a) are
close to the $\Dtriangle$ one, hence we expect the method to perform
similarly over the three of them.
This is confirmed by Fig.~\ref{fig:performance1}: $\Dmaze$ and
$\Dstar$ are almost coincident and very close to $\Dtriangle$ until
the very last meshes, especially in the $L^2$ plots.

\smallskip
The $\Djenga$ plot in Fig.~\ref{fig:indicator}(a) suggests a perfect
convergence rate but greater error values with respect to the previous
three, and again this behaviour is respected in
Fig.~\ref{fig:performance1}.\\
The curve relative to $\Dslices$ in Fig.~\ref{fig:indicator}(a) is
quite distant from the ideal value of 1.
Moreover, the curve keeps decreasing from mesh to mesh, even if we may
assume that it would flatten within few more meshes.
Looking at Fig.~\ref{fig:performance1}, we notice that $\Dslices$
produces an error significantly higher than the previous ones
($\Dtriangle, \Dmaze, \Dstar, \Djenga$), and in some cases the $H^1$
error convergence rate is significantly lower than the theoretical
estimate, thus confirming the prediction.\\
Last, the $\varrho$ values in Fig.~\ref{fig:indicator}(a) predict huge
errors and a completely wrong convergence rate for $\Dulike$.
In fact, this dataset is the one with the worst performance in
Fig.~\ref{fig:performance1}, where the VEM does not even always
converge (for example in the case $k=1$ with the $H^1$ seminorm).

\smallskip
As far as multiply refined datasets are concerned, we note that, since
$\varrho$ only depends on the geometry of the elements, it is not
affected by numerical errors.\\
The $\varrho$ plot for $\DjengaM$ in Fig.~\ref{fig:indicator}(b) is
very similar to the one relative to $\Djenga$ in
\ref{fig:indicator}(a), therefore we should expect $\DjengaM$ in
Fig.~\ref{fig:performance4} to perform similarly to $\Djenga$ in
Fig.~\ref{fig:performance1}.
This is actually the case, at least until the finest mesh for $k=3$
where numerical problems appear which $\varrho$ is not able to
predict.\\
Also $\DslicesM$ and $\Dslices$ are similar in
Fig.~\ref{fig:indicator}, but $\varrho$ decreases faster on the first
than on the second one, reaching a $\varrho$ value of $\sim0.2$
instead of $\sim0.34$ within a smaller number of meshes.
Again, this information is correct because the method performs
similarly over $\DslicesM$ and $\Dslices$ until condition numbers
explode, in the last two meshes for every value of $k$.\\
Last, the $\varrho$ plot of $\DulikeM$ is significantly worse than any
other (including the one of $\Dulike$), both in terms of distance from
$y=1$ and slope.
In Fig.~\ref{fig:performance4} we can observe how, even if $\DulikeM$
does not properly explode (as it suffers less from numerical problems,
cf. Table~\ref{table:numerical_performance}), the approximation error
and the convergence rate are the worst among all the considered
datasets.

\medskip
Summing up these results, we conclude that the indicator $\varrho$ is
able, up to a certain accuracy, to predict the behaviour of the VEM
over the considered datasets, both in terms of error magnitude and
convergence rate.
Looking at the $\varrho$ values it is possible to estimate in advance
if a mesh or a dataset is going to be more or less critical for the
VEM, and it is possible to compare the values relatives to different
meshes and datasets to understand which one is going to perform
better.
The prediction may be less accurate in presence of very similar
performance (e.g. the case of $\Dmaze$ and $\Dstar$) or in extremely
pathological situations, where numerical problems become so
significant to overcome any influence that the geometrical features of
the mesh may have on the performance (e.g. the last meshes of
$\DjengaM$ and $\DslicesM$).
\section{\texttt{PEMesh} benchmarking tool\label{sec:PEMesh-GUI}}
\label{sec:pemesh}


Main existing tools for the numerical solution of PDEs include
\emph{VEMLab}~\cite{vemlab}, which is an open source MATLAB library
for the virtual element method and \emph{Veamy}, which is a free and
open source C++ library that implements the virtual element method
(C++ version of~\cite{vemlab}).
This library allows the user to solve 2D linear elasto-static problems
and the 2D Poisson problem~\cite{veamy}.
Additional tools include the 50-lines MATLAB implementation of the
lowest order virtual element method for the two-dimensional Poisson
problem on general polygonal meshes~\cite{vem-50-lines} and the MATLAB
implementation of the lowest order Virtual Element Method
(VEM)~\cite{Mascotto:2018}.
\begin{figure}[htbp]
  \centering
  \begin{tabular}{c}
    \includegraphics[width=.9\textwidth]{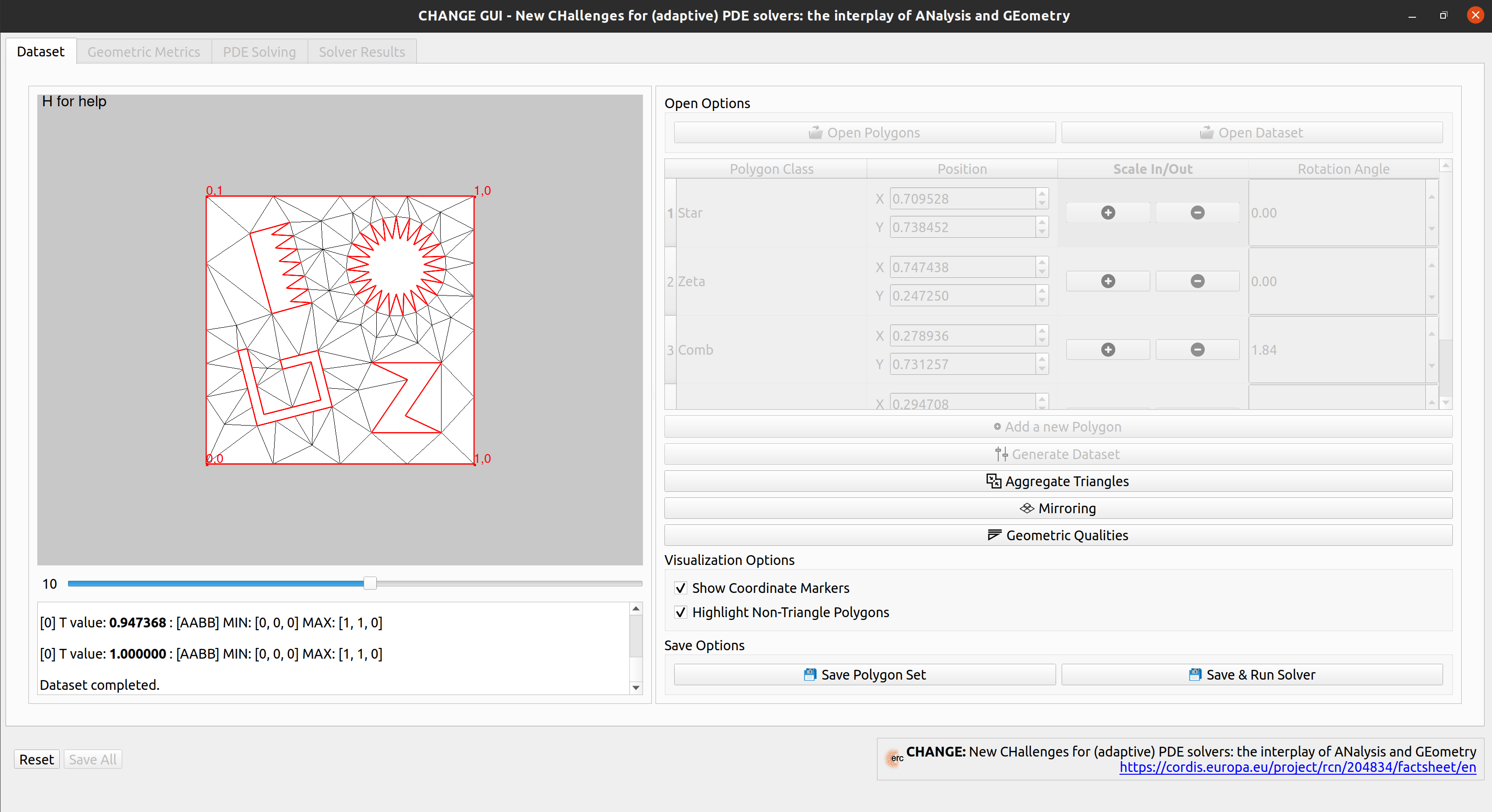}
    \\(a)\\
    \includegraphics[width=.9\textwidth]{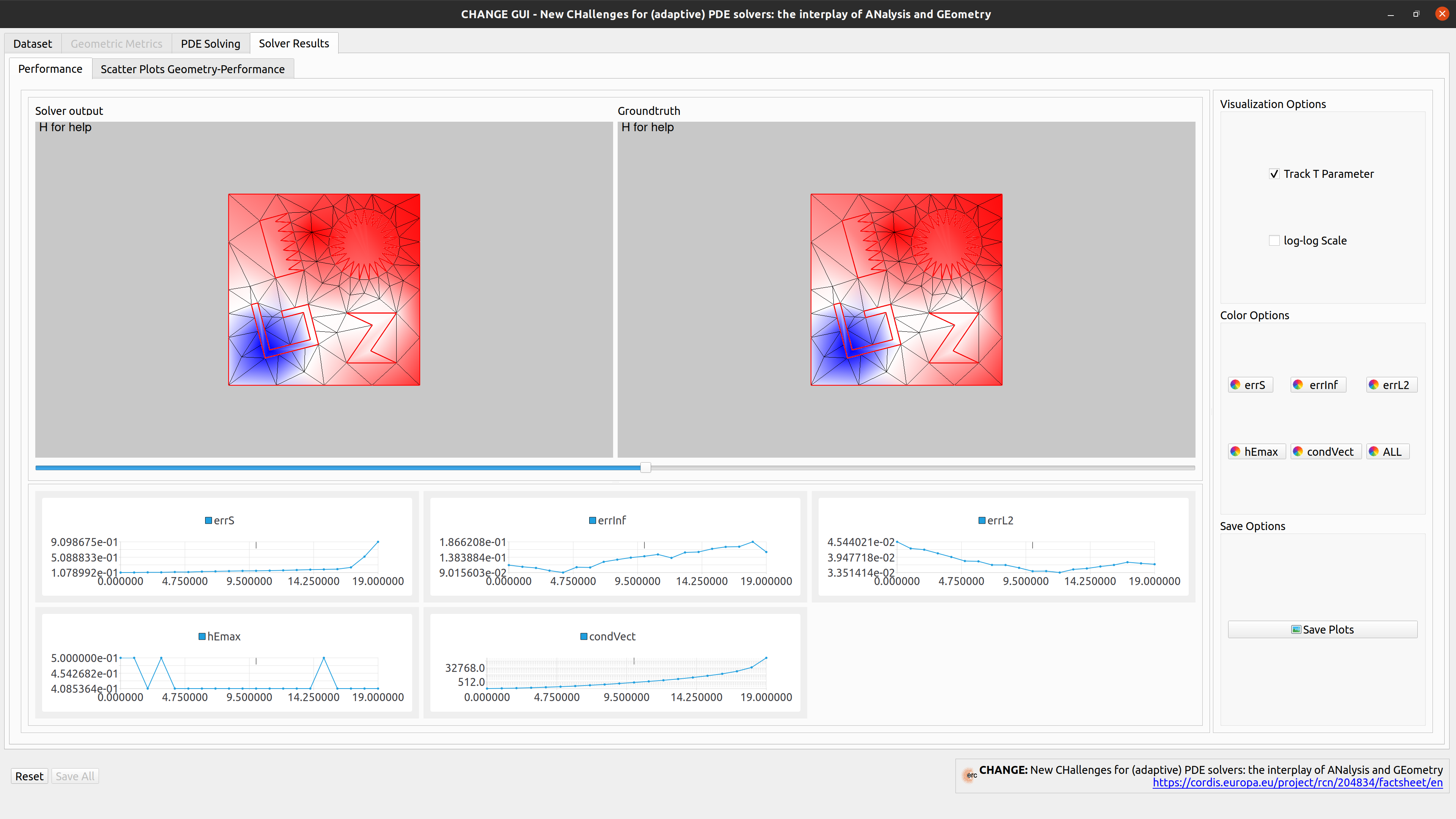}
    \\(b)
  \end{tabular}
  \caption{(a) Main window of \texttt{PEMesh} and (b) examples of PEM
    solver results.
    On the top, both solver output and ground truth are color-mapped
    on the input polygon meshes, while on the bottom a set of linear
    plots show how solver performances vary in the data
    set.\label{fig:app}}
\end{figure}

In this context, \texttt{PEMesh}~\cite{PAPER-CHANGE-GUI-2020}
(Fig.~\ref{fig:app}) is an open-source software tool for experiments
on the analysis and design of polytopal meshes for PEM
solvers.
\texttt{PEMesh} is an advanced graphical tool that seamless integrates
the geometric design of 2D domains and PEM simulations.
It supports the design and generation of complex input polygonal
meshes by stressing geometric properties and provides the possibility
to solve PEMs on the generated meshes.
Furthermore, \texttt{PEMesh} allows us to correlate several geometric
properties of the input polytopal mesh with the performances of PEM
solvers, and to visualise the results though customisable and
interactive plots.
\texttt{PEMesh} also supports to evaluate the dependence of the
performances of a PEM solver on geometrical properties of the input
polygonal mesh. \texttt{PEMesh} includes four main modules, which are
briefly described in the following.

\smallskip
\emph{Polygon mesh generation~$\&$ loading} allow us to load one or
more existing meshes or to generate new ones from scratch by either
exploiting a set of polytopal elements available in the module or
providing an external polytopal element.
In this way, the user can easily generate a large set of options and
parameters.\\
We used this module for the generation of most of the datasets
considered in this work: $\Dtriangle, \Dstar, \Dmaze$ and all datasets
from Section~\ref{subsec:benchmark:results}.
The mirroring and multiple mirroring datasets defined in
Section~\ref{subsec:violating:datasets} have been generated separately
and then loaded into the software to be used in the subsequent
modules.

\begin{figure}[htbp]
  \centering
  \begin{tabular}{c c}
    \includegraphics[width=.47\textwidth, height=.25\textheight]{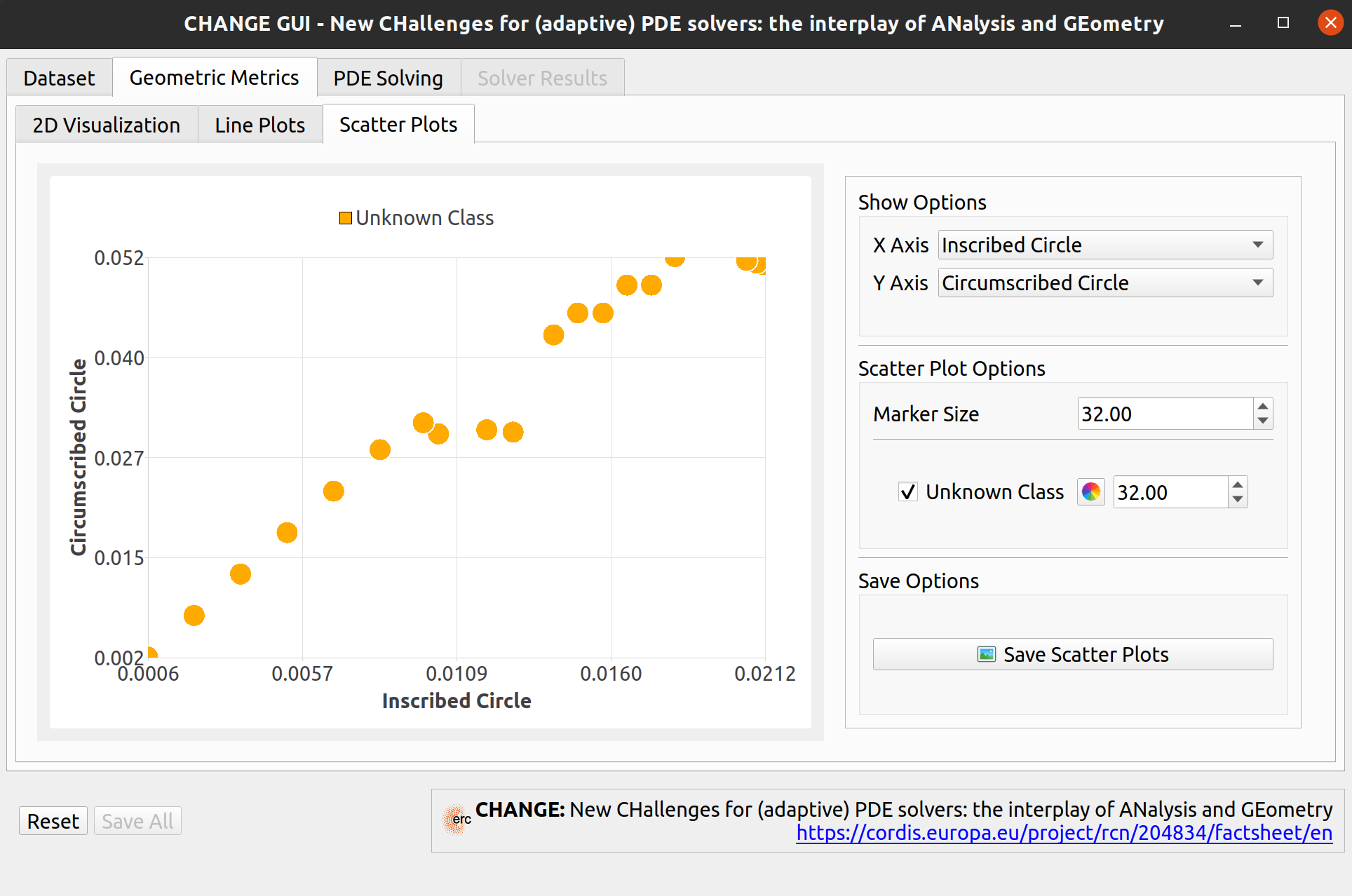} &
    \includegraphics[width=.47\textwidth, height=.25\textheight]{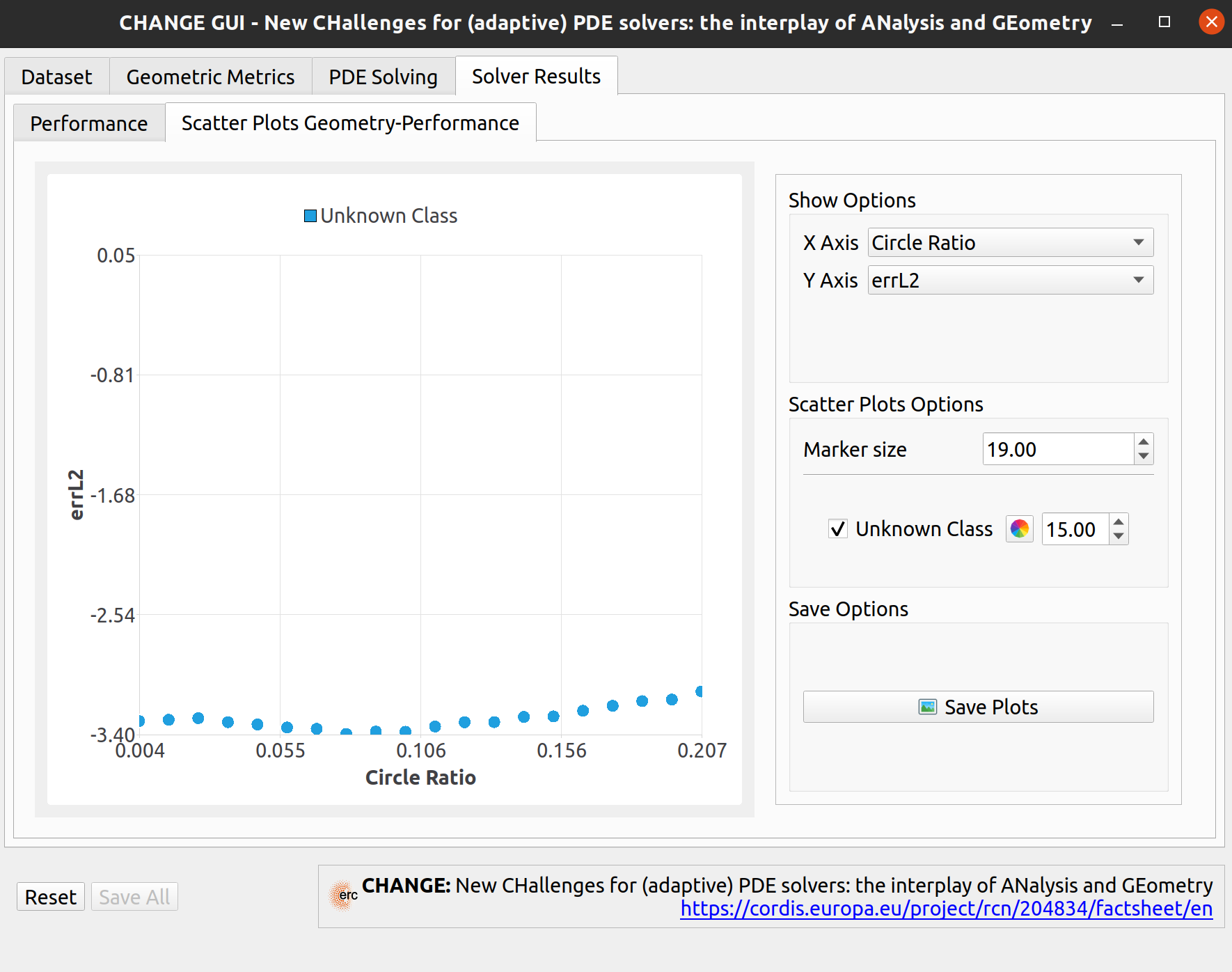}\\
    (a) & (b)\\
  \end{tabular}
  \caption{PEMesh window specialized on the visualization of the
    correlation between geometric metrics (a) and the correlation of 
    geometric metrics and PEM performances (b).
    The right side of the window provides some visualization options
    to enable the possibility to select which metrics to be visualized
    and how (i.e. size, color, ...). 
    \label{fig:geometry-pem-scatter}}
\end{figure}

\smallskip
\emph{Geometric analysis} allows the user to perform a deep analysis
of geometric properties of the input polygonal meshes, to correlate
each of them with the others, and to visualise the results of this
analysis though scatter plots (see
Fig.~\ref{fig:geometry-pem-scatter}(a)).

\smallskip
\emph{PEM solver} allows the user to run a PEM solver and to analyse
its performances on input polygonal meshes.
Any PEM solver, which is run from command line, can be integrated in
the software together with its output.
The PEM and ground-through (if any) solutions of the PEM are shown
directly on the meshes, while the performances of the solver are
visualized through linear plots.\\
This module has been used for solving the problem and analyzing the
performance of the VEM over all datasets in this work.

\smallskip
\emph{Correlation visualization} allows the user to analyse the
correlation between geometric properties of the polygonal meshes and
numerical performances of the PEM solver. Then, the results are
visualised as customizable scatter plots, as shown in
Fig.~\ref{fig:geometry-pem-scatter}(b).

\medskip
\texttt{PEMesh} provides the possibility to solve PEMs on a polygonal
mesh and to visualize the performances of any PEM solver.
To this end, the PEM solver is not part of the tool, but is handled as
an external resource that is called by \texttt{PEMesh} graphical
interface.
Through a general approach, \texttt{PEMesh} assumes that the selected
PEM solver takes an input mesh and returns both the solution and the
ground-truth (if any) of a PDE, together with statistics (e.g.,
approximation accuracy, condition number of the mass and stiffness
matrices) of the PEM solver.
Then, the results of the PEM solver are visualised in a specialized
window, and the numerical and geometric metrics are shown as linear
plots (Fig.~\ref{fig:app}b).
Finally, a set of visualization options are available for the
customisation of color-maps and linear plots.
Through a modular and customisable system, \texttt{PEMesh} allows us
to run any PEM solver, thus providing an easy way to generate complex
discrete polytopal meshes and to study the correlation between
geometric properties of the input mesh and numerical PEMs solvers.

\section*{Acknowledgement}
The Authors acknowledge the financial support of the ERC Project
CHANGE, which has received funding from the European Research Council
under the European Union’s Horizon 2020 research and innovation
program (grant agreement no.~694515).


\end{document}